\documentclass[a4paper, 12pt]{article}
\usepackage{enumerate, theorem}
\usepackage{amsmath, amsfonts, amssymb}
\usepackage[height=22.5cm, width=15cm]{geometry}
\usepackage[all]{xy}
\usepackage{mathrsfs, yfonts}

\DeclareMathOperator{\id}{id}

\DeclareMathOperator{\C}{\mathbb{C}}

\newcommand{\A}{\tilde{\mathcal{A}}}

\newcommand{\parag}[1]{\paragraph{\sc{#1.}}}

\newtheorem{thm}{Th\'eor\`eme}[subsection]
\newtheorem{defn}[thm]{D\'efinition}
\newtheorem{cor}[thm]{Corollaire}
\newtheorem{prop}[thm]{Proposition}
\newtheorem{lemma}[thm]{Lemme}

\begin{document}

\title{ Contruction of holomorphic parameters invariant by  change of variable in the Gauss-Manin connection of an holomorphic map to a disc.}

\author{Daniel Barlet\footnote{Barlet Daniel, Institut Elie Cartan UMR 7502  \newline
Nancy-Universit\'e, CNRS, INRIA  et  Institut Universitaire de France, \newline
BP 239 - F - 54506 Vandoeuvre-l\`es-Nancy Cedex.France. \newline
e-mail : barlet@iecn.u-nancy.fr}.}

\date{14/01/11}

\maketitle

\section*{Abstract.}

When we consider a proper holomorphic map \ $\tilde{f }: X \to C$ \ of a complex manifold \ $X$ \ on a smooth complex curve \ $C$ \  with a critical value at a point \ $0$ \ in \ $C$, the choice of a local coordinate near this point allows to dispose of an holomorphic function \ $f$. Then we may construct, using this function,  an (a,b)-modules structure on the cohomology sheaves of the formal completion (in \ $f$) \ of the complex of sheaves \ $(Ker\, df^{\bullet},d^{\bullet})$. These (a,b)-modules represent a filtered version of the Gauss-Manin connection of \ $f$. The most simple example of this construction is the Brieskorn module (see [Br.70]) of a function with an isolated singular point. See [B.08] for the case of a 1-dimensional critical locus. \\
But it is clear that this construction depends seriously on the choice of the function \ $f$ \ that is to say on the choice of the local coordinate near the critical point \ $0$ \ in the complex curve \ $C$.\\
The aim of the present paper is to study the behaviour of such constructions when we make a change of local coordinate near the origin. We consider the case of \ $[\lambda]-$primitive frescos, which are monogenic geometric (a,b)-modules corresponding  to a minimal filtered differential equation associated to a relative de Rham cohomology class on \ $X$ \ (see [B.09-a] and [B.09-b]).  \\
 An holomorphic parameter is a function on the set of isomorphism classes of frescos which behave holomorphically in an holomorphic family of frescos. In general, an holomorphic parameter is not invariant by a change of variable, but we prove a theorem of stability of holomorphic families of frescos by a  change of variable and it implies that an holomorphic parameter gives again an holomorphic parameter by a change of variable.\\
  We construct here  two different kinds of holomorphic parameters which are (quasi-)\\invariant by change of variable. The first kind is associated to Jordan blocks of the monodromy with size  at least two. The second kind is associated to the semi-simple part of the monodromy and look like some "cross ratio" of eigenvectors. \\
 They allow, in the situation describe above, to associate to a given  (vanishing) relative de Rham cohomology class  some numbers, which will depend holomorphically of our data, and are independant of the choice of the local coordinate  near \ $0$ \ to study the Gauss-Manin connection of this degeneration of compact complex manifolds.

\parag{AMS Classification} 32 S 25, 32 S 40, 32 S 50.

\parag{Key words} Degenerating family of complex manifolds, Relative de Rham cohomology classes, filtered Gauss-Manin connection, Theme, Fresco, (a,b)-module, asymptotic expansion, vanishing period.

\tableofcontents

\section{Introduction.}

When we consider a proper holomorphic map \ $\tilde{f }: X \to C$ \ of a complex manifold \ $X$ \ on a smooth complex curve \ $C$ \  with a critical value at a point \ $0$ \ in \ $C$, the choice of a local coordinate near this point allows to dispose of an holomorphic function \ $f$. Then we may construct, using this function,  an (a,b)-modules structure on the cohomology sheaves of the formal completion (in \ $f$) \ of the complex of sheaves \ $(Ker\, df^{\bullet},d^{\bullet})$. These (a,b)-modules represent a filtered version of the Gauss-Manin connection of \ $f$. The most simple example of this construction is the Brieskorn module (see [Br.70]) of a function with an isolated singular point. See [B.08] for the case of a 1-dimensional critical locus. \\
But it is clear that this construction depends seriously on the choice of the function \ $f$ \ that is to say on the choice of the local coordinate near the critical point \ $0$ \ in the complex curve \ $C$.\\
The aim of the present paper is to study the behaviour of such constructions when we make a change of local coordinate near the origin. We consider the case of \ $[\lambda]-$primitive frescos, which are monogenic geometric (a,b)-modules corresponding  to a minimal filtered differential equation associated to a relative de Rham cohomology class on \ $X$ \ (see [B.09-a] and [B.09-b]).  \\
 An holomorphic parameter is a function on the set of isomorphism classes of frescos which behave holomorphically in an holomorphic family of frescos. In general, an holomorphic parameter is not invariant by a change of variable, but we prove a theorem of stability of holomorphic families of frescos by a  change of variable and it implies that an holomorphic parameter gives again an holomorphic parameter by a change of variable.\\
  We construct here  two different kinds of holomorphic parameters which are (quasi-)\\invariant by change of variable. The first kind is associated to Jordan blocks of the monodromy with size  at least two. The second kind is associated to the semi-simple part of the monodromy and look like some "cross ratio" of eigenvectors. \\
 They allow, in the situation describe above, to associate to a given  (vanishing) relative de Rham cohomology class  some numbers, which will depend holomorphically of our data, and are independant of the choice of the local coordinate  near \ $0$ \ to study the Gauss-Manin connection of this degeneration of compact complex manifolds.

\section{Change of variable in the ring \ $\hat{A}$.}

\subsection{Definition and first properties.}

We shall work with the \ $\C-$algebra
\begin{equation}
\hat{A} : = \sum_{\nu = 0} ^{+\infty} \  P_{\nu}(a).b^{\nu} \quad {\rm where} \quad P_{\nu} \in \C[[a]] \tag{1}
\end{equation}
with the product law defined by the following two conditions
\begin{enumerate}
\item \ $a.b - b.a = b^2$ ;
\item The right and left  multiplications by any \ $T \in \C[[a]]$ \ is continuous for the \ $b-$adic filtration of \ $\hat{A}$.
\end{enumerate}
The first condition implies  the identities 
 $$a.b^n = b^n.a + n.b^{n+1} \quad {\rm and} \quad  a^n.b = b.a^n + nb.a^{n-1}.b$$ 
  and with the second condition we obtain
\begin{align*}
& a.S(b) = S(b).a + b^2.S'(b) \quad {\rm for \ any} \quad S \in \C[[b]] \quad {\rm and} \tag{2a}\\
& T(a).b = b.T(a) + b.T'(a).b \quad {\rm for \ any} \quad T \in \C[[a]] \tag{2b}
\end{align*}

\begin{lemma}\label{debut}
We have the following properties of the \ $\C-$algebra \ $\hat{A}$.
\begin{enumerate}
\item For any \ $x,y$ \ in \ $\hat{A}$ \ we have \ $x.y - y. x \in \hat{A}.b^2$.
\item For each \ $n \in \mathbb{N}$ \ we have  \ $b^n.\hat{A} = \hat{A}.b^n$ \ and it is a two-sided ideal in \ $\hat{A}$.
\item We have \ $a.\hat{A} + b.\hat{A} = \hat{A}.a + \hat{A}.b$ \ and any element \ $x \in \hat{A}$ \ is invertible if and only if \ $x \not\in a.\hat{A} + b.\hat{A}$.
\end{enumerate}
\end{lemma}

\parag{Proof} We let the reader prove \ $1$ \ as an exercice. Let me prove \ $2$ \ by induction on \ $n \geq 0$. Assume that this is true for \ $n$ \ and consider \ $b^{n+1}.x$, for some \ $x \in \hat{A}$. We may write \ $b^n.x = y.b^n$ \ and then \ $b^{n+1}.x = b.y.b^n$. Now, using the fact that \ $b.y = y.b + z.b^2$  \ from \ $1$ \ we obtain \ $b^{n+1}.x = y.b^{n+1} + z.b^{n+2} = (y + z.b).b^{n+1}$.\\
Let \ $a.x + b. y $ \ for \ $x,y \in \hat{A}$, and write \ $x = x_0(a) + x'.b$ \ and \ $b.y = y'.b$ \ using \ $2$ \ for \ $n = 1$. Then \ $a.x + b.y = x_0(a).a + (a.x' + y').b$ \ and \ $a.\hat{A} + b.\hat{A} \subset  \hat{A}.a + \hat{A}.b$. Consider now  \ $u.a + v.b$ \ with \ $u,v \in \hat{A}$. Put \ $u = u_0(a) + b.u' $ \ using \ $(2b)$ \ and \ $v.b = b.v'$ \ using again \ $2$ \ for \ $n =1$. Then we have \ $u.a + v.b = a.u_0(a) + b.(u'.a + v')$.\\
Now consider  \ $x =  \sum_{\nu = 0} ^{+\infty} \  P_{\nu}(a).b^{\nu} $ \ with \ $P_0(0) \not= 0$. Then \ $P_0(a)$ \ is an invertible element of \ $\C[[a]]$. Now \ $P_0(a)^{-1}.x = 1 - b.y$ \ for some \ $y \in \hat{A}$. But now the serie  \ $ \sum_{n=0}^{+\infty} \ (b.y)^n $ \ converges for the \ $b-$adic filtration of \ $\hat{A}$ \ as \ $(b.y)^n$ \ lies in \ $b^n.\hat{A}$ \ thanks to \ $1$ \ and \ $2$. Then \ $1 -b.y$ \ is invertible in \ $\hat{A}$, and so is\ $x$.\\
Conversely, if \ $x$ \ is invertible, then so is the image of \ $x$ \ in \ $\C[[a]] = \hat{A}\big/b.\hat{A}$, and  we have \ $x_0(0) \not= 0$ \ and \ $x = x_0(0) + a.y + b.z $.
 Then \ $x$ \ does not lies in \ $a.\hat{A} + b.\hat{A}$. $\hfill \blacksquare$

\begin{lemma}\label{chgt. variable 0}
Let \ $\theta \in \C[[a]]$ \ such that \ $\theta(0) = 0$ \ and \ $\theta'(0) \not= 0$. Then the elements \ $\alpha : = \theta(a)$ \ and \ $\beta : = b.\theta'(a)$ \ satisfy the commutation relation in \ $\hat{A}$ :
$$ \alpha.\beta - \beta.\alpha = \beta^2 .$$
\end{lemma}

\parag{Proof} This is a simple computation using the identity in \ $(2b)$
\begin{align*}
& \alpha.\beta - \beta.\alpha  = \theta(a).b.\theta'(a) - b.\theta'(a).\theta(a) = b.\theta'(a).b.\theta'(a) = \beta^2
\end{align*}$ \hfill \blacksquare$\\

So for any such \ $\theta$ \ there exists an unique automorphism of \ $\C-$algebra \ $\Theta : \hat{A} \to  \hat{A} $ \ such that \ $ \Theta(1) = 1, \Theta(a) = \alpha , \Theta(b) = \beta $ \ which is continuous for the \ $b-$adic filtration and left and right \ $\C[[a]]-$linear. We have, when \ $x = \sum_{\nu = 0} ^{+\infty} \ P_{\nu}(a).b^{\nu} $
$$ \Theta(x) = \sum_{\nu = 0} ^{+\infty} \ P_{\nu}(\alpha).\beta^{\nu}. $$

\begin{defn}\label{chgt. variable 1}
We shall say that  \ $\theta \in \C[[a]]$ \ such that \  $\theta(0) = 0$ \ and \ $\theta'(0) \not= 0$ \ is  a {\bf change of variable} in \ $\hat{A}$. For any left \ $\hat{A}-$module \ $E$, we define a new \ $\hat{A}-$module \ $\theta_*(E)$, called {\bf the change of variable of \ $E$ \ for \ $\theta$}, by letting \ $\hat{A}$ \ act on \ $E$ \ via the automorphism \ $\Theta$.\\
Explicitely, as a \ $\C-$vector space, we have \ $\theta_*(E) = E$ \ and for \ $x \in \hat{A}$ \ and \ $e \in \theta_*(E) \simeq E$ \ we put \ $x.e : = \Theta(x)._Ee$, where on the right handside \ $\hat{A}$ \ acts on \ $e$ \ as an element of \ $E$.
\end{defn}

\parag{Convention} When we consider an element \ $x \in \theta_*(E)$ \ we write \ $a.x$ \ and \ $b.x$ \ for the action of \ $a$ \ and \ $b$ \ on the \ $\hat{A}-$module \ $\theta_*(E)$. Considering now \ $x$ \ as an element of \ $E$, this means that we look in fact at \ $\alpha.x$ \ and \ $\beta.x$ \ respectively. So we have to distinguish carefully if we are looking at \ $x$ \ as an element of \ $\theta_*(E)$ \ or as an element in \ $E$.

\begin{lemma}\label{chgt. variable 2}
For any \ $\hat{A}-$module \ $E$ \ and any change of variable \ $\theta$ \ we have
$$ b^n.\theta_*(E) = \beta^n.E =  b^n.E \quad \forall \ n \in \mathbb{N}.$$
\end{lemma}

\parag{Proof} We shall prove the statement by induction on \ $n \geq 0$. As the case \ $n = 0$ \ is clear, assume that the equality \ $\beta^n.E = b^n.E$ \ is proved and we shall prove it for \ $n+1$. The inclusion \ $\beta^{n+1}.E \subset b^{n+1}.E$ \ is easy because \ $\beta = b.\theta'(a)$, $b^n.E$ \ is stable by \ $\hat{A}$ \ as \ $b^n.\hat{A} = \hat{A}.b^n$ \ and our induction hypothesis allows to conclude.\\
Assume that \ $x \in E$. As \ $\theta'(0)$ \ is not \ zero, $\theta'(a)$ \ is invertible in \ $\hat{A}$, the element  \ $\theta'(a)^{-1}.b^n.x$ \ is in \ $b^n.E = \beta^n.E$ \ and so there exists \ $y \in E$ \ such that \ $\beta^n.y = \theta'(a)^{-1}.b^n.x$. Then \ $\beta^{n+1}.y = b.\theta'(a).\theta'(a)^{-1}.b^n.x = b^{n+1}.x$, which gives the desired  inclusion \ $b^{n+1}.E \subset \beta^{n+1}.E. \hfill \blacksquare$\\

\begin{cor}\label{chgt. variable 3}
If we have a \ $\hat{A}-$linear map \ $f : E \to F$ \ between two left \ $\hat{A}-$modules, then the {\bf same map} induces a \ $\hat{A}-$linear map between \ $\theta_*(E)$ \ and \ $\theta_*(F)$.
\end{cor}

\parag{Proof} The map \ $f$ \ is \ $\C[[a]]-$linear and also \ $\C[[b]]-$linear. So it is \ $\C[[\alpha]]-$linear. It is also \ $\C[\beta]-$linear because we have
$$ f(\beta.x) = f(b.\theta'(a).x) = b.\theta'(a).f(x)= \beta.f(x) .$$
It is now enough to see that \ $f$ \ is continuous for the \ $\beta-$adic filtration to conclude that \ $f$ \ is \ $\C[[\beta]]-$linear. But as the previous lemma shows that the \ $\beta-$adic filtration is equal to the \ $b-$adic filtration for which \ $f$ \ is continuous, the map \ $f$ \ is \ $\hat{A}-$linear from \ $\theta_*(E)$ \ to \ $\theta_*(F)$. $\hfill \blacksquare$

\parag{Remarks}
\begin{enumerate}
\item Thanks to the previous corollary, if we have a \ $\hat{A}-$linear map \ $f : E \to F$ \ between two left \ $\hat{A}-$modules, then the same map induces a \ $\hat{A}-$linear map between \ $\theta_*(E)$ \ and \ $\theta_*(F)$.  It is then obvious that any change of\ variable  $\theta_*$ \ defines a functor on the category of left \ $\hat{A}-$module into itself, which is the identity on morphisms, and it is an automorphism of this category. \\
Warning : it is not true in general that \ $\theta_*(E)$ \ is isomorphic to \ $E$ \ (see section 4.4).
\item Remark also that for any left \ $\hat{A}-$module \ $E$ \  the sub$-\hat{A}-$modules of the \ $\hat{A}-$module \ $\theta_*(E)$ \ are the same sub-vector spaces of the  given  vector space \ $E = \theta_*(E)$.
\item If the left \ $\hat{A}-$module \ $E$ \ is finitely generated, then so is \ $\theta_*(E)$. This a easy consequence of the fact that \ $\theta_*(\hat{A})$ \ is isomorphic to \ $\hat{A}$ \ as a left \ $\hat{A}-$module via the map \ $\Theta$.
\end{enumerate}

\begin{lemma}\label{chgt. variable 4}
Let \ $E$ \ be a left \ $\hat{A}-$module. Then \ $T_b(E)$, the \ $b-$torsion sub-$\C[[b]]-$module of \ $E$,  is a left  \ $\hat{A}-$submodule. For any change of variable \ $\theta$ \ we have the equality
 $$ \theta_*(T_b(E)) = T_b(\theta_*(E)) = T_{\beta}(E).$$
 If \ $E$ \ is finitely generated on \ $\C[[b]]$ \ then \ $\theta_*(E)$ \ is finitely generated on \ $\C[[b]]$
 \end{lemma}
 
\parag{Proof} First remark that for any left \ $\hat{A}-$module \ $E$ \ and for each \ $n \in \mathbb{N}$ \ the sub-vector space \ $Ker \, b^n$ \ is an \ $\hat{A}-$submodule of \ $E$. Let me prove this by induction on \ $n$. For \ $n = 0$ \ this is clear, For \ $n = 1$ \ it is an immediate  consequence of the formula \ $b.T(a) = T(a).b + b.T'(a).b$. So assume that \ $n \geq 2$ \ and the assertion true for \ $Ker \, b^{n-1}$. Let \ $x \in Ker \, b^n$ \ and \ $T \in \C[[a]]$. Then we have\ $b^n.T(a).x = b^{n - 1}.T(a).b.x + b^n.T'(a).b.x$, and \ $b.x$ \ is in \ $Ker \, b^{n-1}$ \ and also \ $T(a).b.x$ \ and \ $T'(a).b.x$ \ because our inductive assumption. We conclude that \ $T(a).x$ \ is in \ $Ker\, b^n$.
Now let us show that we have \ $Ker \, \beta^n = Ker \, b^n$ \ for each \ $n \in \mathbb{N}$. Again the case \ $n = 0$ \ is obvious and for \ $n=1$ \ we have \ $Ker \, \beta \subset Ker \, b$ \ because \ $\theta'(a)$ \ is invertible in \ $\hat{A}$ \ and \ $Ker \, b \subset Ker \, \beta $ \ because \ $Ker \, b$ \ is stable by \ $\C[[a]]$. Assume \ $n \geq 2$ \ and that  the equality \ $Ker \, \beta^{n-1} = Ker \, b^{n-1}$ \ is proved. If \ $b^n.x = 0$ \ then \ $b.x $ \ is in  \ $Ker \, b^{n-1} = Ker \, \beta^{n-1}$. Now \ $\beta^n.x = \beta^{n-1}.b.\theta'(a).x$ \ and \ $\theta'(a).b.x $ \ is in \ $Ker \, b^{n-1} = Ker \, \beta^{n-1}$ \ because we know that \ $Ker \, b^{n-1}$ \ is \ $\C[[a]]-$stable. So \ $Ker \, b^n \subset Ker \, \beta^n$.\\
Conversely, if \ $\beta^n.x = 0$ \ we have \ $\beta.x \in Ker \, \beta^{n-1} = Ker \, b^{n-1}$ \ and then 
$$b^{n-1}.\beta.x = b^n.\theta'(a).x = 0 .$$
 So \ $\theta'(a).x$ \ is in \ $Ker \, b^n$ \ which is \ $\C[[a]]-$stable, and \ $\theta'(a)$ \ is invertible in \ $\C[[a]]$. So \ $x$ \ is in \ $Ker \, b^n$, and we have the equality \ $Ker \, b^n = Ker \, \beta^n$.\\
Then the \ $b-$torsion of \ $E$ \ co{\"i}ncides with the \ $b-$torsion of \ $\theta_*(E)$.\\
As we know that \ $b.E = \beta.E$, the vector spaces \ $E\big/b.E$ \ and \ $  \theta_*(E)\big/b.\theta_*(E)$ \ co{\"i}ncide. So the finite dimension of \ $E\big/b.E $ \ implies that \ $\theta_*(E)$ \ is finitely generated as a \ $\C[[b]]-$module. $\hfill \blacksquare$\\

We introduce in [B.08] (see also the appendix of [B.11]) the notion of a small \ $\A-$module on the algebra
$$ \A : = \{ \sum_{\nu = 0}^{+\infty} \ P_{\nu}(a).b^{\nu} \} \quad {\rm where} \quad P_{\nu} \in \C[a] ,$$
with conditions analoguous for \ $\A$ \ analoguous to conditions 1 and 2  for \ $\hat{A}$.\\
In fact we were mainly interested in the case where these modules are naturally \ $\hat{A}-$modules. The following definition is the \ $\hat{A}-$analog. \\
First define for any \ $\hat{A}-$module \ $E$ 
$$  \hat{T}_a(E) : = \{ x \in E \ / \  \C[[b]].x \subset  T_a(E) \} $$
where \ $T_a(E)$ \ is the \ $a-$torsion of \ $E$. Note that \ $T_a(E)$ \ is not stable by \ $a$ \ in general, so it is not a sub$-\A-$module (and "a fortiori" not a sub$-\hat{A}-$module).

\begin{defn}\label{small}
We shall say that a \ $\hat{A}-$module\ $E$ \  is {\bf small} when it satisfies the following conditions :
\begin{enumerate}
\item  \ $E$ \ is a finitely generated \ $\C[[b]]-$module.
\item We have the inclusion \ $T_b(E) \subset \hat{T}_a(E)$.
\item  There exists \ $N \in \mathbb{N}$ such \ $a^N. \hat{T}_a(E) = 0 $.
\end{enumerate}
\end{defn}

\begin{prop}\label{change small}
The class of small \ $\hat{A}-$modules is stable by any change of variable.
\end{prop}

The proof of this proposition will be an easy consequence of the lemma \ref{chgt. variable 4} and the next lemma.

\begin{lemma}\label{torsion a}
For any \ $\hat{A}-$module \ $E$, $\hat{T}_a(E)$ \ is the maximal sub$-\hat{A}-$module contained in \ $T_a(E)$.
\end{lemma}

\parag{Proof} It is clear that any sub$-\hat{A}-$module of \ $E$ \ contained in \ $T_a(E)$ \ is contained in \ $\hat{T}_a(E)$. So the only point to prove is that \ $\hat{T}_a(E)$ \ is a sub$-\hat{A}-$module. Let \ $x$ \ be an element in \ $\hat{T}_a(E)$ \ and \ $z$ \ be an element of \ $\hat{A}$. We shall prove that \ $z.x$ \ is in \ $T_a(E)$ \ which is enough to conclude.\\
As \ $x$ \ is in \ $T_a(E)$, there exists \ $N_0 \in \mathbb{N}$ \ with \ $a^{N_0}.x = 0$. Now write
$$ z = \sum_{\nu = 0}^{+\infty} \ b^{\nu}.P_{\nu}(a) $$
where \ $P_{\nu}$ \ is in \ $\C[[a]]$. Put for each \ $\nu \geq 0 $
$$ P_{\nu}(a) = Q_{\nu}(a).a^{N_0} + R_{\nu}(a) $$
with \ $Q_{\nu} \in \C[[a]]$ \ and \ $R_{\nu} \in \C[a]$ \ a polynomial of degree at most \ $N_0 - 1$. Then we have 
 $$z.x = \sum_{\nu = 0}^{+\infty} \ b^{\nu}.R_{\nu}(a).x = \sum_{j=0}^{N_0-1} \ S_j(b).a^j.x  .$$
 But we may write in \ $\hat{A}$
 $$  \sum_{j=0}^{N_0-1} \ S_j(b).a^j =  \sum_{j=0}^{N_0-1} \ a^j.T_j(b) $$
 where each \ $T_j$ \ is in \ $\C[[b]]$. As \ $T_j(b).x$ \ is in \ $T_a(E)$, there exists an integer \ $N_1$ \ such that \ $a^{N_1}.T_j(b).x = 0 $ \ for each \ $j \in [0,N_0-1]$. Then we get
 $$a^{N_1}.z.x =  \sum_{j=0}^{N_0-1} \ a^{N_1+j}.T_j(b).x = 0 $$
 and \ $z.x$ \ is in \ $T_a(E)$. $\hfill \blacksquare$\\
 
 As \ $T_b(E)$ \ is a sub$-\hat{A}-$module of \ $E$, we have equivalence between the two inclusions
 \begin{enumerate}
 \item \  $ T_b(E) \subset T_a(E)$ 
 \item \ $ T_b(E) \subset \hat{T}_a(E) $.
 \end{enumerate}
 For any change of variable \ $\theta$ \ the equality \ $T_a(E) = T_{\alpha}(E) = T_a(\theta_*(E))$ \ is obvious. So the previous lemma implies the equality \ $ \hat{T}_a(E) = \hat{T}_a(\theta_*(E))$ \ as \ $E$ \ and \ $\theta_*(E)$ \ have the same sub$-\hat{A}-$modules, and so the conditions 1,  2 and 3 of the definition \ref{small} are stable by any  change of variable, thanks to the lemma \ref{chgt. variable 4}. This proves the proposition \ref{change small}. $\hfill \blacksquare$
 
 \bigskip

 \begin{defn}\label{regular}
 We shall say that a left \ $\hat{A}-$module \ $E$ \ is a {\bf simple pole (a,b)-module} when \ $E$ \ is free and finitely generated over \ $\C[[b]]$ \ and satisfies \ $a.E \subset b.E$.\\
 We shall say that a left \ $\hat{A}-$module \ $E$ \ is a {\bf regular (a,b)-module} when it is contained in a simple pole (a,b)-module. 
 \end{defn}
 
 \parag{Remarks}
  \begin{enumerate}
 \item A  regular (a,b)-module is free and finitely generated on \ $\C[[b]]$ \ by definition. 
 \item The previous terminology is compatible with the terminology on (a,b)-modules (see [B.93]) because an (a,b)-module (which is not "a priori" a \ $\C[[a]]-$module) is complete for the \ $a-$adic filtration when it is regular. So the action of \ $\hat{A}$ \ is well defined on a regular (a,b)-module.
 \end{enumerate}
  
\begin{lemma}\label{chgt. variable 5}
Let \ $E$ \ a simple pole (a,b)-module (resp. a regular (a,b)-module). Then for any change of variable \ $\theta$ \ the left \ $\hat{A}-$module \ $\theta_*(E)$ \ is a simple pole (a,b)-module (resp. a regular (a,b)-module).
\end{lemma}

\parag{Proof} It is enough to prove the assertion in the simple pole case because the change of variable of a sub$-\hat{A}-$module is a sub$-\hat{A}-$module, because the fact that \ $E$ \ has no \ $b-$torsion is preserved by change of variable and also the finite generation over \ $\C[[b]]$,  thanks to the  lemma \ref{chgt. variable 4}. Assume that \ $a.E \subset b.E$. Then we have \ $\alpha.E \subset b.E = \beta E$. So \ $\theta_*(E)$ \ has a simple pole. $\hfill \blacksquare$\\

Recall now that, for a given regular (a,b)-module \ $E$, there exists an unique "smallest" simple pole (a,b)-module \ $E^{\sharp}$ \ containing \ $E$. The precise meaning of this sentence  is that the \ $\hat{A}-$linear inclusion \ $i : E \to E^{\sharp}$ \ factorises in an unique way any  \ $\hat{A}-$linear map \ $f : E \to F$ \ where \ $F$ \ is a simple pole (a,b)-module. The construction of \ $E^{\sharp}$ \ is easy by considering the saturation of \ $E$ \ by \ $b^{-1}.a$ \ in any embbeding of \ $E$ \ in a simple pole (a,b)-module. A  canonical construction of \ $E^{\sharp}$ \ is obtained by the saturation of \ $E$ \ by \ $b^{-1}.a$ \ in \ $E \otimes_{\C[[b]]} \C[[b]][b^{-1}] $.

\begin{defn}\label{Bernstein 0}
Let \ $E$ \ be a regular (a,b)-module. We define the {\bf Bernstein polynomial} of \ $E$ \ as the minimal polynomial of the action of \ $-b^{-1}.a$ \ on the finite dimensional vector space \ $E^{\sharp}\big/b.E^{\sharp}$.
\end{defn}

\begin{prop}\label{Bernstein 1}
Let \ $E$ \ be a regular (a,b)-module and \ $\theta$ \ a change of variable. Then we have a canonical \ $\hat{A}-$linear isomorphism  \ $ \theta_*(E^{\sharp}) \to \theta_*(E)^{\sharp} $ \ and the Bernstein polynomial of \ $E$ \ and \ $\theta_*(E)$ \ co{\"i}ncide.
\end{prop}

\parag{Proof} As \ $\theta_*(E^{\sharp})$ \ is a simple pole (a,b)-module and \ $\theta_*(E)$ \ is regular, thanks to the previous lemma, the universal property of \ $ i : \theta_*(E) \to \theta_*(E)^{\sharp}$ \ implies, as the change of variable gives an \ $\hat{A}-$linear injection \ $\theta_*(i) : \theta_*(E) \to \theta_*(E^{\sharp})$, that there is an \ $\hat{A}-$linear factorization for \ $\theta_*(i) $
$$ \theta_*(E) \overset{i}{\to}  \theta_*(E)^{\sharp} \overset{j}{\to}  \theta_*(E^{\sharp}) .$$
It is then easy to see (using \ $\theta^{-1}$) \ that the injective map \ $\theta_*(i)$ \ has the same universal property than \ $i$. So \ $j$ \ is an isomorphism.\\
Now we have a natural isomorphism of vector spaces 
 $$ \theta_*(E)^{\sharp}\big/b.\theta_*(E)^{\sharp} \simeq \theta_*(E^{\sharp})\big/b.\theta_*(E^{\sharp}).$$
 To conclude, it is now enough to prove that the endomorphisms \ $b^{-1}.a$ \ on both sides are compatible with this isomorphism. This is a obvious consequence of the following fact : if \ $F$ \ is a simple pole (a,b)-module then \ $b^{-1}.a$ \ and \ $\beta^{-1}.\alpha$ \ induce the same endomorphism in \ $F\big/b.F$, because of  the relation
 $$ \beta^{-1}.\alpha = \theta'(a)^{-1}.(b^{-1}.a).(\theta(a)\big/a) $$
 and the fact that \ $a$ induces the \ $0$ \ map on \ $F\big/b.F$. $\hfill \blacksquare$

\begin{cor}\label{chgt. variable 6}
Let \ $E$ \ be rank 1 regular (a,b)-module. Then for any change of variable\ $\theta$ \ we have \ $\theta_*(E) \simeq E$.
\end{cor}

\parag{Proof} As any regular rank 1 (a,b)-module is isomorphic to \ $\hat{A}\big/(a - \lambda.b).\hat{A}$ \ for some \ $\lambda \in \C$ \ and as the corresponding Bernstein polynomial is \ $z + \lambda$, the previous proposition implies this result. $\hfill \blacksquare$\\

Remark that the isomorphism between \ $E_{\lambda}$ \ and \ $\theta_*(E_{\lambda})$ \ is unique up to a non zero constant, because \ $Aut_{\hat{A}}(E_{\lambda}) \simeq \C^*$. A precise description of this isomorphism will be given later on (see lemma \ref{rk 1 polyn. depend.} ).

\subsection{Some facts on \ $[\lambda]-$primitive frescos.}

\begin{defn}\label{fresco}
We call a geometric monogenic (a,b)-module a {\bf fresco}.
\end{defn}

First recall that is \ $F \subset E$ \ is a normal\footnote{A submodule \ $F \subset E$ \ is normal if \ $E\big/F$ \ is again free on \ $\C[[b]]$. Note that for \ $G \subset F \subset E$ \ it is equivalent to ask that \ $G$ \ is normal in \ $F$ \ and \ $F$ \ normal in \ $E$ \ or to ask that \ $G$ \ and \ $F$ \ are normal is \ $E$.}  submodule of a fresco \ $E$, then \ $F$ \ and \ $E\big/F$ \ are frescos. Moreover, if \ $E$ \ is \ $[\lambda]-$primitive, so are \ $F$ \ and  \ $E\big/F$.\\

We recall some basic facts on \ $[\lambda]-$primitive frescos which will be used later.

\parag{The principal Jordan-H{\"o}lder sequence} As any regular (a,b)-module, a \ $[\lambda]-$primitive fresco \ $E$ \ admits a Jordan-H{\"o}lder sequence, which is a sequence of normal submodules 
 $$0 = F_0 \subset F_1 \subset \cdots \subset F_{k-1} \subset F_k = E $$
 such that the quotients \ $F_{j+1}\big/F_j$ \ have rank 1. Then, for each \ $j$ \ we have \ $F_{j+1}\big/F_j \simeq E_{\lambda_j} \simeq \hat{A}\big/\hat{A}.(a - \lambda_j.b)$. \\
  It is proved in [B.11]  proposition 3.1.2  that any \ $[\lambda]-$primitive fresco admits  a J-H. sequence such that we have \ $ \lambda_1+1 \leq \lambda_2+2 \leq \cdots \leq \lambda_k + k$. Moreover such a J-H. sequence is unique (so the sub-modules \ $F_j$ \ are uniquely determined). It is called the {\bf principal J-H. sequence of \ $E$}. \\
  For any J-H. sequence of a \ $[\lambda]-$primitive  fresco \ $E$ \ the sequence of  numbers  \ $\mu_1, \dots, \mu_k$ \ associated to the consecutive quotients are such that the (unordered) set  \ $\{\mu_1+1, \dots, \mu_k+k\}$ \ is independant of the choice of the J-H. sequence. So writing these numbers in increasing order gives the previous numbers \ $\lambda_1+1, \dots, \lambda_k+k$ \ deduced from the principal J-H. sequence.\\
  
 We call {\bf fundamental invariants} of the rank \ $k$ \ $[\lambda]-$primitive fresco \ $E$ \ the ordered sequence of  numbers \ $\lambda_1, \dots, \lambda_k$ \ corresponding to the principal J-H. sequence. It is often convenient to put \ $\lambda_{j+1} : = \lambda_j + p_j -1$ \ for \ $j \in [1,k-1]$ \ and to give the fundamental invariants via the numbers \ $\lambda_1, p_1, \dots, p_{k-1}$. Note that \ $p_j$ \ are natural integers and that \ $\lambda_1$ \ is rational number such that \ $\lambda_1 > k-1$. \\
 For a \ $[\lambda]-$primitive fresco \ $E$ \ the Bernstein polynomial \ $B_E$ \ is the characteristic polynomial of the action of \ $-b^{-1}.a$ \ on \ $E^{\sharp}\big/b.E^{\sharp}$. So its degree is equal to the rank of \ $E$ \ over \ $\C[[b]]$.\\
  The roots of the  Bernstein polynomial of a \ $[\lambda]-$primitive fresco \ $E$ \ are related to fundamental invariants  \ $\lambda_1, \dots, \lambda_k$ \ by the formula
 $$ B_E(z) = (z + \lambda_1+1-k).(z + \lambda_2 + 2-k) \dots (z + \lambda_k + k - k) .$$
 
\parag{Remark} The uniqueness of \ $F_1$ \ for the principal J-H. sequence of a \ $[\lambda]-$primitive fresco \ $E$ \ implies that the dimension of \ $Ker\, (a - \lambda_1.b)$ \ is a one dimensional vector space and that \ $F_1  = \C[[b]].e_1$ \ where \ $e_1$ \ is any non zero vector in this kernel.

\subsection{Embedding frescos.}

\begin{defn}\label{asympt. 0}
Let \ $\lambda$ \ be a rational number in \ $]0,1]$ \ and \ $N$ \ an integer. We define the left \ $\hat{A}-$module \ $\Xi_{\lambda}^{(N)}$ \ as the free \ $\C[[b]]$ \ module generated by \ $e_0, \dots,e_N$ \ with the action of \ $a$ \ define by the following rules :
\begin{enumerate}
\item \ $ a.e_0 = \lambda.b.e_0 $ ;
\item For \ $j \in [1,N]$ \  \ $a.e_j = \lambda.b.e_j + b.e_{j-1} $ ;
\item The left action of \ $a$ \ is continuous for the \ $b-$adic topology of \ $\Xi_{\lambda}^{(N)}$. As \ $a.\Xi_{\lambda}^{(N)} \subset b.\Xi_{\lambda}^{(N)}$ \ the action of \ $a$ \ extends to an action of \ $\C[[a]]$\\
\end{enumerate}
\end{defn}

It is an easy computation to check that \ $a.b - b.a = b^2$ \ on \ $\Xi_{\lambda}^{(N)}$, so we have defined a simple pole (a,b)-module of rank \ $N+1$.\\
Remark that \ $\Xi_{\lambda}^{(N)}$ \ corresponds to the formal asymptotic expansions of the type
$$ \sum_{j = 0}^N \  S_j(b).s^{\lambda-1}.\frac{(Log\, s)^j}{j!} = \sum_{j = 0}^N \  T_j(s).s^{\lambda-1}.\frac{(Log\, s)^j}{j!} $$
where the action of \ $a$ \ is the multiplication with \ $s$ \ and \ $b$ \ is the "primitive without constant".

Now the following embedding theorem is an elementary by-product of the theorem 4.2.1 in [B.09-b] (A more precise version for \ $[\lambda]-$primitive frescos is given in [B.11].)

\begin{thm}\label{embed.}
Let \ $E$ \ be any rank \ $k$ \ \ $[\lambda]-$primitive fresco. Then there exists an integer \ $N \leq k$, a complex vector space \ $V$ \ of dimension \ $\leq k$ \ and an embedding of \ $E$ \ in the left \ $\hat{A}-$module \ $\Xi_{\lambda}^{(N)} \otimes V$.
\end{thm}

Note that \ $a$ \ and \ $b$ \ acts on \ $\Xi_{\lambda}^{(N)} \otimes V$ \ as \ $a \otimes \id_V$ \ and \ $b \otimes \id_V$.\\

Let me recall two basic definitions and the general structure theorem proved in [B.11]  for \ $[\lambda]-$primitive frescos.

\begin{defn}\label{theme + ss.}
Let \ $E$ \ be a \ $[\lambda]-$primitive fresco. We say that \ $E$ \ is a {\bf theme} if there exists an embedding of \ $E$ \ in \ $\Xi_{\lambda}^{(N)}$ \ for some integer \ $N$.\\
We say that \ $E$ \ is {\bf semi-simple} when any \ $\hat{A}-$linear map \ $\varphi : E \to \Xi_{\lambda}^{(N)}$ \ has rank at most equal to\ $1$.
\end{defn}

\parag{Remarks} 
\begin{enumerate}
\item Of course any rank 1 fresco is at the same time a theme and semi-simple. But it is clear that in rank \ $\geq 2$ \ a fresco cannot be at the same time a theme and semi-simple.
\item Themes and semi-simple frescos have stability properties. See the [B.09-b] and [B.11].
\end{enumerate}

\begin{thm}[{See prop. 4.1.6 and lemma 4.2.1 of [B.11]} ]\label{structure}
Let \ $E$ \ be a \ $[\lambda]-$primitive fresco. There exists a biggest semi-simple submodule \ $S_1(E)$ \ in \ $E$. It is normal and the quotient \ $E\big/S_1(E)$ \ is a theme.\\
Dualy there also exists a smallest submodule \ $\Sigma^1(E)$ \ in \ $E$ \ which is normal and such that the quotient \ $E\big/\Sigma^1(E)$ \ is semi-simple. The fresco \ $\Sigma^1(E)$ \ is a theme.\\
The isomorphism class of \ $E$ \ determines the isomorphism classes of \ $S_1(E), \Sigma^1(E)$, $ E\big/S_1(E)$ \ and \ $E\big/\Sigma^1(E)$.
\end{thm}

\parag{Some more facts on frescos} We defined in \ $[B.11]$ \ the {\bf ss-depth} \ $d(E)$ \ of a \ $[\lambda]-$primitive fresco \ $E$ \ as the maximal  rank of a normal sub-theme in \ $E$. It is also the maximal rank of a quotient theme of \ $E$. It is shown that the rank of \ $S_1(E)$ \ is \ $rk(E)-d(E) + 1$ \ and the rank of \ $\Sigma^1(E)$ \ is \ $d-1$. So \ $E$ \ is a theme if and only if \ $d(E) = rk(E)$, and \ $E$ \ is semi-simple if and only if \ $d(E) = 1$. It is also shown that \ $E\big/S_1(E)$ \ and \ $\Sigma^1(E)$ \ themes.\\

The following subsection will give some complements to the results in [B.11]  which will be usefull later for the study  of holomorphic families of \ $[\lambda]-$primitive frescos.\\

\subsection{Complements.}

\begin{defn}\label{socle}
Let \ $E$ \ be a \ $[\lambda]-$primitive fresco of rank \ $k$. Then define \ $L(E)$ \ as \ $0$ \ if \ $E$ \ is semi-simple, and as \ $L(E) : = S_1(E) \cap \Sigma^1(E)$ \ when \ $d(E) \geq 2$.
\end{defn}

\parag{Remark}
If \ $E$ \ is not semi-simple, as \ $\Sigma^1(E)$ \ is a normal rank \ $(d-1)-$theme and \ $S_1(E)$ \ is normal semi-simple, $L(E)$ \ is a rank \ $1$ \ normal submodule of \ $E$.\\

\begin{lemma}\label{d decroit}
Let   \ $E$ \ be a \ $[\lambda]-$primitive fresco of rank \ $k$. Assume \ $d(E) \geq 2$. Then we have
$$ d(E\big/L(E)) = d(E) - 1 .$$
\end{lemma}

\parag{proof} We shall argue by contradiction. Let \ $T$ \ be a normal sub-theme of rank \ $d : = d(E)$ \ in \ $E\big/L(E)$ \ and \ $\tilde{T}$ \ its pull-back  in \ $E$. Because \ $\tilde{T} \subset E,  d(\tilde{T}) \leq d$ \ and because \ $T$ \ is a rank \ $d$ \  quotient theme of \ $\tilde{T}$ \ we have   \ $d(\tilde{T}) \geq  d$ \ and so \ $d(\tilde{T}) = d$. Thanks again  to [B.11],  there exists a normal sub-theme \ $T' \subset \tilde{T}$ \ of rank \ $d$. Now we shall prove that  \ $E\big/T'$ \ is semi-simple. If the rank of \ $E$ \ is \ $d+1$ \ this is obvious. So we may assume that \ $rk(E) \geq d+1$.
If \ $E\big/T'$ \ is not semi-simple, choose  a J-H. sequence of \ $E$ \ which begins by the J-H.sequence of \ $T'$ \ and such that \ $F_{d+2}\big/F_d$ \ is a rank \ $2$ \ theme. As \ $F_d = T'$ \ and \ $E\big/T'$ \ is not semi-simple, this is possible thanks to theorem 5.1.1 in [B.11]. Then the number of non commuting index in this J-H. sequence is \ $\geq d$ \ as \ $j = 1, \dots, d-1$ \ and \ $j = d+1$ \ are commuting indices. So, thanks to lemma 4.3.4 in [B.11], we have \ $d(E) \geq d+1$ \ which is a contradiction.\\
Now, as \ $E\big/T'$ \ is semi-simple, we have \ $\Sigma^1(E) \subset T'$ \ and \ $L(E) = F_1(T')$. As \ $d \geq 2$ \ we have \ $F_1(T') \subset F_{d-1}(T')$ \ and so \ $\tilde{T}\big/F_{d-1}(T') \to T$ \ is injective ; we obtain that \ $\tilde{T}\big/F_{d-1}(T')$ \ is a rank \ $2$ \ theme. Then \ $\tilde{T}$ \ is a rank \ $d+1$ \  theme thanks to \ref{} (le theme) and this contradicts the fact that \ $d(E) = d$. $\hfill \blacksquare$

\parag{Remarks}
\begin{enumerate}[i)]
\item  Conversely, if we have a normal rank \ $1$ \ submodule \ $F \subset E$ \ in a \ $[\lambda]-$primitive fresco \ $E$ \ such that \ $d(E) \geq 2$ \ satisfying \ $d(E\big/F) = d(E) - 1$ \ then we have \ $F = L(E)$. Indeed let \ $T \subset E$ \ be a normal rank \ $d$ \ theme in \ $E$. We have \ $L(E) = F_1(T)$ \ (see the remark above) and \ $T\big/F\cap T \to E\big/F$ \ is injective. As \ $E\big/F$ \ cannot contain a rank \ $d$ \ theme this implies \ $F \cap T = F$. But \ $T$ \ as an unique normal rank \ $1$ \ sub-module which is \ $F_1(T) = L(E)$. So \ $F = L(E)$.
\item For any normal sub-theme \ $T$ \ of rank \ $\geq 2$ \ in \ $E$, we have  the injection
$$ T\big/(T \cap \Sigma^1(E)) \to E\big/\Sigma^1(E) $$
between a theme and a semi-simple  fresco. So its rank is \ $\leq 1$ \ and then we have \ $T \cap  \Sigma^1(E) = F_1(T) = F_1(\Sigma^1(E))$. This means that \ $L(E) = F_1(T)$ \ for any such \ $T$.\\
\end{enumerate}

 The next lemma is very easy but will be usefull later on.
 
 \begin{lemma}\label{descente}
 Let \ $E$ \ be a \ $[\lambda]-$primitive rank \ $k$ \  fresco such that \ $d(E) \geq 2$. Let \ $F$ \ be a normal sub-module of \ $E$,  and assume that \ $d(F) \geq 2$. Then we have \ $L(E) = L(F)$.
 \end{lemma}
 
  \parag{proof} Let \ $T \subset F$ \ be a normal rank \ $2$ \ theme. Then  \ $F_1(T) = L(F)$, but also \ $F_1(T) = L(E)$. $\hfill \blacksquare$.

 \parag{Remarks}
 \begin{enumerate}
 \item So any normal sub-module which is not semi-simple (i.e. not contained in \ $S_1(E)$) satisfies \ $L(F) = L(E)$.
 \item When \ $d(E) \geq 3$ \ the hypothesis of the previous lemma is always satisfied by any corank \ $1$ \ normal sub-module \ $F$ \ in \ $E$.
 \item When the hypothesis of the previous lemma is not satisfied  by some normal corank \ $1$ \ sumodule \ $F$ \ in \ $E$ \ we have \ $d(E) = 2$ \ and \ $F = S_1(E)$.
 \end{enumerate}

We shall give now a new proof of the theorem \ref{embed.} 

\begin{prop}\label{delta}
Let \ $E$ \ be a \ $[\lambda]-$primitive fresco with rank \ $k$ \ and fundamental invariants \ $\lambda_1, \dots, \lambda_k$. Then for \ $\mu \in [\lambda], \mu \geq \lambda_k + k-1$ \ the number
$$ \delta(E) : = \dim_{\C}\big[ Ker\,(a -\mu.b) : E \to E \big]  $$
is independant of the choice of \ $\mu$ \ and is equal to \ $k - d(E) + 1 $ \ where \ $d(E)$  \ is the ss-depth of \ $E$ \ (see [B.11]).
\end{prop}

\parag{Proof} If we have an \ $x \in E$ \ such that \ $(a - \mu.b).x = 0$ \ for some \ $\mu \in [\lambda]$ \ and \ $x \not\in b.E$, then there exists \ $j \in [1,k]$ \ with \ $\mu = \lambda_j+j-1$. This is because \ $\C[[b]].x \simeq E_{\mu}$ \ is then a rank \ $1$ \ normal submodule of \ $E$ \ and so is the first term of a J-H.sequence of \ $E$. So we have \ $\mu \leq \lambda_k+k-1$.\\
So if \ $x$ \ is in \ $Ker\,(a -\mu.b)$ \ with \ $\mu \geq \lambda_k+k$ \ we have \ $x \in b.E$ \ and so there exists \ $q \geq 1$ \ such that \ $x$ \ lies in \ $b^q.Ker\,(a -(\mu-q).b)$. This implies that for any \ $q \geq 0$ \ we have \ $b^q : Ker\,(a - (\lambda_k+k-1).b) \to Ker\,(a -(\lambda_k+k+q-1).b)$ \ is surjective. As \ $b$ \ is injective, it is an isomorphism an this proves our first assertion.\\
To prove that \ $\delta(E) = k-d(E) + 1$ \ we shall use the structure theorem 5.1.1  of [B.11] thanks to the following lemma.

\begin{lemma}\label{delta + 1}
Let \ $E$ \ be a \ $[\lambda]-$primitive fresco and let \ $F_j, j \in [1,k]$ \ be a J-H. sequence of \ $E$. Assume that for some \ $j \in [0,k-2]$ \ the quotient \ $F_{j+2}\big/F_j$ \ is a theme. Then we have
$$ \delta(F_{j+1}) = \delta(F_{j+2}) .$$
\end{lemma}

\parag{Proof} Fix \ $\mu \in [\lambda]$ \ large enough, and consider \ $x \in F_{j+2}$ \ such that \\
 $(a - \mu.b).x = 0$. Let \ $q \in \mathbb{N}$ \ such that \ $x = b^q.y$ \ and \ $y \not\in b.F_{j+2}$. Then \ $\C[[b]].y$ \ is a normal rank \ $1$ \ submodule of \ $F_{j+2}$. So its image in the rank \ $2$ \ theme \ $F_{j+2}\big/F_j$ \ is contained in \ $F_{j+1}\big/F_j$ \ which is the only normal rank \ $1$ \ submodule of \ $F_{j+2}\big/F_j$. So \ $y$ \ is in \ $F_{j+1}$ \ and also \ $x$ \ because \ $F_{j+1}$ \ is normal. This proves that \ $\delta(F_{j+1}) = \delta(F_{j+2}) . \hfill \blacksquare$\\

\parag{End of proof of proposition \ref{delta}} Using the theorem 5.11 of [B.11] we have a J-H. sequence \ $F_j, j \in [1,k]$ \ for \ $E$ \ such that \ $S_1(E) = F_{k-d(E)+1}$ \ and such that \ $F_{j+2}\big/F_j$ \ is a rank \ $2$ \ theme for each \ $j \geq k-d(E)$. Then the previous lemma gives
$$ \delta(F_{k-d(E)+1}) = \delta(F_{k-d(E)+2}) = \cdots = \delta(F_k) = \delta(E).$$
To conclude it is enough to show that for a semi-simple fresco we have \ $\delta(E) = rk(E)$. But this is clear because assuming \ $E$ \ semi-simple with fundamental invariants \ $\lambda_1, \cdots, \lambda_k$ \   we know that for each \ $j \in [1,k]$ \ we have an \ $x_j \in E \setminus b.E$ \ such that \ $(a - (\lambda_j+j-1).b).x_j = 0$, and we also know that the sequence \ $\lambda_j + j$ \ is strictly increasing (see corollary 4.1.3 and proposition 4.1.4 of [B.11]). It is then easy to see that the element \ $b^{q_j}.x_j, j \in [1,k]$ \ are linearly independant in \ $Ker\,(a - \mu.b)$ \ for \ $\mu \gg 1, \mu \in [\lambda]$ \ where we put \ $\lambda_j+j-1+q_j = \mu $. Now, as \ $F_{k-d(E)+1} = S_1(E)$ \ is semi-simple, we obtain the announced formula.$\hfill \blacksquare$\\

\parag{Remark} In the embedding theorem \ref{embed.} recalled above the number \ $k - d(E) + 1$ \ is the minimal dimension of a complex vector space \ $V$ \ such we may embed the \ $[\lambda]-$primitive fresco \ $E$ \ in \ $\Xi^{(N)}_{\lambda} \otimes V$ \ for \ $N $ \ large enough. But if we have an embedding \ $ E \to \Xi^{(N)}_{\lambda} \otimes V$ \ it induces an injective linear map \ $Ker\, (a - \mu.b) \to s^{\mu-1}\otimes V$. So we obtain immediately that \ $\dim_{\C} V \geq \delta(E)$. But the converse needs  more work, as we have seen.\\

\begin{lemma}\label{tech.}
Let \ $E$ \ be a rank \ $k$  \ $[\lambda]-$primitive fresco such that \ $d(E) = 2$. Assume \ $E \big/F_{k-2}$ \ is a theme, where  \ $F_{j}$ \ denotes  the rank \ $j$ \ sub-module of the principal J-H. sequence of \ $E$. Then we have \ $F_{k-1} = S_1(E)$.
\end{lemma}

\parag{proof} First we shall prove that \ $F_{k-2}$ \ is semi-simple, so contained in \ $S_1)E)$. If this is not the case, we may find a J-H . sequence of \ $F_{k-2}$ \ with at least one non commuting index. When we complete  this J-H . sequence of \ $F_{k-2}$ \ with \ $F_{k-1}$ \ and \ $F_k : = E$ \ we obtain a J-H . sequence of \ $E$ with at least two non commuting indices as \ $E\big/F_{k-2}$ \ is a rank \ $2$ \ theme. Then the proposition ?? of [B.11] gives \ $d(E) \geq 3$ \ contradicting our assumption.\\
So \ $F_{k-2} \subset S_1(E)$ \ and the rank \ $1$ \ quotient \ $S_1(E)\big/F_{k-2}$ \ injects in the rank \ $2$ \ theme  \ $E\big/F_{k-2}$. So its image is in \ $F_1(E\big/F_{k-2}) = F_{k-1}\big/F_{k-2}$. Then we have the inclusion \ $S_1(E) \subset F_{k-1}$. But these two sub-modules are normal and of same rank \ $k-1$. So they are equal. $\hfill \blacksquare$\\

We recall here for the convenience of the reader two results about commutation in the ring \ $\hat{A}$\footnote{in fact in the subring \ $\C[[b]][a]$} that we shall use (with holomorphic parameters) later on. The reader will find proofs respectively in [B.09-a] lemma 3.5.1 and in [B.11] lemma 3.2.1 and corollary 3.2.2.

\begin{lemma}[Commuting lemma]\label{ com.lem.}
Let \ $\lambda \not= \mu$ \ be a complex numbers and \ $S \in \C[[b]]$ \ such that \ $S(0) = 1$. If the difference \ $\delta = \lambda-\mu$ \ is a positive integer, assume that the coefficient of \ $b^{\delta}$ \ in \ $S$ \ is zero. Then we have in \ $\hat{A}$ \ the equality :
$$ (a - \mu.b).S^{-1}.(a - (\lambda-1).b) = U^{-1}.(a - \lambda.b).U.S^{-1}.U.(a - (\mu-1).b).U^{-1} .$$
\end{lemma}

\begin{prop}[Standard computation]\label{stand. comp.}
Let \ $p_1, p_2$ \ be positive integers and \ $\lambda_1$ \ be a rational number \ $\lambda_1 > 2$ ; put \ $\lambda_{j+1} = \lambda_j + p_j - 1$ \ for \ $j = 1,2$. For \ $S_1, S_2 \in \C[[b]]$ \ with \ $S_1(0) = S_2(0) = 1$ \ define 
$$ P : = (a - \lambda_1.b).S_1^{-1}.(a -\lambda_2.b).S_2^{-1}.(a - \lambda_3.b) .$$
Assume that the coefficient of \ $b^{p_1} $ \ in \ $S_1$ \ vanishes and that the coefficient  \ $\alpha$ \ of \ $b^{p_2}$ \ does not vanishes. Then there exists an unique solution \ $U \in \C[[b]]$ \ of the differential equation \ $b.U' = p_1.(U - S_1)$ \ such that  the coefficient of \ $b^{p_1+p_2}$ \ in \ $U.S_2$ \ vanishes. Then we have (applying two times the commuting lemma)
$$ P = U^{-1}.(a - (\lambda_2+1).b).S_1^{-1}.U^2.V^{-1}.(a - (\lambda_3+1).b).U^{-1}.S_2^{-1}.V^2.( a-(\lambda_1-1).b).V^{-1}.$$
Moreover  the coefficient of \ $b^{p_2}$ \ in  \ $S_.U^{-2}.V$ \ is equal to \ $(p_1+p_2).\alpha.p_1$.
\end{prop}

For the convenience of the reader we recall here the computation of the coefficient of \ $b^{p_2}$ \ in  \ $S_.U^{-2}.V$ : from the two differential equations we deduce, with \ $Z : = U^{-1}$
\begin{align*}
& b.Z' = -p_1.(Z - S_1.U^{-2}) \quad {\rm and } \quad b.Z'.V = -p_1.(Z.V - S_1.U^{-2}.V) \\
& b.V'.Z = (p_1+p_2).(V.Z - S_2) \quad {\rm then \ adding} \\
& b.(Z.V) - p_2.Z.V = p_1.S_1.U^{-2}.V - (p_1+p_2).S_2
\end{align*}
and, as the left hanside of the last equality has no term in \ $b^{p_2}$, the right handside also. $\hfill \blacksquare$

\subsection{Change of variable for themes and semi-simple frescos.}

Our next proposition shows that the (a,b)-module of \ $[\lambda]-$primitive asymptotic expansions \ $\Xi_{\lambda}^{(N)}$ \ is stable by change of variable. The proof will use the following (more or less classical) lemma.

\begin{lemma}\label{tech.lin.}
Let \ $N $ \ be the endomorphism of \ $\C^n$ \ define by \ $N(e_j) = e_{j+1}$ \ for \ $j \in [1,n-1]$ \ and \ $N(e_n) = 0$. Then for any \ $\nu \in \C^*$ \ the  linear map
$$ f :  End(\C^n) \to End(\C^n) $$
given by \ $f(S) : = S.N - N.S + \nu.S $ \ is an isomorphism.
\end{lemma}

\parag{proof} We shall prove by induction on \ $n \in \mathbb{N}^*$ \ that \ $f$ \ is injective. The case \ $n = 1$ \ is obvious, assume the result proved for \ $n-1 \geq 1$ \ and we shall prove it for \ $n$.
Consider \ $S \in End(\C^n)$ \ in \ $Ker\, f$. First remark that \ $-N(S(e_n)) + \nu.S(e_n) = 0 $, and, as \ $N$ \ has no non zero eigenvalue, it implies that \ $S(e_n) = 0$. So using the induction via the isomorphism \ $\C^n\big/\C.e_n \simeq \C^{n-1}$, the image of \ $S$ \ is contained in \ $\C.e_n$. Put \ $S(e_j) = \lambda_j.e_n$ \ for \ $j \in [1,n]$. Then using again the condition \ $f(S) = 0$, we obtain\ $\lambda_{j+1}.e_n + \nu.\lambda_j.e_n = 0 $. So we obtain \ $\lambda_j = (-\nu)^{j-1}.\lambda_1$. Now the vanishing of \ $\lambda_n$ \ implies that all \ $\lambda_j$ \ vanish, and so \ $S = 0$. $\hfill \blacksquare$\\

\begin{cor}\label{uniq. dev.}
Let \ $E$ \ be a simple pole (a,b)-module of rank \ $n$ \ such that the map induced by \ $a$ \ on \ $E\big/b^2.E$ \ is given in a suitable basis by
$$ a.e = b.( \lambda.Id_n + N).e $$
where \ $N$ \ is the principal nilpotent endomorphism of the previous lemma. Then there exists a \ $\C[[b]]-$basis \ $\varepsilon : = (\varepsilon_1, \dots, \varepsilon_n)$ \  of \ $E$ \ in which we have
$$ a.\varepsilon = b.( \lambda.Id_n + N).\varepsilon .$$
\end{cor}

\parag{Proof} We begin with a \ $\C[[b]]-$basis \ $e : = (e_1, \dots, e_n)$ \  of \ $E$ \ such that we have, thanks to our hypothesis,
$$ a.e = b.(\lambda.Id_n + N).e + b^2.Z.e $$
where \ $Z$ \ is in \ $End(\C^n) \otimes \C[[b]]$. We look for \ $S \in End(\C^n) \otimes \C[[b]]$ \ such that \ $S(0) = Id_n$ \ and that the \ $\C[[b]]-$basis \ $\varepsilon : = S.e$ \ satisfies our requirement. This gives the equation
\begin{equation*}
a.S.e = S.a.e + b^2.S'.e = S.(\lambda.Id_n + N).e + S.b^2.Z.e + b^2.S'.e = b.( \lambda.Id_n + N).S.e .
\end{equation*}
So we want to solve the equation
\begin{equation*}
S.N - N.S + b.S.Z + b.S' = 0 \tag{1}
\end{equation*}
Writting \ $S : = \sum_{\nu \geq 0}^{+\infty} \ S_{\nu}.b^{\nu}$ \ and \ $Z : = \sum_{\nu \geq 0}^{+\infty} \ Z_{\nu}.b^{\nu}$ \ we want to solve the recursion system
\begin{equation*}
S_{\nu}.N - N.S_{\nu} + \nu.S_{\nu} = - \sum_{j = 0}^{\nu-1} \ S_{\nu -j-1}.Z_j \quad \forall \nu \geq 1 \quad {\rm with} \quad \ S_0 = Id_n. \tag{2}
\end{equation*}
The previous lemma gives existence and uniqueness of the solution \ $S$. $\hfill \blacksquare$\\

\begin{prop}\label{chang. dev.}
Fix  \ $\lambda \in ]0,1] \cap \mathbb{Q}$ \ and \ $p \in \mathbb{N}$. Then for any change of variable \ $\theta$ \ the (a,b)-module \ $\theta_*(\Xi^{(p)}_{\lambda})$ \ is isomorphic to \ $\Xi^{(p)}_{\lambda}$.
\end{prop}

\parag{Proof} As the (a,b)-module \ $E : = \Xi^{(p)}_{\lambda}$ \   is a simple pole (a,b)-module of rank \ $n : = p+1$\ such that the map  induced by \ $a$ \ on \ $E\big/b^2.E$ \ is given in a suitable basis by
$$ a.e = b.( \lambda.Id_n + N).e $$
where \ $N$ \ is the principal nilpotent endomorphism of the lemma above, it is enough to prove that the map induced by \ $\alpha$ \ on \ $E\big/b^2.E \simeq E\big/\beta^2.E$ \ is again given by \ $\beta.( \lambda.Id_n + N)$. But we have \ $\alpha = \theta_1.a + b^2.E$ \ and \ $\beta = \theta_1.b + b^2.E$, which allow to conclude. $\hfill \blacksquare$.\\

Our next  result  gives the the stability of the notions of theme and of semi-simple fresco by any change of variable.

\begin{prop}\label{stability}
Let \ $E$ \ be a \ $[\lambda]-$primitive theme (resp. of a \ $[\lambda]-$primitive semi-simple fresco) and let \ $\theta$ \ be any change of variable. Then \ $\theta_*(E)$ \ is a \ $[\lambda]-$primitive theme (resp. of a \ $[\lambda]-$primitive semi-simple fresco).
\end{prop}

\parag{Proof} In the case of a \ $[\lambda]-$primitive theme, the result is an easy consequence of the fact that we have an isomorphism (see proposition \ref{chang. dev.})
$$ \theta_*(\Xi_{\lambda}^{(N)}) \simeq \Xi_{\lambda}^{(N)}.$$
In the case of a semi-simple \ $[\lambda]-$primitive fresco, consider an \ $\hat{A}-$linear map 
$$\varphi : \theta_*(E) \to  \Xi_{\lambda}^{(N)} .$$
Applying the change of variable \ $\theta^{-1}$ \ we obtain an \ $\hat{A}-$linear map
$$ \psi : E \to (\theta^{-1})_*(\Xi_{\lambda}^{(N)}) \simeq \Xi_{\lambda}^{(N)}  $$
which must have rank \ $\leq 1$. So \ $\varphi$ \ has also a rank \ $\leq 1$ \ and \ $\theta_*(E)$ \ is semi-simple. $\hfill \blacksquare$\\

\begin{cor}\label{chgt. variable struct.}
Let \ $E$ \ be a \ $[\lambda]-$primitive fresco and \ $\theta$ \ a change of variable. Then we have natural isomorphisms
$$ u : \theta_*(S_1(E)) \to S_1(\theta_*(E))\quad {\rm and} \quad v :  \Sigma^1(\theta_*(E)) \to \theta_*(\Sigma^1(E)) .$$
So the numbers \ $d(E)$ \ and \ $\delta(E)$ \ are invariant by any change of variable.
\end{cor}

\parag{Proof} As we know that any change of variable of a \ $[\lambda]-$primitive theme (resp. of a \ $[\lambda]-$primitive semi-simple fresco) is a \ $[\lambda]-$primitive theme (resp. \ $[\lambda]-$primitive semi-simple fresco) we have natural maps \ $u $ \ and \ $v$. But extremality properties allow immediately to conclude. $\hfill \blacksquare$\\

Note also that  for any change of variable \ $\theta$ \  we have also natural isomorphisms
$$ \theta_*(E\big/S_1(E)) \simeq \theta_*(E)\big/\theta_*(S_1(E)) \quad {\rm and} \quad   \theta_*(E\big/\Sigma^1(E)) \simeq \theta_*(E)\big/\theta_*(\Sigma^1(E)) .$$

\section{Holomorphic families of \ $[\lambda]-$primitive  frescos.}

\subsection{Holomorphic families and J-H. principal sequence.}

\begin{defn}\label{Hol. 1}
Let \ $X$ \ be a reduced complex space. We define the sheaf \ $\Xi_{\lambda, X}^{(N)}$ \ as the free sheaf of \ $\mathcal{O}_X[[b]]-$module with basis \ $e_0, \dots, e_N$. We define the action of \ $\hat{A}$ \ on this sheaf by the condition 1. 2. and 3. of the definition \ref{asympt. 0}. We obtain in this way a structure of left \ $\hat{A}-$module on this sheaf. Will shall say that a sheaf on \ $X$ \ with these two compatible structures,  an \ $\mathcal{O}_X[[b]]-$module and an \ $\hat{A}-$left module, is a \ $\hat{A}_X-$sheaf.\\
For any finite dimensional complex vector space we define the \ $\hat{A}_X-$sheaf \ $\Xi_{\lambda, X}^{(N)}\otimes V$ \ as the direct sum of \ $\dim_{\C}(V)$ \ copies of the sheaf \ $\Xi_{\lambda, X}^{(N)}$. So the action of \ $a$ \ and \ $b$ \ on \ $\Xi_{\lambda, X}^{(N)}\otimes V $ \ are given by \ $a \otimes \id_V$ \ and \ $b\otimes \id_V$.
\end{defn}

\parag{remarks}
\begin{enumerate}
\item Of course the action of \ $\mathcal{O}_X$ \ and of \ $\hat{A}$ \ commutes on \ $\Xi_{\lambda, X}^{(N)}\otimes V$.
\item For each point in \ $X$ \ we have an evaluation map
$$ ev_x : \Xi_{\lambda, X}^{(N)}\otimes V \to \Xi_{\lambda}^{(N)}\otimes V$$
which is \ $\hat{A}-$linear and surjective. For any subsheaf \ $\mathbb{E}$ \ of \ $ \Xi_{\lambda, X}^{(N)}\otimes V $ \ we shall denote by \ $\mathbb{E}(x)$ \ its image by this map. \\
Note that \ $\mathbb{E}(x)$ \ is the image of the fiber at \ $x$ \ of \ $\mathbb{E}$ \ in the fiber at \ $x$ \ of \ $\Xi_{\lambda, X}^{(N)}\otimes V$. So it is a quotient of its fiber at \ $x$ \ and the sheaf \ $\mathbb{E}$ \ may be free of finite type on \ $\mathcal{O}_X[[b]]$ \ and this quotient may not be an isomorphism.
\end{enumerate}

Because of this last remark, the following easy lemma will be required later on.

\begin{lemma}\label{easy coh.}
Let \ $X$ \ be a reduced complex space and \ $\mathcal{F} \subset \mathcal{O}_X^p$ \ be a coherent subsheaf such that for each \ $x$ \ the evaluation map \ $ev_x : \mathcal{F} \to \C^p \simeq \mathcal{O}X^p\big/\frak{M}_x.\mathcal{O}X^p$ \ has rank r. Then \ $\mathcal{F}$ \ is locally free rank \ $r$.
\end{lemma}

\parag{proof} Locally on \ $X$ \ consider  \ $M : \mathcal{O}_X^q \to \mathcal{O}_X^p $ \  an \ $\mathcal{O}_X-$map which image \ $\mathcal{F}$. Considering \ $M$ \ as an holomorphic  map to \ $End(\C^q,\C^p)$ \ we see that our hypothesis implies that the rank of \ $M(x)$ \ is constant and equal to \ $r$. So the image of \ $M$ \ is a rank \ $r$ \ vector bundle and \ $\mathcal{F}$, which is the sheaf of holomorphic sections of this rank \ $r$ \  vector bundle is then locally free of rank \ $r$. $\hfill \blacksquare$\\

\begin{defn}\label{Hol. 2}
A section \ $\varphi$ \ of the sheaf \ $\Xi_{\lambda, X}^{(N)}\otimes V$ \ will be called {\bf \ $k-$admissible} when the sub$-\hat{A}_X-$sheaf generated by \ $\varphi$ \ in \ $\Xi_{\lambda, X}^{(N)}\otimes V$ \ is a free rank \ $k$ \ $\mathcal{O}_X[[b]]-$module with basis \ $\varphi, a.\varphi, \dots, a^{k-1}.\varphi$ \ such that, for each \ $x \in X$,  $\mathbb{E}(x)$ \ is a rank \ $k$ \ fresco, where we denote  \ $\mathbb{E}(x)$ \ the image of \ $\mathbb{E}$ \ by the evalution map \ $ev_x$.\\
A \ $\hat{A}_X-$sheaf \ $\mathbb{E}$ \ on \ $X$ \ will be called an {\bf holomorphic family of rank \ $k$ \ $[\lambda]-$primitive frescos} parametrized by \ $X$, when the sheaf \ $\mathbb{E}$ \  is locally isomorphic on \ $X$ \  to a \ $\hat{A}_X-$sheaf generated by an admissible section of \ $ \Xi_{\lambda, X}^{(N)}\otimes V $, for some \ $N$ \ and some finite dimensional complex vector space \ $V$.\\
\end{defn}

Of course the sheaf \ $\mathbb{E}$ \ on \ $X$ \ defines the family \ $\mathbb{E}(x), x \in X$ \ of rank \ $k$ \ $[\lambda]-$primitive frescos. And conversely, such a family will be holomorphic when it is given by the evaluations maps \ $ev_x$ \ from a \ $\hat{A}_X-$sheaf \ $\mathbb{E}$ \ as in the previous definition.

\parag{remarks}
\begin{enumerate}[i)]
\item The sub$-\hat{A}-$sheaf generated by a \ $k-$admissible section \ $\varphi$ \ of \ $\Xi_{\lambda, X}^{(N)}\otimes V$ \ is stable by \ $\C[[a]]$ \ by definition. So we may write locally
$$ a^k.\varphi = \sum_{j=0}^{k-1} \ T_j.a^j.\varphi $$
where \ $T_j$ \ are unique (local) section of \ $\mathcal{O}_X[[b]]$. So there is unique monic degree \ $k$ \ polynomial \ $P : = a^k - \sum_{j=0}^{k-1} \ T_j.a^j$ \ in \ $\mathcal{O}_X[[b]][a]$ \ such that \ $P.\varphi = 0$. An easy argument of division shows that the left ideal in the sub-algebra \ $\mathcal{O}_X[[b]][a]$ \ of \ $\hat{A}_X$ \ generated by \ $P$ \ is the annihilator of \ $\varphi$.
\item As long as we work inside a simple pole \ $\hat{A}_X-$sheaf (i.e. \ $a.\mathcal{F}\subset b.\mathcal{F}$ \ which is a locally free finite type \ $\mathcal{O}_X[[b]]$ \ module (as \ $\Xi_{\lambda,X}^{(N)}$), we dont need to consider the sheaf algebra \ $\hat{A}_X$ \ to define the \ $\hat{A}_X-$structure :  the \ $\hat{A}_X-$structure is completely defined by the \ $\mathcal{O}_X[[b]]-$structure and the action of \ $a$ \ (i.e. the left module structure on the sub-algebra \ $\mathcal{O}_X[[b]][a]$) because the \ $\mathcal{O}_X[[b]]-$completion implies the \ $\mathcal{O}_X[[a]]$ \ completion.\\
\item Our definition of an \ $1-$admissible section  is not satisfied by the following :
$$\varphi : =  z.s^{\lambda_1-1} \in \Xi_{\lambda,\C}^{(0)}$$
where \ $z$ \ is the coordinate on \ $X : = \C$ \ and \ $\lambda_1 \in \mathbb{Q}^{+*}$ \ is in \ $[\lambda]$ \ because the evaluation at \ $z = 0$ \ of the \ $\hat{A}_{\C}-$sheaf generated by \ $\varphi$ \ is \ $0$. But of course the "abstract" sheaf \ $\hat{A}_{\C}.\varphi$ \ is isomorphic to the subsheaf of \ $\Xi_{\lambda,\C}^{(0)}$ \ generated by the \ $1-$admissible section \ $\psi : = s^{\lambda_1-1}$ \ and so it defines an holomorphic family of rank \ $1$ \ frescos parametrized by \ $\C$.\\
\end{enumerate}

\begin{lemma}\label{Hol. 3}
Let \ $\mathbb{E}$ \ be an holomorphic family of rank \ $k$ \ $[\lambda]-$primitive frescos. Then for each \ $x$ \ in \ $X$ \ the \ $\hat{A}-$module \ $\mathbb{E}(x)$ \ is a rank \ $k$ \ $[\lambda]-$primitive fresco and  the Bernstein polynomial of \ $\mathbb{E}(x)$ \ is locally constant on \ $X$.
\end{lemma}

\parag{Proof} It is enough to prove the result for the \ $\hat{A}_X-$sheaf generated by an admissible section \ $\varphi$ \ of the sheaf \ $\Xi_{\lambda, X}^{(N)}\otimes V$. Then write
\begin{equation*}
 a^k.\varphi = \sum_{j=1}^k \ S_j(b,x).a^{k-j}.\varphi  \tag{@}
 \end{equation*}
 
where  \ $S_1, \dots, S_k$ \ are sections of the sheaf \ $\mathcal{O}_X[[b]]$. For each given \ $x$, we know that \ $\mathbb{E}(x)$ \ is a fresco of rank \ $ k$. This implies that the valuation in \ $b$ \ of \ $S_j$ \ is \ $\geq j$. Then we have \ $S_j(b,x) = b^j.s_j(x) + b^{j+1}.\mathcal{O}_X[[b]] $ \ where \ $s_j$ \ is an holomorphic function on \ $X$. As the Bernstein element\footnote{the Bernstein element  of a fresco \ $E$ \  with fundamental invariants \ $\lambda_1, \dots, \lambda_k$ \ is the element \ $P_E : = (a - \lambda_1.b).\dots (a - \lambda_k.b)$ \  in \ $\hat{A}$ ; see [B.09-a].} of \ $\mathbb{E}(x)$ \ is equal to \ $a^k + \sum_{j=1}^k \ s_j(x).b^j.a^{k-j} $ \ and has rational coefficients, we conclude that \ $s_j$ \ is locally constant on \ $X$. $\hfill \blacksquare$\\

\begin{prop}\label{chgt. var. hol.}
Let \ $\mathbb{E}$ \ be an holomorphic family of rank \ $k$ \ $[\lambda]-$primitive frescos parametrized by a reduced  complex space \ $X$\  and let \ $\theta \in \C[[a]]$ \ be a change of variable. Then the family \ $\theta_*(\mathbb{E}(x)), x \in X$ \ is holomorphic.
\end{prop}

\parag{Proof} The problem is local on \ $X$ \ and we may assume that \ $\mathbb{E}$ \ is generated by a \ $k-$admissible section \ $\varphi$ \ of some sheaf \ $\Xi_{\lambda, X}^{(N)}\otimes V $. As the \ $\hat{A}-$module \ $\theta_*(\Xi_{\lambda, X}^{(N)}\otimes V)$ \ is isomorphic to \ $\Xi_{\lambda, X}^{(N)}\otimes V $, thanks to the  propsition \ref{chang. dev.}, it is enough to prove that \ $\theta_*(\varphi)$ \ is a \ $k-$admissible section generating \ $\theta_*(\mathbb{E})$. But this is obvious. $\hfill \blacksquare$\\

\begin{lemma}\label{easy}
Let \ $\varphi$ \ be a \ $k-$admissible section of some sheaf \ $ \Xi_{\lambda,X}^{(N)} \otimes V$ \ and \ $P : = (a - \lambda_1.b).S_1^{-1} \dots S_{k-1}^{-1}.(a - \lambda_k.b) $ \ a generator of the annihilator of \ $\varphi$, not necessarily corresponding to the principal J-H. sequence. Then \ $(a - \lambda_k.b).\varphi$ \ is \ $(k-1)-$admissible and define a normal corank \ $1$ \ holomorphic sub-family of the holomorphic family defined by \ $\mathbb{E}$.
\end{lemma}

\parag{proof} As \ $(a - \lambda_k.b).\varphi(x)$ \ generates  for each \ $x \in X$ \ a rank \ $(k-1)-$fresco which is normal in \ $\A.\varphi(x)$, the assertion is obvious. $\hfill \blacksquare$\\

\begin{thm}\label{J-H. hol.}
Let \ $X$ \ be a reduced complex space and \ $\mathbb{E}$ \ an holomorphic family of rank \ $k$ \ $[\lambda]-$primitive frescos parametrized by \ $X$. Then for each \ $j$ \ in \ $[1,k]$ \ there exists a sub$-\hat{A}_X-$sheaf \ $\mathbb{F}_j \subset \mathbb{E}$ \ which is an holomorphic family of rank \ $j$ \ $[\lambda]-$primitive frescos parametrized by \ $X$ \ such that for each \ $x$ \ in \ $X$ \ the submodules \ $\mathbb{F}_j(x), j \in [1,k]$ \ give the principal J-H. sequence of the fresco \ $\mathbb{E}(x)$.\\
Moreover, for each \ $j \in [1,k]$ \ the quotient sheaf \ $\mathbb{E}\big/\mathbb{F}_j$ \ is an holomorphic family of frescos.
\end{thm}

\parag{proof} In fact we shall prove by induction on \ $k \geq 1$ \ that for any \ $k-$admissible section \ $\varphi$ \ of the sheaf \ $\Xi_{\lambda, X}^{(N)}\otimes V$ \ corresponding to an holomorphic family of rank \ $k$ \ frescos with fundamental invariants \ $\lambda_1, \dots, \lambda_k$, there exists locally on \ $X$, invertible sections \ $S_1, \dots, S_k$ \ of the sheaf \ $\mathcal{O}_X[[b]]$ \ such that
\begin{equation*}
 (a - \lambda_1.b).S_1^{-1}.(a - \lambda_2.b) \dots S_{k-1}^{-1}.(a - \lambda_k.b).S_k^{-1}.\varphi \equiv 0. \tag{@@}
 \end{equation*}
Then it will not be difficult to prove that for each \ $j \in [0,k-1]$ \ the section
$$ \varphi_j : = (a - \lambda_{k-j+1}.b).S_{k-j+1}^{-1} \dots S_{k-1}^{-1}.(a - \lambda_k.b).S_k^{-1}.\varphi $$
is \ $(k-j)-$admissible and generates the sheaf \ $\mathbb{F}_{k-j}$ \ corresponding to the family \ $\mathbb{F}_{k-j}(x), x \in X$ \ locally.\\

Let us begin by the case \ $k = 1$. It is contained in the following lemma.

\begin{lemma}\label{Hol. rank 1}
Let \ $X$ \ be a connected reduced complex space and let  \ $\varphi$ \ be an $1-$admissible section of the sheaf \ $\Xi^{(N)}_{\lambda,X}\otimes V$. Then there exists locally on \ $X$ \ an invertible section \ $S$ \ of the sheaf \ $\mathcal{O}_X[[b]]$ \ and a non vanishing holomorphic function vector \ $v_0 : X \to  V$ \ such that 
$$ \varphi = S(b,x).s^{\lambda_1-1} \otimes v_0(x)$$
where \ $\lambda_1$ \ is a positive rational number in \ $[\lambda] \in \mathbb{Q}\big/\mathbb{Z}$.
\end{lemma}

\parag{proof} By definition the \ $\mathcal{O}_X[[b]]-$module \ $\mathbb{E}$ \  generated by \ $\varphi$ \ is free with rank \ $1$ \ and stable by \ $a$. Write \ $a.\varphi = T(b,x).\varphi$. As for each \ $x$ \ in \ $X$ \ the fresco \ $\mathbb{E}(x)$ \ is rank \ $1$, it has a simple pole and \ $T(b,x)$ \ has valuation \ $\geq 1$ \ in \ $b$. Moreover if we put \ $T(b,x) = t(x).b + b^2.T_1(b,x)$ \ the Bernstein element  of \ $\mathbb{E}(x)$ \ is \ $a - t(x).b$. So, as in the lemma \ref{Hol. 3} this implies that the holomorphic function \ $t$ \ is constant on \ $X$ ; let \ $\lambda_1$ \ the positive  rational value of \ $t$. We may write 
 $$a.\varphi =  b.\big(\lambda_1 +  b.T_1(b,x)\big).\varphi .$$
 Now the section \ $U : = \lambda_1 +  b.T_1(b,x)$ \ of \ $\mathcal{O}_X[[b]]$ \ is invertible and  we may solve the differential equation (with holomorphic parameter in \ $X$)
 $$b.S' + S.(U - \lambda_1) = 0  \quad S(0,x) \equiv 1 $$
 where \ $S$ \ is a section of \ $\mathcal{O}_X[[b]]$, and obtain the new generator \ $\psi : = S.\varphi$ \ for \ $\mathbb{E}$ \ with the following property
 \begin{align*}
 & a.S.\varphi = S.a.\varphi + b^2.S'.\varphi = S.b.U.\varphi + b^2.S'.\varphi =   \lambda_1.b.S.\varphi \quad {\rm and \  so} \\
 &a.\psi = \lambda_1.b.\psi 
  \end{align*}
  and this implies that there exists an holomorphic function \ $v_0 : X \to   V $ \ such that \ $\psi = s^{\lambda_1-1}\otimes v_0$. To conclude the proof it is enough to remark that if \ $v_0(x) = 0$ \ then \ $\mathbb{E}(x) = \{0\}$ \ contradicting our rank \ $1$ \ assumption. So \ $v_0$ \ does not vanish. $\hfill \blacksquare$\\
  
  So now assume that we have proved our assertion \ $(@@)$ \ for \ $k-1 \geq 1$, and consider now an \ $k-$admissible section \ $\varphi$ \ of the sheaf \ $ \Xi_{\lambda, X}^{(N)}\otimes V$ \ corresponding to an holomorphic family of rank \ $k$ \ frescos with fundamental invariants \ $\lambda_1, \dots, \lambda_k$. Consider the subsheaf \ $\mathbb{F} : = \big(\mathcal{O}_X.s^{\lambda_1-1}\otimes V\big) \cap \mathbb{E}$ \ of \ $\mathbb{E}$. It is \ $\mathcal{O}_X-$coherent as it is contained in the \ $\mathcal{O}_X-$coherent subsheaf  of \ $\Xi^{(0)}_{\lambda,X} \otimes V$ \ where the degree in \ $b$ \ is bounded by some  \ $n \gg 1$. As \ $\mathbb{F}(x) = F_1(\mathbb{E}(x)$ \ has dimension 1 for each \ $x$ \ in \ $X$, it is a locally free rank \ $1$ \ $\mathcal{O}_X-$module, thanks to lemma \ref{easy coh.} ; then we may find, at least locally, an holomorphic section \ $v_0 : X \to V$ \ which is a local basis for it, so \ $v_0$ \ does not vanish. The subsheaf \ $\mathcal{O}_X[[b]].\psi$ \ where \ $\psi : = s^{\lambda_1-1}\otimes v_0$ \ has a  simple pole \ because $a.\psi = \lambda_1.\psi$. Then \ $\psi$ \ is admissible an define the holomorphic family \ $\mathbb{F}_1(x)_{x \in X}$. \\
  
To prove the second assertion of the theorem, remark first that if \ $\chi$ \ is a section of \ $\mathbb{E}$ \ of the form \ $\chi : = S(b,x).s^{\mu-1}\otimes v_0 $ \ where we assume that the generic valuation of \ $S$ \ in \ $b$ \ is \ $0$ ;  then \ $\chi$ \ is a section of \ $\mathcal{O}_X[[b]].\psi$ \ because  we have \ $\mu  \geq \lambda_1$ : indeed,  near a generic point \ $x$ \ the section \ $S$ \ will be invertible and we find \ $s^{\mu-1}\otimes v_0$ \ in \ $\mathbb{E}$. But for each \ $x$ \ in \ $X$ \ there is no non zero element in \ $\mathbb{E}(x)$ \ which is annihilated by\ $a - \mu.b$ \ with \ $\mu < \lambda_1$. This proves our claim.\\
Now write locally  \ $\mathcal{O}_X\otimes V = (\mathcal{O}_X \otimes W) \oplus \mathcal{O}_X.v_0 $ \ where \ $W$ \ is an hyperplane in \ $V$, and consider the map
$$ \sigma : \Xi_{\lambda,X}^{(N)} \otimes V \to  \Xi_{\lambda,X}^{(N)} \otimes V $$
given by the identity map on \ $ \Xi_{\lambda,X}^{(N)} \otimes W$ \ and given on \ $\Xi_{\lambda,X}^{(N)} \otimes v_0$ \ by the quotient map by \ $ \Xi_{\lambda,X}^{(0)} \otimes v_0$ :
$$ \Xi_{\lambda,X}^{(N)} \otimes v_0 \to  \Xi_{\lambda,X}^{(N-1)} \otimes v_0$$
composed with the obvious inclusion 
$$  \Xi_{\lambda,X}^{(N-1)} \otimes v_0 \to  \Xi_{\lambda,X}^{(N)} \otimes v_0.$$
As we  proved above  that the intersection of the kernel of \ $\sigma$ \ with \ $\mathbb{E}$ \ is \ $\mathbb{F}_1$ \  the quotient sheaf \ $\mathbb{E}\big/\mathbb{F}_1$ \ is isomorphic to the subsheaf of \ $\Xi_{\lambda,X}^{(N)} \otimes V$ \ generated by the section \ $ \sigma(\varphi)$.\\
Moreover we have\ $\hat{A}_X.\sigma(\varphi)(x) \simeq \mathbb{E}(x)\big/\mathbb{F}_1(x)$. We shall prove that the section  \ $\sigma(\varphi)$ \ is \ $(k-1)-$admissible. First we shall prove 
 that \ $\sigma(\varphi), a.\sigma(\varphi), \dots, a^{k-2}.\sigma(\varphi)$ \ is locally free on \ $\mathcal{O}_X[[b]]$ \ in the sheaf \ $ \Xi_{\lambda,X}^{(N)} \otimes V $. This is equivalent to show that the corresponding classes are \ $\mathcal{O}_X-$free modulo \ $b$. We shall use the following lemma.

\begin{lemma}\label{tech}
In the situation above write
$$ \psi = s^{\lambda_1-1} \otimes v_0 = \sum_{j=0}^{k-1} \ U_j.a^j.\varphi $$
then the \ $b-$valuation of  \ $U_j $ \ is at least \ $b^{k-j-1}$, for \ $j \in [0,k-1]$.
\end{lemma}

\parag{proof} Write \ $(a- \lambda_1.b).\psi = 0$. It gives the equations
\begin{align*}
& U_{k-1}.S_j + b^2.U'_j - \lambda_1.b.U_j + U_{j-1} = 0 
\end{align*}
where we have written  \ $a^k.\varphi = \sum_{j=0}^{k-1} S_j.a^j.\varphi$ \ with the convention \ $U_{-1} \equiv 0 $.\\
Recall that the \ $b-$valuation of \ $S_j$ \ is at least \ $k-j$. Then an easy induction conclude the proof.$\hfill \blacksquare$

\parag{End of the proof of the theorem \ref{J-H. hol.}} Now remember that \ $\psi$ \ does not belongs to \ $b.\mathbb{E}$ \ because \ $\mathbb{F}_1(x)$ \ is normal in \ $\mathbb{E}(x)$ \ for each \ $x$. So \ $U_{k-1}$ \ is invertible in \ $\mathcal{O}_X[[b]]$ \ and we have
$$ U_{k-1}.a^{k-1}.\sigma(\varphi) = - \sum_{j=0}^{k-2} \ U_j.a^j.\sigma(\varphi).$$

 \ $\psi = U_{k-1}.a^{k-1}.\varphi + b.\mathbb{E}$. So we easily deduce that \ $\sigma(\varphi), a.\sigma(\varphi), \dots, a^{k-2}.\sigma(\varphi)$ \ is locally free on \ $\mathcal{O}_X[[b]]$ \ in the sheaf \ $ \Xi_{\lambda,X}^{(N)} \otimes V $ \ and generate a sub-sheaf stable by \ $a$. So \ $\sigma(\varphi)$ \ is \ $(k-1)-$admissible with corresponding fundamental invariants \ $\lambda_2, \dots, \lambda_k$ \ and the induction hypothesis gives \ $S_2, \dots, S_k$ \ invertible in \ $\mathcal{O}_X[[b]]$ \ such that
$$ (a - \lambda_2.b).S_2^{-1} \dots S_{k-1}^{-1}.(a - \lambda_k.b).S_k^{-1}.\sigma(\varphi) \equiv 0 $$
in \ $\mathbb{E}\big/\mathbb{F}_1$. So we get
$$ (a - \lambda_2.b).S_2^{-1} \dots S_{k-1}^{-1}.(a - \lambda_k.b).S_k^{-1}.\varphi = T.\psi $$ with \ $T \in \mathcal{O}_X[[b]].\psi .$
This equality modulo \ $b$ \ gives an invertible element \ $f$ \  of \ $\mathcal{O}_X$ \ such that \ $f.a^{k-1}.\varphi = T(0).\psi \quad modulo \ b\mathbb{E} $. \\
But we already know that 
 $\psi = U_{k-1}.a^{k-1}.\varphi + b.\mathbb{E}$, with \ $U_{k-1}(0)$ \ invertible in \ $\mathcal{O}_X$, so we conclude that we have 
$$ (f - T(0).U_{k-1}(0))a^{k-1}.\varphi \in b.\mathbb{E} .$$
As we know that for each \ $x \in X$ \ we have \ $a^{k-1}.\varphi(x) \not\in b.\mathbb{E}(x)$ \ it implies  that \ $f = U_{k-1}(0).T(0)$. So \ $T$ \ is an invertible element in \ $\mathcal{O}_X[[b]]$ \ and we may define \ $S_1 : = T$ \ to conclude our induction.\\

Now to complete the proof of the theorem it is enough to prove that for each \\
 $j \in [1,k-1]$ \ the section   \ $\varphi_j$ \ is \ $(k-j)-$admissible. But it is clear that, modulo \ $b.\mathbb{E}$, the classes of \ $\varphi_j, a.\varphi_j, \dots, a^{k-j-1}.\varphi_j$ \ co{\"i}ncide with \ $a^j.\varphi, \dots, a^{k-1}.\varphi$. So they are independant on \ $\mathcal{O}_X[[b]]$. The identity
$$ (a - \lambda_1.b).S_1^{-1} \dots (a - \lambda_{k-j}.b)S_{k-j}^{-1}.\varphi_j \equiv 0 $$
shows that \ $a^{k-j}.\varphi_j$ \ is in the \ $\mathcal{O}_X[[b]]-$module with basis \ $\varphi_j, a.\varphi_j, \dots, a^{k-j-1}.\varphi_j$, and then \ $\varphi_j$ \ is \ $(k-j)-$admissible for each \ $j$.$\hfill \blacksquare$\\

\subsection{Basic theorems on holomorphic families.}

\begin{thm}\label{suite exacte}
Let \ $X$ \ be a reduced complex space and consider an exact sequence of \ $\hat{A}_X-$sheaves :
$$ 0 \to \mathbb{F} \to \mathbb{E} \to \mathbb{G} \to 0 .$$
Assume that \ $ \mathbb{F}$ \ and \ $ \mathbb{G}$ \ are homorphic families of \ $[\lambda]-$primitive frescos parametrized by \ $X$.
 Assume also  that for each \ $x \in X$ \ we have an exact sequence of frescos given by the fibers at \ $x$
$$0 \to \mathbb{F}(x) \to \mathbb{E}(x) \to \mathbb{G}(x) \to 0 $$for each \ $x \in X$, $\mathbb{E}(x)$ \ is a fresco. Then the family \ $\mathbb{E}$ \ is holomorphic.
\end{thm}

\parag{Proof} We shall show that it is enough to prove the case where \ $\mathbb{F}$ \ has rank \ $1$.
Assume that \ $\mathbb{F}$ \ has rank \ $k-1 \geq 2$ \ and that the theorem has been proved for rank \ $\leq k-2$ \ in the first sheaf of the exact sequence. Let \ $\mathbb{F}_1$ \ be the sheaf given by the rank \ $1$ \ term in the principal J-H. sequence of \ $\mathbb{F}$. Then we have an exact sequence of \ $\hat{A}_X-$sheaves :
$$  0 \to \mathbb{F}\big/\mathbb{F}_1 \to \mathbb{E}\big/\mathbb{F}_1 \to \mathbb{G} \to 0 $$
and now the first sheaf has rank \ $k-2$ \ and is an holomorphic family thanks to \ref{J-H. hol.}. So the induction hypothesis implies the holomorphy of the family \ $ \mathbb{E}\big/\mathbb{F}_1$. So applying the rank \ $1$ \ case to the exact sequence of \ $\hat{A}_X-$sheaves :
$$ 0 \to \mathbb{F}_1 \to \mathbb{E} \to \mathbb{E}\big/\mathbb{F}_1 \to 0 $$
we conclude that \ $\mathbb{E}$ \ is an holomorphic family.\\

Now we consider the case where \ $\mathbb{F}$ \ is rank \ $1$.  Let \ $k$ \ be the rank of \ $\mathbb{E}$.The assertion is local on \ $X$ \ so we may assume that \ $\mathbb{G}$ \ is generated by a \ $(k-1)-$admissible section \ $\psi$ \ of some sheaf \ $\Xi_{\lambda,X}^{(N)} \otimes V$. We may find (locally on \ $X$) \ a section \ $\varphi$ \ of \ $\mathbb{E}$ \ which lift \ $\psi$. Let \ $Q$ \ be a section of the sheaf \ $\sum_{j=0}^{k-1} \ \mathcal{O}_X[[b]].a^j $, monic of degree \ $k-1$ \ in \ $a$ \ such that \ $\mathbb{G} \simeq \hat{A}_X\big/\hat{A}_X.Q $. So \ $Q$ \ generates the annihilator\footnote{ see the remark i) following the definition \ref{Hol. 2}.} of \ $\psi$ \ in \ $\hat{A}_X$.
Then \ $Q.\varphi$ \ is a section of \ $\mathbb{F}$ \ which generates \ $\mathbb{F}(x)$ \ for each \ $x \in X$, and as we may assume at least locally on \ $X$, that \ $\mathbb{F} \simeq \mathcal{O}_X[[b]].s^{\lambda-1}$, there exists an invertible section \ $S_1$ \ of the sheaf \ $\mathcal{O}_X[[b]]$ \ such that \ $P : = (a - \lambda_1.b).S_1^{-1}.Q$ \ annihilates \ $\varphi$. So we may (locally) assume that \ $\varphi$ \ is a section of some sheaf \ $\Xi_{\lambda,X}^{(N')}\otimes V'$ \ such that \ $\mathbb{E}(x) \simeq \hat{A}.\varphi(x)$.\\
To conclude, it is then enough to show that \ $\varphi, a.\varphi, \dots, a^{k-1}.\varphi$ \ is a \ $\mathcal{O}_X[[b]]-$basis of \ $\mathbb{E}$. As this is true for each given \ $x$ \ in \ $X$, the only point to prove is the fact that they generate the sheaf \ $\mathbb{E}$. Let \ $\sigma$ \ be a local section of  \ $\mathbb{E}$. We may write, locally on \ $X$, the image of \ $\sigma$ \ by the quotient map \ $\pi :  \mathbb{E} \to \mathbb{G}$ : 
$$  \pi(\sigma) = \sum_{j=0}^{k-2} \ T_j.a^j.\psi $$
where \ $T_j$ \ are local sections of \ $\mathcal{O}_X[[b]]$. Then \ $\sigma - \sum_{j=0}^{k-2} \ T_j.a^j.\varphi$ \ is a section of \ $\mathbb{F}$ \ and it may be written as \ $T.Q.\varphi$. This gives the conclusion. $\hfill \blacksquare$\\

\begin{thm}\label{hol. dual.}
Let \ $X$ \ be a reduced complex space and \ $\mathbb{E}$ \ be an holomorphic family of rank \ $k$ \ $[\lambda]-$primitive frescos parametrized by \ $X$. Assume that for some \ $\delta \in \mathbb{Q}$ \ and for each \ $x \in X$ \ the \ $\hat{A}-$module \ $\mathbb{E}(x)^*\otimes E_{\delta}$ \ is a fresco. Then the family \ $\mathbb{E}^*\otimes E_{\delta}$ \ is holomorphic.
\end{thm}

\parag{Remark} Assuming \ $X$ connected, let  $k$ \ be the rank of \ $ \mathbb{E}$ \ and  \ $\lambda_1, \dots, \lambda_k$ \ the fundamental invariants  of \ $\mathbb{E}$ ; then the inequality \ $\delta > \lambda_k + k-1$ \ implies that \ $\mathbb{E}(x)^*\otimes E_{\delta}$ \ is a fresco for each \ $x \in X$.

\parag{Proof} We shall make an induction on the rank of the frescos. Remark that we may assume \ $X$ \ connected as the problem is clearly local on \ $X$, so the rank is well define.\\
In rank \ $1$ \ we know that the family may be locally given, assuming that the Bernstein element is \ $a - \lambda_1.b$ \ with \ $\lambda_1 \in \mathbb{Q}^{+*}$,  by the sheaf \ $\mathcal{O}_X[[b]].s^{\lambda_1-1} \subset \Xi_{\lambda,X}^{(0)} $ \ and the (twisted) dual is then defined (locally) by the sheaf \ $\mathcal{O}_X[[b]].s^{\delta - \lambda_1 -1} $ \ where the condition on \ $\delta \in \mathbb{Q}$ \ is that \ $\delta - \lambda_1 > 0 $. So the rank \ $1$ \ is clear.\\
Assume now that for \ $k-1 \geq 1$ \ the theorem is proved and consider  an holomorphic family \ $\mathbb{E}$ \  of rank \ $k$ \ $[\lambda]-$primitive frescos parametrized by \ $X$. We have a exact sequence, thanks to \ref{J-H. hol.},
$$ 0 \to \mathbb{F}_{k-1} \to \mathbb{E} \to \mathbb{E}\big/\mathbb{F}_{k-1} \to 0 $$
where \ $\mathbb{F}_{k-1}(x)$ \ is the family of the \ $(k-1)$ \ terms of the principal J-H. sequence of \ $\mathbb{E}(x)$ \ for each \ $x \in X$. Then we have an exact sequence of sheaves
$$ 0 \to ( \mathbb{E}\big/\mathbb{F}_{k-1})^* \to \mathbb{E}^* \to \mathbb{F}_{k-1} ^* \to 0  $$
and we know that the families \ $( \mathbb{E}\big/\mathbb{F}_{k-1})^* \to \mathbb{E}^*$ \ and \ $\mathbb{F}_{k-1} ^* $ \ are holomorphic, using the induction hypothesis. Then the theorem \ref{suite exacte} gives the holomorphy of \ $ \mathbb{E}^*$. $\hfill \blacksquare$\\

Combining these two theorems we obtain the following result.

\begin{thm}\label{deux sur trois}
Let \ $X$ \ be a reduced complex space and consider an exact sequence of \ $\hat{A}_X-$sheaves :
$$ 0 \to \mathbb{F} \to \mathbb{E} \to \mathbb{G} \to 0 .$$
Assume that two of these are homorphic families of \ $[\lambda]-$primitive frescos parametrized by \ $X$ \ and  that for each \ $x \in X$ \ we have an exact sequence of frescos given by the fibers at \ $x$
$$0 \to \mathbb{F}(x) \to \mathbb{E}(x) \to \mathbb{G}(x) \to 0 $$
 then the third  family is holomorphic.
\end{thm}

\parag{Proof} The  statement in the case where \ $\mathbb{F}$ \ and \ $\mathbb{E}$ \ are holomorphic is easily reduced to  the rank \ $1$ \ case for \ $\mathbb{F}$ : if \ $\mathbb{F}$ \ has rank \ $r \geq 2$ \ let \ $\mathbb{F}_1$ \ be the first term is the principal J-H. sequence for \ $\mathbb{F}$. Then we have  the exact sequence
$$ 0 \to \mathbb{F}\big/\mathbb{F}_1 \to \mathbb{E}\big/\mathbb{F}_1 \to \mathbb{G} \to 0 $$
where the first two sheaves are holomorphic families, thanks to the rank \ $1$ \ case, as we know that \ $\mathbb{F}_1$ \ is an holomorphic applying the theorem \ref{J-H. hol.}.\\
So we have to prove the first case in the rank \ $1$ \ case for \ $\mathbb{F}$. But, as the statement in local on \ $X$ \ we may assume that \ $\mathbb{E}$ \ is generated by a \ $k-$admissible section \ $\varphi$ \ of some sheaf \ $\Xi_{\lambda,X}^{(N)}\otimes V$. Then  \ $\mathbb{F}$ \ is locally generated by a section \ $\psi$ \ of \ $ \Xi_{\lambda,X}^{(N)}\otimes V$ \ satisfying \ $(a - \mu.b).\psi = 0$ \ where \ $-\mu$ \ is the root of the Bernstein polynomial of \ $\mathbb{F}$. If we show that \ $\psi(x)$ \ never vanishes then  we may argue as in the proof of the theorem \ref{J-H. hol.} when we proved that \ $\mathbb{E}\big/\mathbb{F}_1$ \ is holomorphic. But if \ $\psi(x) = 0$ \ for some \ $x \in X$, its means that the map \ $\mathbb{F}(x) \to \mathbb{E}(x)$ \ vanishes\footnote{here \ $\psi(x)$ \ means the evaluation at \ $x$ \ of the section \ $\psi$ \ of \ $ \Xi_{\lambda,X}^{(N)}\otimes V$, which is the image of the class of \ $\psi$ \ in the fiber \ $F(x) : = \mathbb{F}\big/\frak{M}_x.\mathbb{F}$ \ by the inclusion map \ $i$ \ of \ $\mathbb{F}$ \ in \ $ \Xi_{\lambda,X}^{(N)}\otimes V$ ; but  the map \ $i$ \ after tensorization by \ $\mathcal{O}_X\big/\frak{M}_x$ \ may be not injective. See the example given in remark iii) following the definition \ref{Hol. 2}.} so we contradict our assumption
because our \ $\mathbb{F}(x)$ \ has rank \ $1$ \ for any \ $x$.\\
The second case  is consequence of the first one using duality and  the theorem \ref{hol. dual.}. \\
We proved  the third case in the theorem \ref{suite exacte}. $\hfill \blacksquare$

\subsection{Existence of versal families of \ $[\lambda]-$primitive frescos.}

Let me give an  easy corollary of the previous theorem, thanks to the theorem 3.4.1 of [B.09-b].

\begin{thm}\label{versal}
Let \ $\lambda_1, p_1, \dots, p_{k-1}$ \ be the fundamental invariants of a \ $[\lambda]-$primitive fresco. There exists a finite dimensional  affine complex manifold \ $F(\lambda_1, p_1, \dots, p_{k-1})$ \  and an holomorphic family of rank \ $k$ \ $[\lambda]-$primitive frescos \ $\mathbb{E}$ \ on it  which is locally versal near each point in \ $F(\lambda_1, p_1, \dots, p_{k-1})$.
\end{thm}

\parag{Proof} The key point is the following easy generalization of the lemma 3.4.2 of [B.09-a] :

\begin{lemma}\label{polyn. inv.}
Let \ $X$ \ be a reduced complex space and let \ $\mathbb{E}$ \ be the sheaf generated by a $k-$admissible section \ $\varphi$ \ of the sheaf \ $\Xi_{\lambda,X}^{(N)} \otimes V$. Assume that we have an \ $\hat{A}_X-$surjective map \ $\pi : \mathbb{E} \to \mathbb{E}_{\mu}$ \ where \ $\mathbb{E}_{\mu}$ \ is the constant family parametrized by \ $X$ \ with value \ $E_{\mu}$ \ where \ $\mu = \lambda+q$ \ for some \ $q \in \mathbb{N}$. Let \ $T_1, \dots, T_k$ \ be invertible sections of \ $\mathcal{O}_X[[b]]$, and let 
$$ \pi' : (a - \lambda_1.b).T_1^{-1}.(a - \lambda_2).T_2^{-1} \dots (a - \lambda_k.b).T_k^{-1}\big[\mathbb{E}\big] \to \mathbb{E}_{\mu} $$
the map induced by \ $\pi$, where \ $\lambda_1, \dots, \lambda_k$ \ are in \ $\lambda+ \mathbb{N}$.\\
Then the image of \ $\pi'$ \ contains \ $b^{k+m+1}.\mathbb{E}_{\mu}$ \ for \ $m$ \ large enough. In fact if we put \ $\mu = \lambda_j - m_j$ \ any non negative integer \ $m \geq \sup \{m_j, j\in [1,k] \}$ \ will be convenient.
\end{lemma}

As they are very similar to the proof of the same results for \ $[\lambda]-$primitive themes given in [B.09-b], we let the proof of this lemma and the proof of the corollary above as exercices for the reader.\\

\begin{thm}\label{Themes et ss}
Let \ $\mathbb{E}$ \ be an holomorphic family of \ $[\lambda]-$primitive frescos parametrized by a reduced complex space \ $X$. For each integer \ $q \geq 1$ \  the set 
 $$Y_q(\mathbb{E}) : = \{ x \in X \ / \ d(\mathbb{E}(x)) \leq q \},$$
  where \ $d(E)$ \ denotes the ss-depth of the fresco \ $E$, is a closed analytic subset in \ $X$.
\end{thm}

\parag{remark} Of course, we have the following special cases :\\
   Let \ $T(\mathbb{E}) : = \{x \in X\ / \ \mathbb{E}(x) \ {\rm is \ a \  theme} \}$ \ and \ $SS(\mathbb{E}) : = \{x \in X \ / \ \mathbb{E}(x) \ {\rm is \ semi-simple} \}$. Then \ $T(\mathbb{E})$ \ is a Zariski open set in \ $X$ \ and \ $SS(\mathbb{E})$ \ is a closed analytic subset in \ $X$.\\

The proof of theorem is an easy consequence of the following  two lemmas. 

\begin{lemma}\label{ordre suff.}
Let \ $E$ \ be a \ $[\lambda]-$primitive fresco with fundamental invariants \ $\lambda_1, \cdots, \lambda_k$, and let \ $\mu \in [\lambda]$. Then for \ $M > \mu - \lambda_1$ \ the linear map 
$$ Ker \big[(a - \mu.b) : E \to E\big] \to  Ker \big[(a - \mu.b) : E\big/b^M.E  \to E\big/b^M.E\big] $$
is injective.
\end{lemma}

\parag{Proof} Note that \ $Ker \, (a - \mu.b)$ \ is the vector complex space \ $Hom_{\hat{A}}(E_{\mu},F)$ \ for any \ $\hat{A}-$module \ $F$. Now the kernel of the map
$$ \varphi : Hom_{\hat{A}}(E_{\mu},E) \to Hom_{\hat{A}}(E_{\mu}, E\big/b^M.E)$$
is the vector space \ $Hom_{\hat{A}}(E_{\mu}, b^M.E) \simeq Ker\,(a - \mu.b) \cap b^M.E $. But if we have \ $x = b^M.y$ \ such that \ $(a - \mu.b).x = 0$ \ then \ $(a - (\mu-M).b).y = 0$ \ and this implies \ $y = 0$ \ if \ $\mu - M < \lambda_1$. So \ $\varphi$ \ is injective for \ $M > \mu - \lambda_1$. $\hfill \blacksquare$\\

\begin{lemma}\label{coherence}
In the situation of the previous theorem, let \ $\mu \in [\lambda]$. Then the subsheaf \ $Ker\, (a - \mu.b) \subset \mathbb{E}$ \ is \ $\mathcal{O}_X-$coherent. More precisely, there exists locally on \ $X$ \  an integer \ $M \gg 1$ \ such that  \ $Ker\, (a - \mu.b) $ \ is isomorphic to the kernel of the \ $\mathcal{O}_X-$endomorphism induced by \ $(a - \mu.b)$ \ on the locally free finite type \ $\mathcal{O}_X-$module \ $\mathbb{E}\big/b^M.\mathbb{E}$.
\end{lemma}

\parag{Proof} First remark that in a \ $[\lambda]-$primitive  fresco \ $E$ \ with fundamental invariants \ $\lambda_1, \cdots, \lambda_k$ \ if we fix \ $\mu \in [\lambda]$, then, thanks to the previous lemma, there exists an integer \ $M$ \ such  that the projection on \ $E\big/b^M.E$ \ of \ $ Ker\, (a -\mu.b)$ \  is an isomorphism. So we have an isomorphism of the sheaf \ $Ker\, (a - \mu.b)$ \ onto  \ $Ker \big[(a -\mu.b) : \mathbb{E}\big/b^M.\mathbb{E} \to  \mathbb{E}\big/b^M.\mathbb{E}\big] $. But the sheaf \ $\mathbb{E}\big/b^M.\mathbb{E}$ \ is a locally free finite rank \ $\mathcal{O}_X-$module \ and the map \ $(a - \mu.b)$ \ is an \ $\mathcal{O}_X$ \ endomorphism of it. So its kernel is \ $\mathcal{O}_X-$coherent. $\hfill \blacksquare$

\parag{Proof of the theorem \ref{Themes et ss}} From proposition \ref{delta} we know that fixing \ $\mu \in [\lambda], \mu \geq \lambda_k+k-1$ \ we have for each \ $x \in X$
$$ \delta(\mathbb{E}(x)) = \dim_{\C}\big[ K(x)\big]  $$
where \ $K(x)$ \ is the image in \ $\mathbb{E}(x)$ \ of  \ $Ker\,(a- \mu.b)(x)$, the fiber at \ $x$ \ of the sheaf \ $Ker\, (a - \mu.b)$ \ which is coherent thanks to lemma \ref{coherence}. Then if we trivialize locally the holomorphic vector bundle \ $\mathbb{E}\big/b^M.\mathbb{E} \simeq X \times C^N$ \ and if \ $M : X \to End_{\C}(\C^N)$ \ is the holomorphic matrix of the endomorphism induced by \ $(a - \mu.b)$ \ in this (local) trivialization, then \ $K(x)$ \ is the kernel of \ $M(x)$. Then its is clear that the set where \ $\dim_{\C}[K(x)]$ \ is bigger than \ $q$ \ is a closed analytic set. Then the theorem is consequence of the formula \ $\delta(\mathbb{E}(x) = k - d(\mathbb{E}(x) + 1$, proved in the proposition \ref{delta}.$\hfill \blacksquare$ \\

\subsection{The theorem of the semi-simple part.}

We shall use the following definition.

\begin{defn}\label{ana.const.}
Let \ $X$ \ be a reduced complex space. We shall say that a closed set \ $Y$ \ in \ $X$ \ is {\bf analytically conctructible} ( {\bf cac} fort short) when there exists a finite collection of subsets \ $(Y_i)_{i \in [1,N]}$ \ in \ $X$ \ with the following properties :
\begin{enumerate}[i)]
\item For each \ $i \in [1,N]$ \ the union \ $\cup_{h \in [1,i-1]} Y_h $ \ is closed in \ $X$.
\item For each \ $i \in [1,N]$ \ \ $Y_{i}$ \ is a closed analytic subset in \ $X \setminus \cup_{h \in [1,i-1]} Y_h $.
\item  $Y = \cup_{i \in [1,N]} Y_i $
\end{enumerate}
\end{defn}

The main result of this section is the following theorem.

\begin{thm}\label{hol. S and Co}
Let \ $\mathbb{E}$ \ be an holomorphic family of  \ $[\lambda]-$primitive frescos parametrized by a reduced complex space \ $X$. Assume that at generic points in \ $X$ \ we have \ $d(\mathbb{E}(x) \geq 2$. Then there exists a closed analytically constructible set \ $Y \subset X$ \ with no interior points in \ $X$ \ such that on \ $X \setminus Y$ \ the families \ $S_1(\mathbb{E}(x))$ \ and \ $\Sigma^1(\mathbb{E}(x))$ \ are holomorphic.
\end{thm}

Of course using the results of section 3.2  we obtain immedialety the analoguous result for \ $S_j(\mathbb{E})$ \ and \ $\Sigma^j(\mathbb{E})$ \ for any \ $j \geq 2$.\\

We give in section 3.5  an example of an holomorphic family with constant ss-depth  parametrized by \ $\C$ \ such that the family \ $\Sigma^1(\mathbb{E}(x)) = L(\mathbb{E}(x))$ \ is not holomorphic at \ $0$. By twist duality we obtain also an example of an holomorphic family \ $\mathbb{E}'$ \ with constant ss-depth such that  the family \ $S_1(\mathbb{E}'(x))$ \ is not holomorphic at \ $0$.\\

 \begin{lemma}\label{induction key 1}
 Let  \ $\mathbb{E}$ \ be an holomorphic family of \ $[\lambda]-$primitives frescos parametrized by an  complex space \ $X$. Assume that we have \ $d(\mathbb{E}(x)) = 2$ \ for each \ $x \in X$ \ and that the holomorphic family \ $\mathbb{E}\big/\mathbb{F}_{k-2}$ \ is a family of themes. Then the family \ $L(\mathbb(E)(x)), x \in X$ \ is holomorphic and is locally isomorphic to the constant family equal to \ $E_{\lambda_{k-1}+k-2}$.
 \end{lemma}
 
 \parag{Proof} The assertion is local and we may assume that \ $\mathbb{E}$ \ is generated by a \ $k-$admissible section \ $\varphi $ \ of some sheaf \ $\Xi_{\lambda,X}^{(N)} \otimes V$. We then find sections \ $S_1, \dots, S_{k-1}$ \ of \ $\mathcal{O}_X[[b]]$ \ with \ $S_1(0) = \dots = S_{k-1}(0) = 1$ \ such that the annihilator ideal of \ $\varphi$ \ is generated by
 $$ P : = (a - \lambda_1.b).S_1^{-1} \dots S_{k-1}^{-1}.(a - \lambda_k.b) .$$
 Now our assumption on \ $\mathbb{E}\big/\mathbb{F}_{k-2}$ \ implies that for each \ $x \in X$ \ the coefficient \ $\alpha_{k-1}(x)$ \ of \ $b^{p_{k-1}}$ \ in \ $S_{k-1}(x)$ \ is not zero.  The fact that \ $d(\mathbb{E}(x)) = 2$ \ for each \ $x \in X$ \ implies that we have \ $\mathbb{F}_{k-1}(x) = S_1(\mathbb{E}(x))$ \ for each \ $x \in X$, thanks to lemma \ref{tech.}. So this implies that we have \ $\alpha_j(x) = 0$ \ for each \ $x \in X$ \ and each \ $j \in [1,k-2]$, where \ $\alpha_j(x)$ \ is the coefficient of \ $b^{p_j}$ \ in \ $S_j(x)$.
Then solving the differential equation\ $x \in X$  \ $b.U' = p_{k-2}.(U - S_{k-2}) $ \  with the holomorphic parameter with may write (using the commuting lemma of [B.09-a] lemma 3.5.1 or the lemma  \ref{com. lem.} at end of the section 2.4 )
$$ (a - \lambda_{k-2}.b).S_{k-2}^{-1}.(a - \lambda_{k-1}.b) = U^{-1}.(a - (\lambda_{k-1}+1).b).[S_{k-2}.U^{-2}]^{-1}.(a - (\lambda_{k-2}-1).b).U^{-1} $$
So writing 
 $$P = P_{k-3}.(a - \lambda_{k-2}.b).S_{k-2}^{-1}.(a - \lambda_{k-1}.b).S_{k-1}^{-1}.(a - \lambda_k.b) $$
 this gives
 $$ P = P_{k-3}.U^{-1}.(a - (\lambda_{k-1}+1).b).[S_{k-2}.U^{-2}]^{-1}.(a - (\lambda_{k-2}-1).b).[U.S_{k-1}]^{-1}.(a - \lambda_k.b).$$
Let \ $U_0$ \ the solution with no term in \ $b^{p_{k-2}}$, and let \ $\beta \in \mathcal{O}_X$ \ be the coefficient of \ $b^{p_{k-2} + p_{k-1}}$ \ in \ $U_0.S_{k-1}$. Choosing \ $U : = U_0 + \rho.b^{p_{k-2}}$ \ the coefficient of \ $b^{p_{k-2}}$ \ in \ $U.S_{k-1}$ \ will be \ $\beta + \rho.\alpha_{k-1}$. So the choice \ $\rho = - \beta.\alpha_{k-1}^{-1} \in \mathcal{O}_X$ \ allows to assume that the coefficient of \ $b^{p_{k-2} + p_{k-1}}$ \ in \ $U.S_{k-1}$ \ vanishes for all \ $x \in X$.\\
Then if \ $V \in \mathcal{O}_X[[b]]$ \ is a solution of \ $b.V' = (p_{k-2}+p_{k-1}).(V - U.S_{k-1}) $ \ we will obtain for \ $P$ \ the value
$$P_{k-3}.U^{-1}.(a - (\lambda_{k-1}+1).b).[S_{k-2}.U^{-2}.V]^{-1}.(a - (\lambda_{k}+1).b).[U.S_{k-1}.V^{-2}]^{-1}.(a - (\lambda_{k-2}-2).b).V^{-1} .  $$
The "standard computation" (see proposition  \ref{stand. comp.}) shows then that the coefficient of \ $b^{p_1}$ \ in   \ $S_{k-2}.U^{-2}.V$ \ is equal to \ $\frac{p_{k-1}+p_{k-2}}{p_{k-2}}.\alpha_{k-1}$. This means that for rank \ $k-1$ \  the holomorphic\footnote{The fact that \ $\chi$ \ is \ $(k-1)-$admissible follows from the lemma \ref{easy}.}family \ $\mathbb{G}$ \ associated to 
 $$ \chi : = [U.S_{k-1}.V^{-2}]^{-1}.(a - (\lambda_{k-2}-2).b).V^{-1}.\varphi$$
 satisfies the same conditions than the initial family \ $\mathbb{E}$ : the normal inclusion \ $\mathbb{G} \subset \mathbb{E}$ \ and the fact the coefficient of \ $b^{p_1}$ \ in  \ $S_{k-2}.U^{-2}.V$ \ implies that we have \ $d(\mathbb{G}(x)) = 2$ \ and \ $L(\mathbb{G}(x)) = L(\mathbb{E}(x))$ \ for any \ $x \in X$. So the only point to finish the proof is to begin the induction on the rank. But in rank \ $2$ \ we have a family of themes and then \ $L(E) = F_1(\mathbb{E})$ \ in this case. So the holomorphy of the family \ $L(\mathbb{E}(x))$ \ is a consequence of the theorem  \ref{J-H. hol.}. $\hfill \blacksquare$\\
 
 \begin{lemma}\label{induction key 2}
 We consider the same situation than in the previous lemma, but we assume that \ $d(\mathbb{F}_{k-1}) =1$ \ and that \ $\mathbb{E}\big/\mathbb{F}_{k-2}$ \ is a family of  semi-simple frescos. Then there exists locally on \ $X$ \ a normal  corank \ $1$ \ holomorphic sub-family \ $\mathbb{G} \subset \mathbb{E}$ \ such that \ $d(\mathbb{G}(x)) = 2$ \ for each \ $x \in X$.
 \end{lemma}
 
 \parag{proof}  We begin the proof as in the previous lemma, but now we use the fact that \ $\alpha_{k-1} \equiv 0$ \ to use the commutation lemma in order to write
 $$ P = P_{k-3}.(a -\lambda_{k-2}.b).[S_{k-2}.W]^{-1}.(a - (\lambda_k+1).b).[S_{k-1}.W^{-2}]^{-1}.(a - (\lambda_{k-1}-1).b).$$
 Then we put \ $\mathbb{G} : = \A.\chi$ \ where \ $\chi : = [S_{k-1}.W^{-2}]^{-1}.(a - (\lambda_{k-1}-1).b).\varphi$. Then \ $\mathbb{G}$ \ is not containes in  \ $\mathbb{F}_{k-1}$ \ because this would imply that \ $(a -( \lambda_{k-1}-1).b).\varphi$ \ lies in \ $\mathbb{F}_{k-1}$ and so \ $[\lambda_k - (\lambda_{k-1}-1)].b.\varphi = p_{k-1}.b.\varphi$ \ lies in \ $\mathbb{F}_{k-1}$. This impossible as \ $F_{k-1}(x)$ \ is normal for each \ $x$ \ and \ $p_{k-1} \not= 0$ \ because \ $\mathbb{E}(x)\big/\mathbb{F}_{k-2}(x)$ \ is semi-simple. So \ $\mathbb{G}_1$ \ is not semi-simple and we have \ $d(\mathbb{G}(x)) = 2$ \ for each \ $x \in X$. $\hfill \blacksquare$\\
 
 Note that we have \ $L(\mathbb{E}) = L(\mathbb{G})$ \ in the situation of the previous lemma, thanks to lemma \ref{descente}.

\begin{prop}\label{rank 2}
Let  \ $\mathbb{E}$ \ be an holomorphic family of \ $[\lambda]-$primitives frescos parametrized by an  complex space \ $X$. Assume that we have \ $d(\mathbb{E}(x)) =  2$ \ for each \ $x \in X$. Then there exists a closed  analytically constructible subset \ $Y$ \ in \ $X$ \ with no interior points in \ $X$ \ such that the family \ $L(\mathbb{E}(x)), x \in X \setminus Y$ \ is holomorphic.
\end{prop}

\parag{proof} We argue by induction on the rank \ $k$ \ of \ $\mathbb{E}$. The case of rank \ $2$ \ is obvious as we have a family of themes, so \ $L(\mathbb{E}) = F_1(\mathbb{E})$ \ in this case.\\
Assume the proposition proved for rank \ $k-1 \geq 2$, and consider the rank \ $k$ \ case. Either the holomorphic family \ $\mathbb{E}\big/\mathbb{F}_{k-2}$ \ is on the complement of a closed analytic set \ $Y \subset X$ \ with no interior point, a family of themes, and we apply the lemma \ref{induction key 1}, either it is a family of semi-simple frescos and we apply the lemma \ref{induction key 2}.

 \begin{prop}\label{hol. $L(E)$}
 Let \ $\mathbb{E}$ \ be an holomorphic family of \ $[\lambda]-$primitives frescos parametrized by an  complex space \ $X$. Assume that  we have \ $d(\mathbb{E}(x)) \geq 2$ \ for each \ $x \in X$. Then there exists a closed analytically constructible subset \ $Y$ \ in \ $X$ \ with no interior points such that the family \ $L(\mathbb{E}(x)), x \in X \setminus Y$ \ is holomorphic.
\end{prop}
 
 \parag{proof} We shall prove the statement by induction on the rank \ $k$ \  of \ $\mathbb{E}$. The first case is \ $k = 2$ \ which is obvious because for each \ $x \in X$ \ $\mathbb{E}(x)$ \ is a theme. So \ $L(\mathbb{E}(x)) = F_1(\mathbb{E}(x))$. But we have already proved that \ $\mathbb{F}_1$ \ is an holomorphic family, where \ $\mathbb{F}_1$ \ is given by the first terms of the principal J-H. (see the theorem \ref{J-H. hol.}).\\
 Assume the result proved in rank \ $k-1$ \ and \ $d(\mathbb{E}(x)) \geq 3$. Consider now the holomorphic family \ $\mathbb{F}_{k-1}$ \ given by the corank 1 terms of the J-H. sequence of \ $\mathbb{E}$. Then we have, thanks to lemma \ref{descente},  we have \ $L(\mathbb{E}(x)) = L(\mathbb{F}_{k-1}(x))$ \ for each \ $x$ \ and so the induction hypothesis allows to conclude. \\
 The only case which is left is when \ $d = 2$ \ which is proved in the proposition \ref{rank 2}.$\hfill \blacksquare$
 
 \parag{Proof of the theorem \ref{hol. S and Co}} We may assume that \ $X$ \ is irreducible. Let \ $d$ \ be the generic value of \ $d(\mathbb{E}(x))$ \ on \ $X$. As the set where \ $d(\mathbb{E}(x)) \leq d-1$ \ is a closed analytic set with no interior point, we may assume that we have \ $d(\mathbb{E}(x)) = d$ \ for each \ $x$ \ in \ $X$. We shall prove now by induction on \ $d$ \ that there exists a cac \ $Y$ \ with no interior point in \ $X$ \ such that the family \ $\Sigma^1(\mathbb{E}(x)), x \in X \setminus Y$ \ is holomorphic. The case \ $d = 1$ \ is obvious and when \ $d = 2$ \ we have \ $\Sigma^1(E) = L(E)$ \ and this case is already obtain in the proposition \ref{hol. $L(E)$}. So assume \ $d \geq 3$ \ and the result proved for \ $d-1$. The proposition \ref{hol. $L(E)$} implies that, up to replace \ $X$ \ by \ $X \setminus Y$ \ we may assume that the family \ $L(\mathbb{E}(x))$ \ is holomorphic. Then the theorem \ref{} gives the holomorphy of the family \ $\mathbb{E}\big/L(\mathbb{E}$ \ and the lemma \ref{}  shows that we may apply the inductive hypothesis to this family. As we have for each \ $x$ \ a canonical isomorphism \ $\Sigma^1(\mathbb{E}(x))\big/L(\mathbb{E}(x)) \simeq \Sigma^1(\mathbb{E}(x)\big/L(\mathbb{E}(x)))$ \ the exact sequence
 $$ 0 \to L(\mathbb{E}) \to \Sigma^1(\mathbb{E}) \to \Sigma^1(\mathbb{E})\big/L(\mathbb{E}) \to 0 $$
 and the theorem \ref{deux sur trois} allow to conclude.\\
 We obtain the case of \ $S_1(\mathbb{E})$ \ by duality using the therem \ref{hol. dual.}. $\hfill \blacksquare$

\subsection{Rank 3  fresco with nci = 0 for the J-H.}

We consider here a rank 3  fresco \ $E$ \ defined as \ $E = \A\big/\A.P $ \ where 
$$ P : = (a - \lambda_1.b).S_1^{-1}.(a - \lambda_2.b).S_2^{-1}.(a - \lambda_3.b) $$
with \ $\lambda_1 > 2$ \ rationnal and \ $p_1, p_2 \in \mathbb{N}^*$,  $\lambda_{j+1} : = \lambda_j + p_j - 1$ \ for \ $j = 1,2$ \  and \ $S_1, S_2 \in \C[[b]]$ \ satisfy \ $S_1(0) = S_2(0) = 1$.\\
We define \ $\alpha$ \ and \ $\beta$ \ the coefficient of \ $b^{p_1}$ \ in \ $S_1$ \ and of \ $b^{p_2}$ \ in \ $S_2$. Assuming \ $\alpha = 0$ \ we denote by \ $U$ \ an element in \ $\C[[b]]$ \ solution of the equation
\begin{equation*}
b.U' = p_1.(U - S_1) \tag{1}
\end{equation*}
We define \ $\gamma$ \ as the coefficient of \ $b^{p_1+p_2}$ \ in \ $U.S_2$. Remark that if we have \ $\beta = 0 $ \ then \ $\gamma$ \ does not depend of the choice of the solution \ $U$ \ of the equation \ $(1)$ :  because if \ $U_0$ \ is the solution without term in \ $b^{p_1}$ \ the general solution is \ $U = U_0 + \rho.b^{p_1}$. And then \ $U.S_2 = U_0.S_2 + \rho.b^{p_1}.S_2$ \ and the coefficient of \ $b^{p_1+p_2}$ \ in \ $b^{p_1}.S_2$ \ is \ $\beta$.\\
Remark also that if \ $\beta \not= 0$ \ we may always choose \ $U$ \ (i.e. \ $\rho$) in an unique way such that the coefficient  of \ $b^{p_1+p_2}$ \ in \ $U.S_2$ \ is \ $0$.\\

\begin{lemma}\label{ss}
With the previous notations, assume that \ $\alpha = \beta = 0$. Then the fresco \ $E$ \ is semi simple if and only iff \ $\gamma = 0$.\\
When \ $\gamma \not= 0$ \ we have \ $ d(E) = 2$ \ and \ $L(E) \simeq E_{\lambda_1}$.
\end{lemma}

\parag{proof}  Using the commuting lemma \ref{ com.lem.} we have as soon as \ $\alpha = 0$
$$ P = U^{-1}.(a -(\lambda+2).b).S_1^{-1}.U^2.(a - (\lambda_1-1).b).(U.S_2)^{-1}.(a - \lambda_3.b). $$
Assuming now that the coefficient of \ $b^{p_1+p_2}$ \ in  \ $U.S_2$ \ is zero  by our  assumption  \ $\beta= \gamma = 0$\  the commuting lemma applies again and  gives
$$ P = U^{-1}.(a - (\lambda_2+1).b).[S_1.U^{-2}.V]^{-1}.(a - (\lambda_3+1).b).[U.S_2.V^{-2}]^{-1}.(a - (\lambda_1-2).b).V^{-1} $$
where \ $V \in \C[[b]]$ \ is solution of the equation
\begin{equation*}
b.V' = (p_1+p_2).(V - U.S_2) \tag{2}
\end{equation*}
 Now if we may apply again the commuting lemma in order to "commute" \ $\lambda_2$ \ and \ $\lambda_3$, we shall obtain a J-H. sequence for \ $E$ \ with a decreasing sequence of numbers \ $\mu_j+j$ \ (explicitely \ $\mu_1 = \lambda_3+2, \mu_2 = \lambda_2, \mu_3 = \lambda_1-2$), so \ $E$ \ will be semi-simple. \\
  The necessary and sufficient to apply the commuting lemma is now that the coefficient of \ $b^{p_2}$ \ in \ $S_1.U^{-2}.V$ \ is \ $0$. Put \ $Z : = U^{-1}$. Then the equation \ $(1)$ \ gives
 $$ b.Z' = -p_1.(Z - S_1.U^{-2}) $$
 and also
 \begin{equation*}
 b.Z'.V = -p_1.(Z.V - S_1.U^{-2}.V) \tag{3}
 \end{equation*}
 The equation \ $(2)$ \ gives
 \begin{equation*}
 b.V'.Z = (p_1+p_2).(V.Z - S_2) \tag{4}
 \end{equation*}
 and adding \ $(3)$ \ and \ $(4)$ \ leads to
 $$ b.(V.Z)' = p_2.(V.Z) + p_1.S_1.U^{-2}.V - (p_1+p_2).S_2 $$
 which implies that \ $p_1.S_1.U^{-2}.V - (p_1+p_2).S_2$ \ has no term in \ $b^{p_2}$. As we know that \ $\beta = 0$, this gives that \ $E$ \ is semi-simple when \ $\gamma = 0$.\\
 
 Assuming that \ $\gamma \not= 0$, we produce a rank 2  quotient theme of \ $E$ \ after the first commutation. So \ $d(E) = 2$ \ in this case. And we obtain that \ $S_1(E)$ \ is the rank \ $2$ \ fresco given by the second step of the principal J-H. sequence of \ $E$.\\
 In order to compute \ $L(E) = \Sigma^1(E)$ \ in the cas where \ $\gamma \not= 0$ \ we must begin by commuting \ $\lambda_2$ \ and \ $\lambda_3$ \ using \ $\beta = 0$. So let \ $X \in \C[[b]]$ \ be a solution of
 \begin{equation*}
 b.X' = p_2.(X - S_2) \tag{5}
 \end{equation*}
 then the commuting lemma gives
 $$ P = (a - \lambda_1.b).[S_1.X]^{-1}.(a - (\lambda_3+1).b).[S_2.X^{-2}]^{-1}.(a - (\lambda_2-1).b).X^{-1} $$
 The claim is now that the coefficient of \ $b^{p_1+p_2}$ \ in \ $S_1.X$ \ is not \ $0$. But \ $(5)$ \ and \ $(1)$ \ give
 \begin{align*}
 & b.U'.X = p_1.(U.X - S_1.X) \\
 & b.X'.U = p_2.(X.U - S_2.U) \quad {\rm and \ adding}\\
 & b.(X.U) = (p_1+p_2).X.U - p_1.S_1.X - p_2.S_2.U
 \end{align*}
 so the coefficient of \ $b^{p_1+p_2}$ \ in \ $S_1.X$ \ is equal to \ $-\frac{p_2}{p_1}.\gamma \not= 0 $. \\
 Then the second step of the J-H. sequence of \ $E$ \ given by the above writing of  \ $P$ \ is a rank 2 theme. So its rank \ $1$ \ sub-theme \ $\simeq E_{\lambda_1}$ \ is equal to \ $L(E)$. $\hfill \blacksquare$\\
 
 \parag{The exemple}
 Choose \ $\lambda_1, p_1, p_2$ \ as above. Let the parameter space be \ $\C$ \ and consider the family of rank \ $3$ \ frescos \ $\mathbb{E}(z), z \in \C$ \ with \ $\mathbb{E}(z) : \A\big/\A.P_z$ \ and where
 $$ P_z : = (a - \lambda_1.b).(1 + b^{p_1+p_2})^{-1}.(a - \lambda_2.b).(1 + z.b^{p_2})^{-1}.(a - \lambda_3.b) .$$
 Then it is clear that \ $\mathbb{E}(z)\big/\mathbb{F}_1(z)$ \ is a rank \ $2$ \ theme for \ $z \not= 0$. For \ $z = 0$ \ we have \ $\alpha = \beta = 0$ \ (as \ $\alpha(z) \equiv 0$ \ and \ $\beta(z) = z$). Moreover we have \ $\gamma = - \frac{p_1}{p_2} \not= 0$ \ because \ $U = \frac{1}{p_1} - \frac{p_1}{p_2}.b^{p_1+p_2} $ \ and \ $S_2(0) =  1$. So \ $L(\mathbb{E}(0) \simeq E_{\lambda_1}$ \ thanks to the previous lemma.
 For the computation of \ $L(\mathbb{E}(z))$ \ when \ $z \not= 0$ \ we will produce a normal rank \ $2$ \ subtheme of \ $\mathbb{E}(z)$. As \ $\alpha = 0$ \ we apply the commuting lemma \ref{} \ but choosing  $$U : =  \frac{1}{p_1} - \frac{p_1}{p_2}.b^{p_1+p_2}  + \frac{p_1}{z.p_2}.b^{p_1} $$
 inorder that the coefficient of \ $b^{p_1+p_2}$ \ in \ $U.S_2$ \ in zero (note that we use here that \ $\beta(z) = z \not= 0 $ \ to choose a convenient solution \ $U$ \ of \ $(1)$). \\
 Then we may apply the commuting lemma a second time because with our choice the equation \ $(2)$ \ has a solution \ $V \in \C[[b]]$. Then we obrain that for \ $z \not= 0$ \ we have 
 $$ P_z = U^{-1}.(a - (\lambda_2+1).b).[S_1.U^2.V]^{-1}.(a - \lambda_3.b).[U.S_2.V^{-2}]^{-1}.(a - (\lambda_1-2).b).V^{-1} .$$
 In order to prove that the second term of the J-H. associated to this writing of \ $E$ \ is a theme we must compute the coefficient of \ $b^{p_2}$ \ in \ $S_1.U^2.V$. But using equations \ $(1),(2)$ \ we have, puttong \ $Z : = U^{-1}$
 \begin{align*}
 & b.Z' = -p_1.(Z - S_1.U^{-2}) \quad {\rm and \ so} \\
 & b.Z'.V = -p_1.(Z.V - S_1.U^{-2}.V) \\
 & b.V'.Z = (p_1+p_2).(V.Z - S_2) \quad {\rm then \ adding} \\
 & b.(V.Z)' = p_2.V.Z  + p_1.S_1.U^{-2}.V - (p_1+p_2).S_2
 \end{align*}
 and this shows that  the coefficient of \ $b^{p_2}$ \ in \ $S_1.U^2.V$ is equal to \ $\frac{p_1+p_2}{p_1}.z \not= 0 $. So we conclude that for \ $z \not= 0$ \ we have
 $$ L(\mathbb{E}(z)) \simeq E_{\lambda_2+1} .$$
 So the family \ $L(\mathbb{E}(z)), z \in \C$ \ is not an holomorphic family.\\
 
\parag{Remark} It is not difficult to see that the holomorphic family \ $\mathbb{F}_2(z)$ \ associated to the second term of the principal J-H. sequence of the \ $\mathbb{E}(z)$ \ is equal to the family \ $S_1(\mathbb{E}(z))$ \ in this example. But of course, using the twisted duality produces an holomorphic family of frescos \ $\mathbb{G}(z), z \in \C$ \  for which the corresponding family \ $S_1(\mathbb{G}(z))$ \ is not holomorphic.

\subsection{d-constant for the principal J-H.}

As a consequence of  the theorem \ref{Themes et ss}, for an holomorphic family \ $\mathbb{E}$ \ of frescos parametrized by a reduced and  irreducible complex space \ $X$, there exists  a Zariski dense open set in \ $X$ \ on which the ss-depth of \ $\mathbb{E}(x)$ \ is maximal. The corollary of the  following proposition shows that this constancy of \ $d(\mathbb{E}(x))$ \ is also true for each term of the principal J-H. sequence on the same open set.

\begin{prop}\label{var. d}
Let \ $\mathbb{E}$ \ be an holomorphic family of rank \ $k$ \ $[\lambda]-$primitive frescos parametrized by an reduced complex space \ $X$. Let \ $\mathbb{F}_{k-1}$ \ be the corank \ $1$ \ sub-family associated to the principal Jordan-H{\"o}lder sequence of \ $E$. Let
$$ Z : = \{ x \in X \ / \  d(\mathbb{F}_{k-1}(x)) = d(\mathbb{E}(x)) \} .$$
Then \ $Z$ \ is a closed analytic subset of \ $X$.
\end{prop}

\parag{proof} Choose \ $\mu \in [\lambda]$ \ large enough and define
\begin{align*}
& \mathcal{E} : = Ker\big[ ( a - \mu.b) : \mathbb{E} \to \mathbb{E} \big]  \\
&  \mathcal{F}_{k-1} : = Ker\big[ ( a - \mu.b) : \mathbb{F}_{k-1} \to \mathbb{F}_{k-1} \big] 
\end{align*}
Then these sheaves are coherent \ $\mathcal{O}_X-$modules because, thanks to the lemma  \ref{ordre suff.},  locally on \ $X$, there exists  \ $N \gg 1$ \  such that they are isomorphic with the analoguous kernels computed respectively in the \ $\mathcal{O}_X-$coherent quotients  \ $\mathbb{E}\big/b^N.\mathbb{E}$ \ and \ $\mathbb{F}_{k-1}\big/b^N.\mathbb{F}_{k-1}$, and the the map \ $(a - \mu.b)$ \ is \ $\mathcal{O}_X-$linear. Now we have an \ $\mathcal{O}_X-$linear inclusion \ $\mathcal{F}_{k-1} \to \mathcal{E} $ \ and the coherent sheaf \ $\mathcal{Q} : =  \mathcal{E}\big/\mathcal{F}_{k-1}$ \ has \ $Z$ \ as support. This proves our assertion. $\hfill \blacksquare$

\parag{Remark} The fiber at \ $x \in X$ \ of the sheaf \ $\mathcal{Q} $ \ has dimension \ $0$ \ when \ $x \not\in Z$ \ and dimension  \ $1$ \ when \ $x$ \ is in \ $X$.\\

\begin{cor}\label{stab. d= cste}
Assume that the holomorphic family \ $\mathbb{E}$ \ is such that for each \ $x $ \ in \ $X$ \ we have \ $d(\mathbb{E}(x)) = d$. Assume also that \ $X$ \ is irreducible. Then for each \ $j \in [1,k]$ \ there exists an integer \ $d_j \in [1,d]$ \ such that for each \ $x\in X$ \ we have
$$ d(\mathbb{F}_j(x)) = d_j .$$
\end{cor}

\parag{proof} Define  the subset \ $Z_j : = \{ x \in X \ / \  d(\mathbb{F}_j(x)) = d(\mathbb{F}_{j+1}(x))\}$ \ for each  \ $j \in [1,k-1]$. Then denote by \ $J : = \{ j \in [1,k-1] \ / \  Z_j \not= X\} $ \ and let \ $l : = Card(J)$. For each \ $j \in J$ \ the subset \ $Z_j$ \ is closed in \ $X$ \ with no interior point. So on the dense open set \ $X \setminus \cup_{j \in J} Z_j$ \ we have 
$$ d(\mathbb{E}(x)) = l+1 = d .$$
So if a point \ $x \in X$ \ lies in \ $h$ \ of the subsets \ $Z_j, j \in J$ \ we have \ $d(\mathbb{E}(x)) = l+1 -h$. This implies \ $h = 0$ \ because, by assumption, we have \ $d(\mathbb{E}(x)) = d$ \ for each \ $x \in X$.\\
Then each \ $Z_j$ \ is  empty for \ $j \in J$. Then \ $d_j = 1 + \delta_j$ \ where \ $\delta_j : = Card(J \cap [1,j-1])$ \ is the number of \ $Z_i$ \ which are empty with \ $i \in [1,j-1]$.$\hfill \blacksquare$

 \section{Two constructions of holomorphic parameters invariant by change of variable.}

\subsection{Holomorphic parameters.}

\begin{defn}\label{Hol. param. 1}
Fix the Bernstein polynomial  \ $B$ \ of a \ $[\lambda]-$primitive rank \ $k$ \ fresco, and let \ $Z(B)$ \ be the set of isomorphism classes of  \ $[\lambda]-$primitive rank \ $k$ \ frescos with Bernstein polynomial \ $B$. A map \ $f : Z(B) \to Y$ \ with values in a reduced complex space \ $Y$ \ will be called an {\bf holomorphic parameter} if for any holomorphic family of \ $\mathbb{E}$ \ of  \ $[\lambda]-$primitive rank \ $k$ \ frescos parametrized by a reduced complex space \ $X$ \  with Bernstein polynomial \ $B$, the function \ $ f\circ I : X \to Y$ \ is holomorphic, where \ $I : X \to Z(B)$ \ is the map which associates to \ $x$ \ the isomorphism class of \ $\mathbb{E}(x)$.
\end{defn}

Remark that, as the Bernstein plynomial is invariant by a change of variable \ $\theta$, such a \ $\theta$ \ induces a bijection \ $\theta_* : Z(B) \to Z(B)$ \ which is "holomorphic" in the sens that the notion of holomorphic family is preserved by change of variable  thanks to the  proposition  \ref{chgt. var. hol.}. \\

The  notion of holomorphic parameter does not look  so interesting if there exists a "well described" universal family of \ $[\lambda]-$primitive frescos with the prescribe Bernstein polynomial \ $B$. But, in general, such an universal family does not exists even locally around a given isomorphism class (see for instance [B.9-b] and use the theorem  \ref{Themes et ss} which gives the fact that themes are open in frescos). Moreover, even in the case where such an universal family exists it is not so easy to describe the set \ $Z(B)$ \ and the complex structure on it. A fortiori holomorphic maps defined on it.\\

As explain in the introduction, the main goal will be to construct in a rather general case holomorphic parameters which are invariant or quasi-invariant by any change of variable. Remark that, even in the case where we have an universal family, the construction of an holomorphic map on the parameter space \ $Z(B)$ \ (with a suitable structure of a reduced complex space)   which is constant on equivalence classes defined by all change of variables is not so obvious because we have to analyse when two given isomorphism classes (with the same \ $B$) may correspond by some change of variable.

\subsection{Invariance of the principal  parameters of  \ $[\lambda]-$primitive themes.}

We shall study in this section the action of the change of variable on the  (universal) family of rank \ $2$ \ $[\lambda]-$primitive themes. So we shall fix a rational number \ $\lambda_1 > 1$ \ and an integer \ $p_1 \geq 1$. In the case \ $p_1 = 0$ \ there exists only one isomorphism class of rank \ $2$ \ $[\lambda]-$primitive theme with these fundamental invariants : 
 $$\hat{A}\big/\hat{A}.(a - \lambda_1.b).(a - (\lambda_1-1).b) .$$
  So any change of variable preserves this isomorphism class. For \ $p \geq 1$ \ the isomorphism class is determined by the principal parameter \ $\gamma$ \ (see [B.09-b]) :
$$ E(\gamma) : = \hat{A}\big/\hat{A}.(a - \lambda_1.b).(1 + \gamma.b^{p_1})^{-1}.(a - (\lambda_1+p_1-1).b).$$

\begin{defn}\label{compatible}
Let \ $ f : E \to F$ \ be a \ $\C-$linear map between two \ $\C[[b]]-$modules. We shall say that \ $f$ \ is \ {\bf $b-$compatible} if for each integer \ $n$ \ we have 
 $$f(b^n.E) \subset b^n.F.$$
\end{defn}

\begin{defn}\label{polyn. depend.}
Let \ $E$ \ be a free finite type \ $\C[[b]]-$module and let \ $f_{\theta} : E \to E$ \ be a \ $\C-$linear map which is \ $b-$compatible for each value of the parameter \ $\theta \in \C^N$. We say that \ $f$ \ {\bf depends polynomially on \ $\theta$} when for each \ $x \in E$ \ and each integer \ $n$ \ the map \ $\C^N \to E\big/b^n.E$ \ given by \ $\theta \mapsto [f_{\theta}(x)] \in E\big/b^n.E$ \ is polynomial.
\end{defn}

\parag{Remarks} 
\begin{enumerate}
\item In the case of an element \ $S_{\theta} \in \C[[b]]$ \ acting by multiplication on a finitely generated \ $\C[[b]]-$module the polynomial dependance on \ $\theta$ \ is equivalent with the fact that each coefficient of \ $S_{\theta}$ \ is a polynomial function of \ $\theta$.
\item If we have two \ $\C-$linear endomorphisms \ $f$ \ and \ $g$ \ of \ $\C[[b]]$ \ which are \ $b-$compatible and depend polynomially on \ $\theta$, the composed endomorphism is again \ $b-$compatible and depends polynomially in \ $\theta$. This comes from the fact that the coefficient of \ $b^n$ \ in \ $g(f(x))$ \ depends only on the coefficients of \ $f(x)$ \ in \ $\C[[b]]\big/b^{n+1}.\C[[b]]$ \ and of the endomorphism induced by \ $g$ \ on this quotient ; but the matrix of this induced  endomorphism are polynomials in \ $\theta$.
\end{enumerate}

\bigskip

We shall need a more precise version of the invariance of the isomorphism class of a rank 1 fresco by a change of variable.

\begin{lemma}\label{rk 1 polyn. depend.}
We consider a change of variable \ $\theta(a) = a + \sum_{j=2}^N \theta_j.a^j$ \ as an element in \ $\C^{N-1}$. For each \ $\lambda \in \C$ \ there exists an unique (invertible) element \ $S^{\lambda}_{\theta} \in \C[[\beta]]$, depending polynomially on \ $\theta \in \C^{N-1}$, such that  the element \ $\varepsilon_{\theta}^{\lambda}$ \ in \ $E_{\lambda} : = \hat{A}\big/\hat{A}.(a - \lambda.b)$ \ defined as \ $\varepsilon_{\theta}^{\lambda} : = S^{\lambda}_{\theta}(\beta).e_{\lambda}$ \  satisfies :
$$ \alpha_{\theta}.\varepsilon_{\theta}^{\lambda} = \lambda.\beta_{\theta}.\varepsilon_{\theta}^{\lambda} $$
where \ $\alpha_{\theta} : = \theta(a)$ \ and \ $\beta_{\theta} : = b.\theta'(a)$.
\end{lemma}

\parag{Proof} We shall prove first that there exists, for each integer \ $n$, an invertible element \ $\chi_n \in \C[[b]]$, depending polynomially on \ $\theta \in \C^{N-1}$, such the following identity holds in \ $E_{\lambda}$ :
\begin{equation*}
b^n.e_{\lambda} = \beta^n.\chi_n(\beta).e_{\lambda} \tag{*}
\end{equation*}
where \ $\beta_{\theta} : = b.\theta'(a)$.\\
Remark that \ $\alpha_{\theta} : = \theta(a)$ \ and \ $\beta_{\theta}$ \ are two \ $\C-$linear \ $b-$compatible endomorphisms  of \ $E_{\lambda}$ \ which depends polynomially on \ $\theta \in \C^{N-1}$. Now the equality 
 $$b^n.\C[[b]].e_{\lambda} = \beta^n.\C[[\beta]].e_{\lambda}$$
 implies the existence and the uniqueness of \ $\chi_n$ \ for each \ $n \in \mathbb{N}$. So it is enough to show that \ $\chi_n$ \ depends polynomially on \ $\theta$.  Fix an integer \ $p \gg 1$ \ and consider the basis of the vector space \ $V_p : = E_{\lambda}\big/b^p.E_{\lambda}$. We have two basis \ $b^q.e_{\lambda}$ \ and \ $\beta^q.e_{\lambda}, q \in [0,p-1]$ \ of \ $V_p$ \ and they satisfy
$$ \beta^q.e_{\lambda} = b^q.e_{\lambda} + b^{q+1}.V_p \quad \forall q \in [0,p-1] .$$
So the corresponding  base change is triangular with entries equal to  $1$  on the diagonal. As \ $\beta$ \ depends polynomially of \ $\theta$, the inverse of this matrix of base change has polynomial entries in \ $\theta$ \ and this allows to conclude.\\

Now we shall prove the assertion of the lemma. We assume \ $\lambda$ \ fixed, and so we shall omit some dependance in \ $\lambda$ \ in what follows.\\
Remark that we have
$$ \alpha.e_{\lambda} = \lambda.b.e_{\lambda} + \sum_{j=2}^N L_j.\theta_j.b^j.e_{\lambda} $$
where \ $L_j : = \lambda.(\lambda+1) \dots (\lambda+j-1) $. So using our previous statement, we may write
$$ \alpha.e_{\lambda} = \lambda.\beta.e_{\lambda} + \beta^2.R_{\theta}(\beta).e_{\lambda}$$
where \ $R_{\theta} \in \C[[b]]$ \ depends polynomially in \ $\theta$. Now the desired equality
$$ \alpha.S_{\theta}.e_{\lambda} = \lambda.\beta.S_{\theta}.e_{\lambda} $$
is equivalent to
\begin{equation*}
 S_{\theta}.\big[ \lambda.\beta.e_{\lambda} + \beta^2.R_{\theta}(\beta).e_{\lambda}\big] + \beta^2.S_{\theta}'.e_{\lambda} =  \lambda.\beta.S_{\theta}.e_{\lambda} 
 \end{equation*}
 and after simplification we find the differential equation :
 $$ S_{\theta}' + R_{\theta}.S_{\theta} = 0 .$$
 Let \ $\tilde{R}_{\theta}$ \ be the primitive without constant term of \ $R_{\theta}$.\\
  Then we obtain \ $S_{\theta} = exp[-\tilde{R}_{\theta}]$, and as \ $R_{\theta}$ \ depends polynomially on \ $\theta$, the same is true for \ $\tilde{R}_{\theta}$ \ and \ $ exp[-\tilde{R}_{\theta}]$. $\hfill \blacksquare$\\

\begin{prop}\label{chgt. var. theme rk 2}
Let \ $E : = \hat{A}\big/\hat{A}.(a - \lambda_1.b).(1 + z.b^p)^{-1}.(a - (\lambda_1+ p-1).b)$ \ where we assume that \ $\lambda_1 >1$ \ is a rational number and \ $p \in \mathbb{N}^*$. Let \ $\theta \in \C[a]$ \ be an unitary polynomial change of variable \ $\theta(a) = a + \sum_{j=2}^N \ \theta_j.a^j$. Then \ $\theta_*(E)$ \ is isomorphic to
$$ \hat{A}\big/\hat{A}.(a - \lambda_1.b).(1 + z_{\theta}.b^p)^{-1}.(a - (\lambda_1+ p-1).b) $$
where \ $z_{\theta}$ \ is a polynomial in \ $\theta = (\theta_2, \dots, \theta_{N})  \in \C^{N-1}$.
\end{prop}

\parag{Proof} Let \ $e_1,e_2$ \ be the standard \ $\C[[b]]-$basis associated to the isomorphism  \ $E : = \hat{A}\big/\hat{A}.(a - \lambda_1.b).(1 + z.b^p)^{-1}.(a - (\lambda_1+ p-1).b)$. Then \ $\varepsilon_2 : = S_{\theta}^{\lambda_2}(\beta).e_2$ \ satifies
$$ (\alpha - \lambda_2.\beta).\varepsilon_2 \in F_1 : = \C[[b]].e_1= \C[[\beta]].e_1 \simeq E_{\lambda_1}.$$
Now \ $z_{\theta}$ \ is the coefficient of \ $\beta^p$ \ when we write
$$  (\alpha - \lambda_2.\beta).\varepsilon_2 = \sum_{\nu=0}^{+\infty} \beta^{\nu}.\varepsilon_1 $$
with \ $\varepsilon_1 : = S_{\theta}^{\lambda_1}(\beta).e_1$. Then the polynomial dependance of the \ $\chi_n$ \ and of \ $(S_{\theta}^{\lambda_i}), i = 1, 2$ \ allow to conclude. $\hfill \blacksquare$\\

\begin{thm}\label{inv. param.}
Let \ $E$ \ be a \ $[\lambda]-$primitive rank \ $k$ \ theme with fundamental invariants \ $\lambda_1, p_1, \dots, p_{k-1}$. Let \ $z_1, \dots, z_{k-1}$ \ the parameters of the rank \ $2$ \ $[\lambda]-$primitive themes \ $F_2, F_3\big/F_1, \dots, F_k\big/F_{k-2}$. Let \ $\theta$ \ a polynomial change of variable \ $\theta : = r.a + \sum_{j=2}^N \theta_j.a^j $ \ with \ $r \not= 0 $. Then \ $\theta_*(E)$ \ is a rank \ $k$ \ $[\lambda]-$primitive theme with fundamental invariants \ $\lambda_1, p_1, \dots, p_{k-1}$. The parameters of the rank \ $2$  \ $[\lambda]-$primitive themes 
 $$\theta_*(F_2), \theta_*(F_3\big/F_1), \dots, \theta_*(F_k\big/F_{k-2})$$
 are the numbers \ $r^{p_1}.z_1, \dots, r^{p_{k-1}}.z_{k-1}$.
\end{thm}

\parag{Proof} It is enough to prove the result in the rank \ $2$ \ case because if we consider a \ $[\lambda]-$primitive theme \ $E$ \ with rank \ $k$ \ and fundamental invariants \ $\lambda_1, p_1, \dots, p_{k-1}$ \ we have the following facts :
\begin{enumerate}
\item For any change of variable \ $\theta$ \ the fresco \ $\theta_*(E)$ \ is again a $[\lambda]-$primitive  rank \ $k$ \ theme  with fundamental invariants \ $\lambda_1, p_1, \dots, p_{k-1}$.
\item The quasi-invariance\footnote{for the character \ $\theta \mapsto \theta'(0)^{p_j}$.} of the parameter of the rank \ $2$ \ quotient theme \
$F_{j+1}\big/F_{j-1}$ \ $  j \in [1,k-1]$ \  is obtained by the rank \ $2$ \ case, as \ $\theta_*(F_j)$ \ is the (unique)  normal rank \ $j$ \ sub-theme of \ $\theta_*(E)$.
\end{enumerate}

To prove the rank \ $2$ \ case, it is enough to consider polynomial \ $\theta$ \ because for \ $N $ \ large enough, it is easy to see that \ $(\theta_N)_*(E)$ \ is isomorphic to \ $\theta_*(E)$ \ where \ $\theta_N$ \ is the expansion up to the order \ $N$ \ of \ $\theta$. Then the result in now an easy consequence of the proposition \ref{chgt. var. theme rk 2} because the only polynomial functions on \ $\C^N$ \ which never vanish are the constant functions. $\hfill \blacksquare$

\subsection{Classification of rank 3 semi-simple frescos.}
Let \ $\lambda > 2$ \ be a rational number and  \ $p_1 \geq 2, p_2 \geq 1$ \ be two integers. Put \ $\lambda_1 : = \lambda, \lambda_{j+1} = \lambda_j + p_j -1$ \ for \ $j = 1,2$ \ and for \ $\gamma \in \C$ \ define
\begin{equation*}
P_{\gamma} : = (a - \lambda_1.b).(1 + \gamma.b)^{-1}.(a - \lambda_2.b).(a - \lambda_3.b). \tag{@}
\end{equation*}
and  \ $E(\gamma) : = \hat{A}\big/\hat{A}.P_{\gamma}$.

\begin{prop}\label{rg 3 ss 1}
Any semi-simple  \ $[\lambda]-$primitive rank \ $3$ \  fresco with fundamental invariants 
 \ $(\lambda_1, p_1, p_2)$\  is isomorphic to some  \ $E(\gamma)$. Moreover, for \ $p_2 \geq 2$,  $E(\gamma)$ \ and \ $E(\gamma')$ \ are not isomorphic for  \ $\gamma \not= \gamma'$. For \ $p_2 = 1$,  each \ $E(\gamma)$ \ is isomorphic to \ $E(0)$.
\end{prop}

\parag{Proof} We begin by the proof that each fresco \ $E(\gamma)$ \ is semi-simple. Using proposition 4.1.4 of [B.11] it is enough to show that \ $E(\gamma)$ \ has a Jordan-H{\"o}lder sequence such that the  numbers corresponding to the rank \ $1$ \ quotients are the numbers  \ $\lambda_3 +  2, \lambda_2 , \lambda_1 - 2$ \  which correspond to is the strictly decreasing order of the  sequence \ $\lambda_1+1, \lambda_2+2, \lambda_3+3$.\\
First the equality \ $(a - \lambda_2.b).(a - \lambda_3.b) = (a - (\lambda_3+1).b).(a - (\lambda_2-1).b) $ \  shows that   \ $E(\gamma) $ \ is isomorphic to \ $ \A\big/\A.(a - \lambda_1.b).(1 +\gamma.b)^{-1}.(a - (\lambda_3+1).b).(a - (\lambda_2-1).b) $. Now we may use the commuting  lemma \ref{ com.lem.} (see also  3.5.1 of [B.09-a])  as  \ $\lambda_3+1 - (\lambda_1-1) = p_1 + p_2 \not= 1$. It implies the equality in \ $\hat{A}$  :
\begin{align*}
&   (a - \lambda_1.b).(1 +\gamma.b)^{-1}.(a - (\lambda_3+1).b).(a - (\lambda_2-1).b) = \\
& \qquad  U^{-1}.(a - (\lambda_3+2).b).U.(1 +\gamma.b)^{-1}.U.(a -(\lambda_1-1).b).U^{-1}.(a - (\lambda_2-1).b)
\end{align*}
where \ $U$ \ is a solution in \ $\C[[b]]$ \ of the differential equation  \ $b.U' -\delta.U = \delta.(1 + \gamma.b)$ \ with \ $\delta = \lambda_3+2 -\lambda_1 = p_1+p_2$.  So we can choose
$$ U = -(1 + \rho.b) \quad {\rm with} \quad  \rho = \frac{(p_1+p_2).\gamma}{p_1 + p_2 + 1} .$$
Then, applying again {\it loc.cit.} to \ $(a -(\lambda_1-1).b).U^{-1}.(a - (\lambda_2-1).b)$, which is possible as \ $\lambda_2 - (\lambda_1-1) = p_1 \not= 1$, we obtain a J-H. sequence for \ $E(\gamma)$ \ with the required property.\\

Consider now a rank \ $3$ \ \ $[\lambda]-$primitive semi-simple fresco \  $E$ \ with fundamental invariants \ $\lambda_1, p_1,p_2$. Then the theorem 3.4.1 of [B.09-a] allows to find an invertible element \  $S_1$ \ in \ $ \C[[b]]$ \ such that 
\begin{enumerate}[i)]
\item \ $E \simeq \A\big/\A.(a - \lambda_1.b).S_1^{-1}.(a - \lambda_2.b).(a - \lambda_3.b) $.
\item \ $S(0) = 1$ 
\item  There is no term in \ $b^{p_1}$ \ or in \ $b^{p_1+p_2}$ \  in \ $S_1$. 
\end{enumerate}
Let me explain the condition iii). If \ $S_1$ \ has a non zero term in \ $b^{p_1}$ \ the normal submodule \ $F_2$ \ in the principal J-H. sequence of \ $E$ \ is a rank \ $2$ \ theme, and this contradicts the semi-simplicity of \ $E$. In an analoguous way, if \ $S_1$ \ has a non zero term in \ $b^{p_1+p_2}$, the equality \ $(a - \lambda_2.b).(a - \lambda_3.b) = (a - (\lambda_3+1).b).(a - (\lambda_2-1).b) $ \ in \ $\hat{A}$ \ implies that \ $E$ \ is isomorphic to \ $\hat{A}\big/\hat{A}.(a - \lambda_1.b).S_1^{-1}.(a - (\lambda_3+1).b).(a - (\lambda_2-1).b) $ \ and as  \ $\lambda_3 + 1-(\lambda_1-1) =  p_1+p_2$, we obtain again a rank \ $2$ \ sumodule of \ $E$ \ which is a theme, contradicting the semi-simplicity of \ $E$.\\
Let \ $e_1, e_2,e_3$ \ be the standard basis of \ $E$ \ which is associated to the above isomorphism. We have 
$$ (a - \lambda_3.b).e_3 = e_2 \quad (a - \lambda_2.b).e_2 = S_1.e_1 \quad {\rm and} \quad (a - \lambda_1.b).e_1 = 0 .$$
Put \ $F_2 : = \C[[b]].e_1 \oplus \C[[b]].e_2$. If we consider \ $\varepsilon_2 : = e_2 + \Sigma.e_1$ \  as the generator of \ $F_2$, with $\Sigma \in \C[[b]]$, we obtain
\begin{align*}
& (a - \lambda_2.b).\varepsilon_2 = S_1.e_1 + b^2.\Sigma'.e_1 +  \Sigma.(\lambda_1 - \lambda_2).b.e_1 \\
& \qquad  = (S_1 + b^2.\Sigma' -(p_1-1).b.\Sigma).e_1
\end{align*}

Then put \ $S_1 : = 1 + \gamma.b + b^2.T $, and  look for \ $\Sigma$ \ such that
$$ b^2.T + b^2.\Sigma' - (p_1-1).b.\Sigma = 0 .$$
As we assumed that \ $S_1$ \ has no term in \ $b^{p_1}$, it implies that \ $b.T$ \ has no term in \ $b^{p_1-1}$. So we may choose a solution \ $\Sigma \in \C[[b]]$ \  with no constant term (note that \ $p_1 \geq 2$).\\
Then defining \ $\varepsilon_1 = e_1$ \ we have
$$ (a - \lambda_2.b).\varepsilon_2 = (1 + \gamma.b).\varepsilon_1 \quad {\rm and}\quad (a - \lambda_1.b).\varepsilon_1 = 0 .$$
Now we look for  \ $\varepsilon_3 = e_3 + V.\varepsilon_2 + W.\varepsilon_1$ \  such that \ $(a - \lambda_3.b).\varepsilon_3 = \varepsilon_2$. This implies the following relations
\begin{align*}
& e_2 + b^2.V'.\varepsilon_2 + (\lambda_2-\lambda_3).b.V.\varepsilon_2 + V.(1 + \gamma.b).\varepsilon_1 + b^2.W'.\varepsilon_1 + (\lambda_1-\lambda_3).b.W.\varepsilon_1 = e_2 + \Sigma.e_1
\end{align*}
So we must have :
\begin{align*}
& b^2.V' - (p_2-1).b.V = 0 \quad \quad  {\rm and} \\
& b^2.W' - (p_1+p_2 -2).b.W = \Sigma - V.(1 + \gamma.b) \tag{@@}
\end{align*}
The first equation gives \ $V = \tau.b^{p_2-1} $, and we deduce that the second equation is equivalent to :
\begin{equation*}
 b.W' - (p_1+p_2 - 2).W = \frac{1}{b}.[\Sigma - \tau.b^{p_2-1}.(1+\gamma.b)] . \tag{@@@}
 \end{equation*}
 Remark that we may divide by \ $b$ \ in \ $\C[[b]]$ \ in the right handside because  \ $p_2 \geq 2$ \ and \ $\Sigma(0) = 0$. We shall have a solution in \ $\C[[b]]$ \ as soon as the right handside has no term in  \ $b^{p_1+p_2-2}$. But \ $S_1$ \ has no term in  \ $b^{p_1+p_2}$, and so \ $b.T$ \ and  \ $\Sigma$ \   have no term in  \ $b^{p_1+p_2-1}$.  So there exists a solution \ $W \in \C[[b]]$. This gives an isomorphism between \ $E$ \ and \ $E(\gamma)$ \ where  \ $\gamma = S_1'(0)$.\\
 For \ $p_2 = 1$, we must choose \ $\tau = 0$ \ in the previous computation and we may conclude as before because we know that \ $\Sigma$ \ has no constant term ($b.T$ \ has no constant term and \ $p_1\geq 2$).\\
 
 Now we shall prove that for  \ $p_2 \geq 2$ \ there is no isomorphism between \ $E(\gamma)$ \ and \ $E(\gamma')$ \ for \ $\gamma \not= \gamma'$. We go back to the previous proof with  \ $S_1 : = 1 + \gamma'.b$. It gives the equation
 $$ b.(\gamma' - \gamma) + b^2.\Sigma' - (p_1-1).b.\Sigma = 0 $$
for which solutions are given by
 $$\Sigma : = \frac{\gamma' - \gamma}{p_1-1} + \sigma.b^{p_1-1} $$
 where \ $\sigma \in \C$ \ is arbitrary choosen. Then we look for \ $\varepsilon_3 = e_3 + V.\varepsilon_2 + W.\varepsilon_1$ \ such that \ $(a - \lambda_3.b).\varepsilon_3 = \varepsilon_2$ \  with \ $\varepsilon_2 : =  e_2 + \Sigma.e_1$. Then we have to verify the conditions \ $(@@)$, which imply  \ $(@@@)$. But in order to have a solution  \ $W \in \C[[b]]$ \ we have to check that \ $\Sigma$ \ has no constant term. As  
 $$\Sigma(0) = \frac{\gamma' - \gamma}{p_1-1} $$
 this forces \ $\gamma = \gamma'$ \ and our statement is proved.\\
 
 For \ $p_2 = 1$ \ we reach
$$\Sigma(0) =  \frac{\gamma' - \gamma}{p_1-1} + \sigma .$$
The choice of \ $\sigma$ \ allows to have a solution for any given \ $\gamma'$ \ as soon as \ $\Sigma$ \ has no term in \ $b^{p_1+p_2-1} = b^{p_1}$. But this is true. To be explicit, put
  \begin{align*}
  & \varepsilon_3 = e_3 + \sigma.e_2 + \frac{\sigma.\gamma}{p_1-1}.e_1 \\
  & \varepsilon_2 = e_2 + \sigma.e_1 \quad {\rm and} \quad \varepsilon_1 = e_1
  \end{align*}
  then we obtain an isomorphism from \ $E(\gamma)$ \ onto \ $E(\gamma')$ \ by choosing the constant  \
   $\sigma =  -\frac{\gamma' - \gamma}{p_1-1} .$ \ $\hfill \blacksquare$

\bigskip

\subsection{Change of variable for the universal family of rank 3 semi-simple frescos.}

We know that the Bernstein polynomial of a fresco  and  the semi-simplicity of a fresco  are  invariant by change of variable. We want to study now the effect of a change of variable of the form 
   $$\alpha : = \theta(a) = a + \rho.a^2, \qquad \beta : = b.\theta'(a) = b + 2\rho.b.a $$
   on the family described above.
    
  \parag{A computation} We look for three elements \ $S,Z_0,Z_1 \in \C[[b]]$ \  such that \\
   $S(0)  = 1$ \ and satisfying the following relation in the algebra \ $\hat{A}$, where \ $\lambda, \rho$ \ are given complex numbers
   \begin{equation*}
  \alpha - \lambda.\beta = (a + \rho.a^2 - \lambda.(b + 2\rho.b.a)).S = (Z_0 + a.Z_1)(a - \lambda.b) .\tag{1} 
  \end{equation*}
 The coefficient of \ $a^2$ \ imposes
  \begin{equation*}
  Z_1 = \rho.S . \tag{2}
  \end{equation*}
  So we get
  \begin{align*}
  & a^2.S - a.S.a = a^2.S - a.(a.S - b^2.S') = a.b^2.S' = b^2.S'.a + b^2.(b^2.S')' \\
  & \qquad \qquad = b^2.S'.a  + 2b^3.S' + b^4.S'' .
   \end{align*}
   Now the sum \ $T$ \ of the terms in \ $\C[[b]].a$ \ in the relation \ $(1)$ \ is given by :
   \begin{align*}
   & T = \rho.b^2.S'.a +  a.S - 2\lambda.\rho.b.a.S - Z_0.a + \lambda.\rho.a.S.b \quad modulo \  \C[[b]] \\
   & T = \Big(\rho.b^2.S' + S -  2\lambda.\rho.b.S - Z_0  + \lambda.\rho.S.b\Big).a  \quad modulo \   \C[[b]] 
   \end{align*}
  and we obtain :
   \begin{equation*}
    Z_0 = (1 - \lambda.\rho.b).S + \rho.b^2.S' . \tag{3}
    \end{equation*} 
   The sum \ $U$ \  of the terms in  \ $\C[[b]]$ \ is given by :
    \begin{align*}
    & U =  \rho.\big[ 2b^3.S' + b^4.S''\big] + b^2.(1 - 2\lambda.\rho.b).S' + \lambda.\rho.b^2.(b.S)' - \lambda.b.S + \lambda.b.Z_0 \\
    & U = \rho.b^4.S'' + b^2.(1 + 2\rho.b - \lambda.\rho.b).S' + \lambda.\rho.b^2.S + \lambda.b.(\rho.b^2.S' - \lambda.\rho.b.S ) \\
    & U = \rho.b^4.S'' + b^2.(1 + 2\rho.b).S' + \lambda.\rho.b^2.(1 - \lambda).S
    \end{align*}
    
    So we obtain for \ $S$ \ the differential equation :
\begin{equation*}
\rho.b^2.S'' + (1 + 2\rho.b).S' + \lambda.\rho.(1 - \lambda).S = 0 . \tag{4}
\end{equation*}

Define the sequence  \ $(\gamma_n^{\lambda})$ \  as follows :  \ $\gamma^{\lambda}_0 = 1$ \ and 
$$ \gamma_{n+1}^{\lambda} = \Big( 1 - \frac{1}{n+1} + \frac{\lambda.(1 - \lambda)}{(n+1)^2}\Big).\gamma_n^{\lambda} .$$
Then we shall have \ $S_{\rho}^{\lambda} = \sum_{n = 0}^{+\infty} s_n^{\lambda}.b^n$ \ with  \ $s_n^{\lambda} = n!(-\rho)^n.\gamma_n^{\lambda} $. This implies
$$ S_{\rho}^{\lambda}(b) = 1 + s_1^{\lambda}.b  \quad modulo \ b^2.\C[[b]] $$
with  \ $s_1^{\lambda} : = \rho.\lambda.(\lambda -1) $.\\

Remark that \ $S^{\lambda}_{\rho}$ \ is Gevrey 1 and so \ $S^{\lambda}_{\rho}$ \ acts on convergent power series germs near \ $0$.

Fix a complex number \ $\mu$. The relation \ $(1)$ \  can be written (without writing dependance in \ $\rho, \lambda, \mu$)
\begin{align*}
&(\alpha - \lambda.\beta).S_{\rho}^{\lambda} = \big(\tilde{Z}_0 + \tilde{Z}_1.(a - \mu.b)\big).(a - \lambda.b) \tag{1bis}
\end{align*}
where  \ $\tilde{Z}_i, i = 0, 1$ \ are in \  $\C[[b]]$. We simply write
$$ a.Z_1 = Z_1.a + b^2.(Z_1)'  = Z_1.(a - \mu.b) + \mu.b.Z_1 + b^2.(Z_1)'  $$
which gives
$$ \tilde{Z}_0 = Z_0 + \mu.b.Z_1 + b^2.(Z_1)' \quad {\rm and} \quad \tilde{Z}_1 = Z_1 .$$
So we have the following formulas :
\begin{align*}
& \tilde{Z}_0 = (1 - \lambda.\rho.b + \mu.\rho.b).S + 2\rho.b^2.S' = 1 + [\rho.(\mu - \lambda) + s_1^{\lambda}].b \quad modulo \ b^2.\C[[b]] \\
& \tilde{Z}_1 = \rho.S = \rho.(1 + s_1^{\lambda}.b) \quad modulo \ b^2.\C[[b]]  \tag{1ter} 
\end{align*}

\subsection{The sharp filtration of a fresco.}

\begin{defn}\label{Filt.fine}
Let \ $E$ \ be a rank \ $k$ \ fresco. We define on \ $E$ \ the following filtration :
For \ $h \in [0. k-1]$ \
 $$  \Phi_{k.n+h} : = a^h.b^n.E + b^{n+1}.E  = \sum_{j=1}^{k-h} \C.b^n.e_j + b^{n+1}.E = b^n.\Phi_h.$$
We shall call it the {\bf sharp filtration} of \ $E$.
\end{defn}

This filtration has the following properties :

\begin{enumerate}[i)]
\item For each  \ $\nu \geq 0$ \quad  $\dim_{\C} \Phi_{\nu}\big/\Phi_{\nu+1} = 1 $.
\item For each  \ $\nu \geq 0  \quad  a.\Phi_{\nu} \subset \Phi_{\nu+1}$ \ and  \quad  $b.\Phi_{\nu} \subset \Phi_{\nu+k}$.
\item For each integer \ $n \quad  \Phi_{k.n} = b^n.E$.
\item Moreover \ $a.\Phi_{k.n+ k-1} \subset \Phi_{(k+1).n +1}$ \ which implies  \ $a^k.\Phi_{\nu} \subset \Phi_{\nu+k+1}$ \ for each \ $\nu \geq 0$.
\item For an\ $\hat{A}-$linear map \ $f : E \to F$ \ where \ $E$ \ and \ $F$ \ are rank \ $k$ \ frescos \ $f$ \ sends \ $\Phi^E_{\nu}$ \ in \ $\Phi^F_{\nu}$, for each \ $\nu \geq 0$.
\item The filtration \ $(\Phi_{\nu})$ \ is invariant by any change of variable given by  \\
 $\alpha : = \theta(a) = \sum_{j=1}^{+\infty} \ \theta_j.a^j $ \ with \ $\theta_1 \not= 0 $ \ and \ $\beta : = b.\theta'(a)$:  it is easy to see that \ $\theta(a) - \theta_1.a \in a^2.\C[[a]]$ \ is of weight \ $2$ \ and that \ $b.\theta'(a) - \theta_1.b \in b.a.\C[[a]]$ \ is of weight \ $k+1$ \ for this filtration.
\end{enumerate}

Fix \ $\lambda_1 > 2$ \ a rational number and \ $p_1, p_2$ \ two integers bigger or equal  to \ $2$. We want now to compute the result of the change of variable associated to  \ $\theta_{\rho}(a) = a + \rho.a^2$ \ on the isomorphism class of the semi-simple rank \ $3$ \ fresco \ $E(\gamma)$ \ (see above). So we have to compute the function  \ $ \varphi : \C^2 \to \C$ \  such that we have \ $(\theta_{\rho})_*(E(\gamma)) \simeq E(\varphi(\rho,\gamma))$, as we know that \ $(\theta_{\rho})_*(E(\gamma))$ \ is a semi-simple rank \ $3$ \ frescos with the prescribed fundamental invariants \ $\lambda_1, p_1, p_2$, and then isomorphic to some \ $E(\gamma')$ \ for some (unique) complex number \ $\gamma'$.\\

Consider \ $E : = E(\gamma)$ \ and \ $\rho \in \C$. Let \ $e_1, e_2, e_3$ \ be the standard basis of \ $E(\gamma)$, that is to say a  \ $\C[[b]]-$basis such that
$$ (a - \lambda_3.b).e_3 = e_2 \quad (a - \lambda_2.b).e_2 = (1 + \gamma.b).e_1 \quad {\rm and} \quad (a - \lambda_1.b).e_1 = 0 .$$

\parag{Analysis} Denote \ $(F_j), j= 1,2,3$ \ the principal J-H. sequence of \ $E$ \ and   \ $\tilde{E}$ \ the fresco obtained from \ $E$ \ by the change of variable \ $\theta_{\rho}$. In fact \ $\tilde{E}$ \ is \ $E$ \ endowed with the operations \ $\alpha,\beta$ \  associated to \ $\theta_{\rho}$. There exists an unique (up to a non zero multiplicative constant) generator   \ $\varepsilon_1 \in F_1$  \ such that  \ $(\alpha - \lambda_1.\beta).\varepsilon_1 = 0$.\\
We shall show that there exists a generator \ $\varepsilon_2^0$ \ of \ $F_2$, unique modulo \ $F_1$, such that \ $(\alpha - \lambda_2.\beta).\varepsilon_2^0$ \ is in \ $F_1$ \ and induces in \ $F_1\big/b.F_1$ \ the same class than \ $\varepsilon_1$.\\
Then we shall show that, choosing \ $\varepsilon_2 = \varepsilon_2^0 + W(\beta).\varepsilon_1$, where \ $W \in \C[[\beta]]$ \ is suitable, we can choose a generator \ $\varepsilon_3$ \  of \ $\tilde{E}$ \ (as a \ $\hat{A}-$module) such that 
 $$(\alpha - \lambda_3.\beta).\varepsilon_3 = \varepsilon_2 .$$
   Then the coefficient of \ $\beta.\varepsilon_1$ \ in \ $(\alpha - \lambda_2.\beta).(\alpha - \lambda_3.\beta).\varepsilon_3$ \ will determine the isomorphism class of  \ $\tilde{E}$ \  because of the relation
\begin{equation*}
(\alpha - \lambda_2.\beta).(\alpha - \lambda_3.\beta).\varepsilon_3 = T(\beta).\varepsilon_1  \tag{0}
\end{equation*}
with \ $T(0) = 1$, and the proof above which gives  the value of \ $\gamma$ \ as \ $T'(0)$ \ in this case.

\parag{Computation of \ $\varphi(\rho,\gamma)$} Put  \ $\varepsilon_1 : = S_{\rho}^{\lambda_1}(b).e_1$. Thanks to \ $(1)$ \  with $\lambda = \lambda_1$ \  we get  \  $(\alpha - \lambda_1.\beta).\varepsilon_1 = 0 $, where \ $\alpha$ \ and \ $\beta$ \ are given by 
$$\alpha : = \theta_{\rho}(a) = a + \rho.a^2, \quad \beta : = b.\theta'(a) =  b + 2\rho.b.a. $$
Remark that  \ $b^2.E = \beta^2.E = \Phi_6$, and that the equality \ $b^2.F_1 = \beta^2.F_1$ \ implies :
\begin{equation*}
\varepsilon_1 = (1 + \rho.\lambda_1(\lambda_1-1).b).e_1  \quad modulo \  b^2.F_1 \tag{5}
\end{equation*}
and so
\begin{equation*}
e_1 =  (1 - \rho.\lambda_1(\lambda_1-1).\beta).\varepsilon_1 \quad modulo \  \beta^2.F_1 \tag{5bis}
\end{equation*}
Now we look for \ $\varepsilon_2^0$ \ which gives with \ $\varepsilon_1$ \ a basis of \ $F_2$ \ on \ $\C[[\beta]]$ \ and satisfies :
\begin{equation*}
(\alpha - \lambda_2.\beta).\varepsilon_2^0  \in F_1 \tag{6}
\end{equation*}
Put
\begin{equation*}
\varepsilon_2^0 : =  S_{\rho}^{\lambda_2}(b).e_2 \tag{7}
\end{equation*}
This implies \ $(6)$ \  thanks to \ $(1)$.\\
Then we look for a generator \ $\varepsilon_3$ \ of \ $\tilde{E}$ \ (as a \ $\hat{A}-$module) such that
\begin{equation*}
(\alpha - \lambda_3.\beta).\varepsilon_3 = \varepsilon_2^0 + W(\beta).\varepsilon_1  \tag{8}
\end{equation*}
and this is equivalent to \ $(\alpha - \lambda_3.\beta).\varepsilon_3 = \varepsilon_2$ \  where we define \ $ \varepsilon_2 : = \varepsilon_2^0 + W(\beta).\varepsilon_1$, with a suitable \ $W \in \C[[\beta]]$ \  to determine.\\
Put
\begin{equation*}
\varepsilon_3 = S_{\rho}^{\lambda_3}(b).e_3 + U(\beta).\varepsilon_2^0 + V(\beta).\varepsilon_1 \tag{9}
\end{equation*}
because in  \ $E\big/F_2$ \ the only generator \ $x$ \ (up to a non zero constant) which satisfies  \ $(\alpha - \lambda_3.\beta).x = 0$ \ is given by the image of \ $S_{\rho}^{\lambda_3}.e_3$.\\
Use the identity  \ $(1bis)$ \ with \ $\lambda = \lambda_3$ \ and \ $\mu = \lambda_2$. We obtain :
\begin{align*}
& (\alpha - \lambda_3.\beta).\varepsilon_3  = Z^{\lambda_3}.e_2 + \beta^2.U'(\beta).\varepsilon_2^0 + (\lambda_2-\lambda_3).\beta.U(\beta).\varepsilon_2^0 \  + \\
& \qquad  \quad\quad  +  U(\beta).Z^{\lambda_2}.(1 + \gamma.b).e_1 + \beta^2.V'(\beta).\varepsilon_1 + (\lambda_1 - \lambda_3).\beta.V(\beta).\varepsilon_1 \tag{10}
\end{align*}
Write
$$ Z^{\lambda_3}.e_2  =  X(\beta).\varepsilon_2^0 + Y(\beta).\varepsilon_1 .$$
The relation \ $(8)$ \  forces :
\begin{align*}
& X(\beta) + \beta^2.U'(\beta) - (p_2 - 1).\beta.U(\beta) = 1  \quad {\rm and}  \tag{11} \\
& Y(\beta).\varepsilon_1 +  U(\beta).Z^{\lambda_2}.(1 + \gamma.b).e_1 + \beta^2.V'(\beta).\varepsilon_1 + (\lambda_1 - \lambda_3).\beta.V(\beta).\varepsilon_1 = W(\beta).\varepsilon_1. 
\end{align*}
But we know "a priori" that a solution \ $U \in C[[\beta]]$ \ of the first equation \ $(11)$ \ exists, and we also know that a corresponding solution \ $V$ \  of the seconde equation \ $(11)$ \ exists, for a suitable choice of \ $U$. We want to compute \ $U(0)$ \ for such a solution, because we have :
\begin{align*}
&  (\alpha - \lambda_2.\beta).(\alpha - \lambda_3.\beta).\varepsilon_3 = (\alpha - \lambda_2.\beta).\varepsilon_2 \\
& \qquad = (\alpha - \lambda_2.\beta).\varepsilon_2^0 + W(0).(\lambda_1 - \lambda_2).\beta.\varepsilon_1 \quad modulo \  \beta^2.F_1 
\end{align*}
and then it will be enough to compute \ $W(0)$ \ and the value of 
$$  (\alpha - \lambda_2.\beta).\varepsilon_2^0 = Z^{\lambda_2}.(1 + \gamma.b).e_1 $$
in \ $F_1\big/\beta^2.F_1$. \\
The formula \ $(1bis)$ \  with \ $\lambda = \lambda_2$ \ and \ $\mu = \lambda_1$ \ gives :
\begin{align*}
&  Z^{\lambda_2}.(1 + \gamma.b).e_1  = \tilde{Z}_0^{\lambda_2}.(1 + \gamma.b).e_1 + \tilde{Z}^{\lambda_2}_1.(a - \lambda_1.b).(1 + \gamma.b).e_1 \\
& \qquad = ( 1 +  \gamma + s_2 +  \rho.(\lambda_1 - \lambda_2).b).e_1 +  b^2.F_1 .
\end{align*}
Using  \ $(1ter)$, we get
\begin{equation*}
\varphi(\rho,\gamma) =  \gamma +  s_2 + \rho.(\lambda_1 - \lambda_2) +  (\lambda_1 - \lambda_2).W(0) \tag{@}
\end{equation*}
Remark that the solution \ $U$ \ is unique up to \ $\C.\beta^{p_2-1}$ ; this shows as \ $p_2  \geq 2$, that the number \ $U(0)$ \ is independant of the choice of the solution \ $U$.\\

To determine these numbers first remark that, computing in  \ $\Phi_1\big/\Phi_2$, we have \ $X(0) = 1$ \  and, thaks to  \ $(11)$
\begin{align*}
& X'(0) = -(p_2-1).U(0) \quad {\rm  and \ also } \\ 
& W(0) = Y(0) + U(0) . \tag{12}
\end{align*}
because the constant term of \ $Z^{\lambda_2}(b)$ \  is \ $1$.\\
Using \ $(1ter)$ \ with  \ $\lambda = \lambda_3$ \ and \ $\mu = \lambda_2$ \ we get
\begin{align*}
 & Z^{\lambda_3}.e_2 = \tilde{Z}_0.e_2 + \tilde{Z}_1.(1 +\gamma.b).e_1 \\
 & \qquad = X(\beta).\varepsilon_2^0 + Y(\beta).\varepsilon_1. \tag{13}
 \end{align*}
 We have
  \begin{align*}
& \varepsilon_2^0 = (1 + s_2.b).e_2  \quad modulo \ b^2.F_2   \quad {\rm and} \\
& \varepsilon_1 = ( 1 + s_1.b).e_1 \quad modulo \ b^2.F_1 
\end{align*}
and, as  \ $(b - \beta -2\rho.\beta.\alpha).F_2 \subset b^2.F_2$, since \ $b.a^2.\Phi_1 \subset \Phi_6$,
\begin{align*}
& e_2 =  (1 - s_2.\beta).\varepsilon_2^0 -  2\rho.s_2.\beta.\alpha.\varepsilon_2^0  \quad modulo \ \beta^2.F_2 \quad {\rm recall \ that} \\
& e_1 =  ( 1 - .s_1.\beta).\varepsilon_1 \quad modulo \ \beta^2.F_1 \quad {\rm and \ so} \quad   b.e_1 = \beta.\varepsilon_1  \quad modulo \ \beta^2.F_1 .
\end{align*}
As we have 
\begin{align*}
& \alpha.\varepsilon_2^0 = (a + \rho.a^2).(1+ \rho.s_2.b).e_2 \quad modulo \ \beta^2.F_2 \\
&  \alpha.\varepsilon_2^0  =  e_1  \quad modulo \ \beta.F_2
\end{align*}
we finally get
\begin{equation*}
 e_2 =  (1 -  s_2.\beta).\varepsilon_2^0 -  2\rho.s_2.\beta.\varepsilon_1
\end{equation*}
Using again  \ $(1ter)$, with \ $\lambda = \lambda_3$ \ and \ $\mu = \lambda_2$ \ we get
\begin{align*}
&\tilde{Z}_0^{\lambda_3} = (1 + \rho.(\lambda_2 - \lambda_3).b).(1 + s_3.b) + b^2.\C[[b]] \\
&  \tilde{Z}_0^{\lambda_3} =  1 + (s_3 +  \rho.(\lambda_2 - \lambda_3)).b +  b^2.\C[[b]]  \\
& \tilde{Z}_1^{\lambda_3}  = \rho.(1 + s_3.b) 
\end{align*}
and then
\begin{align*}
 & Z^{\lambda_3}.e_2 = e_2 + ( s_3 +  \rho.(\lambda_2 - \lambda_3)).b.e_2 + \rho.e_1 + \rho.(\gamma + s_3).b.e_1  \quad modulo \ b^2.F_2  \\
&  Z^{\lambda_3}.e_2 = \varepsilon_2^0 + ( s_3 - s_2 +  \rho.(\lambda_2 - \lambda_3).\beta.\varepsilon_2^0 +   \rho.\varepsilon_1 + \rho.(\gamma + s_3 - 2s_2).\beta.\varepsilon_1 \quad modulo \ \beta^2.F_2 .
\end{align*}
So we get
\begin{align*}
& X(\beta) = 1 +  ( s_3 - s_2 +  \rho.(\lambda_2 - \lambda_3)).\beta \quad modulo \ \beta^2.\C[[\beta]] \\
& Y(\beta) = \rho + \rho.(\gamma + s_3 - 2s_2).\beta \quad modulo \ \beta^2.\C[[\beta]].
\end{align*}
Then, as \ $s_i = \rho.\lambda_i.(\lambda_i+1)$,
\begin{align*}
& X'(0) =  s_3 - s_2 +  \rho.(\lambda_2 - \lambda_3) = - (p_2-1).U(0) \\
& Y(0) = \rho \quad \quad  {\rm and \ so, \ as} \quad  s_3 - s_2 + \rho.(\lambda_2 - \lambda_3) = \rho.(\lambda_3 - \lambda_2).\big[\lambda_2 + \lambda_3 -2\big]  \\
& X'(0) = \rho.(p_2-1).(\lambda_2 + \lambda_3 - 2) \\
& U(0) = -\rho.(\lambda_2 + \lambda_3 - 2) \\
& W(0) =  \rho - \rho.(\lambda_2 + \lambda_3 - 2) = - \rho.[\lambda_2 + \lambda_3 - 3]
\end{align*}
Finally we get 
\begin{align*}
& \varphi(\rho,\gamma) = \gamma + s_2 - \rho.(\lambda_2 - \lambda_1) + (\lambda_1-\lambda_2).W(0) \\
& \varphi(\rho,\gamma) = \gamma + s_2 - \rho.(\lambda_2 - \lambda_1) +  \rho.(p_1-1).[ \lambda_2 + \lambda_3 - 3] \\
&  \varphi(\rho,\gamma) = \gamma + \rho.\big[ \lambda_2^2 - \lambda_2  + (p_1-1).[ \lambda_2 + \lambda_3 - 4] \big]\\
& \varphi(\rho,\gamma) = \gamma + \rho.\big[ \lambda_2^2 - (2p_1 - 3).\lambda_2 +(p_1-1).(p_2-5)\big].
\end{align*}

So we have proved the following proposition where the constant \ $L$ \ is equal to
$$L : = \lambda_2^2 - (2p_1 - 3).\lambda_2 +(p_1-1).(p_2-5) = \lambda_1^2 + \lambda_1 + (p_1-1).(p_1 - p_2 +3).$$

\begin{prop}\label{calcul}
Let \ $\theta(a) : = a + \rho.a^2$ \ for some \ $\rho \in \C$. Then we have \ $\theta_*(E(\gamma)) \simeq E(\gamma + \rho.L)$ \ where \ $L$ \ is a constant depending only on the fundamental invariants \ $\lambda_1, p_1,p_2$.\\
\end{prop}

Our next proposition will allow to reach the general change of variable.\\

\begin{prop}\label{la clef}
Let \ $\theta(a) = a  + \sum_{n=3}^{+\infty} \ \theta_n.a^n $ \ be a formal power serie. Then the fresco  \ $\theta_*(E(\gamma))$ \ is isomorphic to \ $E(\gamma)$.
\end{prop}

The proof will be cut in  two lemmas

\begin{lemma}\label{techn.1}
Let \ $E_{\lambda} : = \hat{A}\big/\hat{A}.(a - \lambda.b) $ \ with \ $\lambda > 0$ \ and denote \ $e_{\lambda}$ \ a generator of \ $E_{\lambda}$ \ which is annihilated by \ $a - \lambda.b$. Let \ $\theta(a) = a + a^3.\tau(a)$ \ be a change of variable and let \ $\alpha : = \theta(a)$ \ and \ $\beta : = b.\theta'(a)$. There exists an unique \ $S \in \C[[\beta]]$ \ such  \ that \ $S(0) = 1$ \ and \ $\varepsilon_{\lambda} : = S(\beta).e_{\lambda}$ \ satisfies \ $(\alpha - \lambda.\beta).\varepsilon_{\lambda} = 0$, and then \ $S'(0) = S''(0) = 0 $.
\end{lemma}

\parag{Proof} The proof of the existence of an unique \ $S \in \C[[\beta]]$ \ with \ $S(0) = 1$ \  such that \ $\varepsilon_{\lambda} : = S(\beta).e_{\lambda}$ \ satisfies \ $(\alpha - \lambda.\beta).\varepsilon_{\lambda} = 0$ \ is analoguous to the proof of the lemma \ref{rk 1 polyn. depend.}. To see that \ $S'(0) = S''(0) = 0$ \ in our case write first
$$ (\alpha - \lambda.\beta).S(\beta).e_{\lambda} = 0. $$
But we have \ $a^3.E_{\lambda} = b.a^2.E_{\lambda} = b^3.E_{\lambda} = \beta^3.E_{\lambda}$. So the equality above gives
$$ (a - \lambda.b).( 1 + S'(0).b + S''(0).\frac{b^2}{2}).e_{\lambda} \in b^4.E_{\lambda} .$$
This forces \ $S'(0) = S''(0) = 0. \hfill \blacksquare$\\

\begin{lemma} Let \ $E(\gamma)$ \ the rank \ $3$ \ (a,b)-module with basis \ $e_1,e_2,e_3$ \ with \ $a$ \ defined by the relations, where \ $\lambda_1, \lambda_2, \lambda_3$ \ are as above, 
\begin{align*}
& (a - \lambda_3.b).e_3 = e_2 \\
& (a - \lambda_2.b).e_2 = (1 + \gamma.b).e_1 \quad {\rm where} \quad \gamma \in \C \\
& (a - \lambda_1.b).e_1 = 0 .
\end{align*}
Let \ $\theta(a) = a + a^3.\tau(a)$ \ be a change of variable and let \ $\alpha : = \theta(a)$ \ and \ $\beta : = b.\theta'(a)$. Then there exists \ $\varepsilon_1, \varepsilon_2, \varepsilon_3$ \ is a \ $\C[[\beta]]-$basis of \ $\theta_*(E(\gamma))$ \ such that
\begin{align*}
& (\alpha - \lambda_3.\beta).\varepsilon_3 = \varepsilon_2
& (\alpha - \lambda_2.\beta).\varepsilon_2 = T(\beta).\varepsilon_1 \\
& (\alpha - \lambda_1.\beta).\varepsilon_1 = 0 .
\end{align*}
with \ $T(0) = 1$ \ and   \ $T'(0) = \gamma$.  So \ $\theta_*(E(\gamma))$ \ is isomorphic to \ $E(\gamma)$\footnote{see the first part of the proof of the proposition \ref{rg 3 ss 1}.}.
\end{lemma}

\parag{Proof} First remark that the principal J-H. sequence \ $F_1 \subset F_2 \subset E$ \ of \ $E : = E(\gamma)$ \ is  invariant   by the change of variable (as vector spaces). But also the quotient \ $\theta_*(F_{j+1})\big/\theta_*(F_j)$ \ is isomorphic to \ $F_{j+1}\big/F_j$ \ for each \ $j$.\\
Using  the lemma \ref{techn.1} we find an unique \ $\varepsilon_1: = S_1(\beta).e_1$ \ such that \ $(\alpha - \lambda_1.\beta).\varepsilon_1 = 0 $ \ holds with \ $S_1(0) = 1$ \ and \ $ S'_1(0) = S_1''(0) = 0 $.\\
Then using again the lemma \ref{techn.1} we find an unique \ $S_2 \in \C[[\beta]], S_2(0) = 1$ \ and \ $ S'_2(0) = S_2''(0) = 0 $ \ such if we define  \ $\varepsilon^0_2 = S_2(\beta).e_2$ \ we have 
 $$(\alpha - \lambda_2.\beta).\varepsilon_2^0 = U(\beta).\varepsilon_1 \in F_1.$$
Now we look for \ $\varepsilon_2 : = \varepsilon^0_2 + W(\beta).\varepsilon_1$ \ and 
$$ \varepsilon_3 = S_3(\beta).e_3 + X(\beta).\varepsilon_2^0 + Y(\beta).\varepsilon_1 $$
where \ $S_3 \in \C[[\beta]]$ \ such \ $S_3(0) = 1, S'_3(0) = S_3''(0) = 0 $ \ is given by the lemma \ref{techn.1} in order that we have
$$(\alpha - \lambda_3.\beta).\varepsilon_3 = \varepsilon_2.$$
We have
\begin{align*}
& (\alpha - \lambda_3.\beta).\varepsilon_3 \   = \  \varepsilon_2 \  =  \ \varepsilon_2^0 +  W(\beta).\varepsilon_1 \\
& \qquad \quad = (a - \lambda_3.b).e_3 + \xi.e_3 + X(\beta).U(\beta).\varepsilon_1 + \beta^2.X'(\beta).\varepsilon_2^0 + (\lambda_2-\lambda_3).\beta.X(\beta).\varepsilon_2^0 +\\
& \qquad\qquad +  (\lambda_1-\lambda_3).\beta.Y(\beta).\varepsilon_1 + \beta^2.Y'(\beta).\varepsilon_1 \\
\end{align*}
where \ $\xi $ \ is in \ $\hat{A}.a^3 + \hat{A}.b.a^2$. We want to show that \ $W(0) = 0$ \ so we shall consider this equation modulo \ $\Phi_5$. So we have
\begin{align*}
& a^3.e_3 = (\lambda_1 + \lambda_2 + \lambda_3).b.e_1 \quad modulo \  \Phi_6 \\
& b.a^2.e_3 \in \Phi_6 
\end{align*}
and the term \ $\xi.e_3$ \ is in \ $\Phi_5$.\\
Moreover this equality implies, looking at the coefficient of \ $\varepsilon_1$,  the relation    
$$W(0) = X(0).U(0)$$
 and as  \ $(a - \lambda_3.b).e_3 = e_2 = S_2^{-1}(\beta).\varepsilon_2^0 $ \ we obtain, looking at the coefficient of \ $\beta.\varepsilon_2^0$ \ and using  \ $S_2(0) = 1$ \ and \ $S_2'(0) = 0$
\begin{align*}
& (\lambda_2 - \lambda_3).X(0) = 0 
\end{align*}
so we conclude that \ $X(0) = 0$ \ and also \ $W(0) = 0$.\\
Now we get, as \ $\varepsilon_2^0 = e_2  \ modulo \  \Phi_6$ \ and \ $(\hat{A}.a^3 + \hat{A}.b.a^2).e_2 $ \ is contained in \ $\Phi_6$
\begin{align*}
&  (\alpha - \lambda_2.\beta).\varepsilon_2 = (\alpha - \lambda_2.\beta).\varepsilon_2^0 +  \beta^2.F_1\\
& \qquad \quad = (a - \lambda_2.b).e_2 \quad modulo \ \Phi_6 \\
&\quad\quad = (1 + \gamma.\beta).\varepsilon_1 \quad modulo \ \Phi_6 \cap F_1 \subset b^2.F_1 .
\end{align*}
This completes the proof of the second lemma and of the proposition. $\hfill \blacksquare$\\

\begin{thm}\label{changement de variable ss}
Let \ $E(\gamma)$ \ be a semi-simple fresco of the previous family and \ $\theta(a) : = \sum_{j=1}^{+\infty} \ \theta_j.a^j $ \ a change of variable. Then we have
$$ \theta_*(E(\gamma)) \simeq E( \delta) \quad {\rm with } \quad \delta : = \frac{\gamma}{\theta_1} + \frac{\theta_2}{\theta_1^2}.L   $$
where  \ $L : =  \lambda_1^2 + \lambda_1 + (p_1-1).(p_1 - p_2 +3) $ \ depends only on the fundamental invariants  \ $\lambda_1, p_1,p_2$ \ of the choosen family.
\end{thm}

\parag{Proof} We have  for \ $\theta(a) = \theta_1.a$ \ with \ $\theta_1 \in \C^*$, the easy isomorphism
$$ \theta_*(E(\gamma)) \simeq E(\gamma/\theta_1).$$
Now  a general \ $\theta$ \ is the composition of the change of variable \ $a \mapsto \theta_1.a$ \  then the change of variable  $$a \mapsto a + \frac{\theta_2}{\theta_1^2}$$
 and finally  some change of variable tangent to the identity at the order at least \ $2$. Then the propositions \ref{calcul} and \ref{la clef} allow to conclude. $\hfill \blacksquare$

\parag{Consequence} The action of the group of change of variables on the family   \ $E(\gamma), \gamma \in \C$ \ reduces to the action of the affine group on \ $\C$  \ via
$$ \theta(\gamma) = \frac{\gamma}{\theta_1} + \frac{\theta_2}{\theta_1^2}.L $$
where \ $L : = \lambda_1^2 + \lambda_1 + (p_1-1).(p_1 - p_2 +3)  $ \ depends only on the fundamental invariants  \ $\lambda_1, p_1,p_2$ \ of the choosen family.\\

 \parag{Remark} The cases where \ $L = 0$ \ are of course of interest. As we have \ $\lambda_1 > 2$ \ we shall have \ $L = 0$ \ when
 $$  \lambda_1 = \frac{\sqrt{4(p_1-1).(p_1-p_2+3) + 1} -1}{2} .$$
 This implies that \ $2\lambda_1 + 1 \in \sqrt{\mathbb{N}}\cap \mathbb{Q} = \mathbb{N}$. So put \ $\lambda_1 = r/2$ \ where \ $r$ \ is an integer. The relation
 $$ r.(r+2) = 4(p_1-1).(p_1-p_2+3) $$
 shows that \ $r$ \ is even, and so \ $\lambda_1$ \ is an integer \ $\geq 3$.\\
 Then we can find many families for which \ $L = 0$, so that the number \ $\gamma$ \ is a quasi-invariant holomorphic parameter (with the character \ $\theta \mapsto \theta'(0)^{-1}$. Here are the first examples :
 \begin{align*}
 & \lambda_1 \geq 3 \quad p_1 = \lambda_1 +1, p_2 = 3 \\
 & \lambda_1 \geq 3 \quad p_1 = \lambda_1 + 2, p_2 = 5 \\
 & \lambda_1 = 2t, t \geq 2 \quad p_1 = 4t+3, p_2 = 3t+6 \\
& \lambda_1 = 2u+1, u \geq 1 \quad p_1 =4u+3, p_2 = 3u+5.
\end{align*}
Of course, using the decomposition of \ $\lambda_1$ \ or \ $\lambda_1+1$ \ we can build many others.

\subsection{Application to holomorphic families of  \ $[\lambda]-$primitive frescos.}

Let us formulate now as a theorem the construction which is an easy consequence of our various results. \\

\begin{thm}\label{param. inv. ss}
Let \ $E$ \ be a \ $[\lambda]-$primitive semi-simple fresco of rank \ $k \geq 5$ \ and let \ $F_j, j \in [1,k]$ \ its principal Jordan-H{\"o}lder sequence.
 Assume that we may choose \ $3 \leq j_1< j_2 < j_3$, such that the integers \ $p_{j_1-2}, p_{j_1-1}, p_{j_2-2}, p_{j_2-1},p_{j_3-2},p_{j_3-1}$ \ are at least \ $2$ \  and let \ $\gamma_1, \gamma_2, \gamma_3$ \ the complex numbers defining respectively the isomorphism class of the rank \ $3$ \ semi-simple \ $[\lambda]-$primitive frescos \ $F_{j_1}\big/F_{j_1-3}, F_{j_2}\big/F_{j_2-3}, F_{j_3}\big/F_{j_3-3}$. Then the number \ $(\gamma_3- \gamma_2)\big/(\gamma_3 - \gamma_1)$ \ is an holomorphic parameter of  \  $E$ \ which is invariant by any base change.
\end{thm}

\parag{Remarks}Observe that, in the semi-simple case, if we restrict ourself to consider change of variable tangent to identity, we need only to consider two rank \ $3$ \ subquotient frescos from the principal J-H. sequence to produce an invariant holomorphic parameter.

\parag{Conclusion} Let \ $\mathbb{E}$ \ be an holomorphic family of \ $[\lambda]-$primitive frescos parametrized par a reduced complex space \ $X$. Let \ $\mathbb{F}_j, j \in [1,k]$ \ its principal J-H. sequence ; then each \ $mathbb{F}_j$ \ is an holomorphic family (thanks to the theorem \ref{}). Then for each \ $j \in [1, k-1]$ \ the function on \ $X$ \ defined as the parameter of the rank \ $2$ \ fresco \ $\mathbb{F}_{j+1}\big/\mathbb{F}_{j-1}$ \ is holomorphic and is invariant by change of variable.\\

Assume \ $X$ \ irreducible to simplify the discussion. If the generic value of \ $d(\mathbb{E}(x))$ \ is at least equal to \ $3$, then we may find a closed analytically constructible set \ $Y \subset X$ \ such that the family \ $\Sigma^1(\mathbb{E}(x)), x \in X \setminus Y$ \ is holomorphic, and we produced on \ $X \setminus Y$ \ non trival holomorphic functions using parameters of the rank \ $2$ \ subquotients of the J-H. sequence of this family of themes  (with rank \ $\geq 2$).\\

When the generic value of \ $d(\mathbb{E}(x))$ \ is arbitrary, we also used the holomorphic family (on \ $X \setminus Y$ \ of semi-simple frescos \ $S_1(\mathbb{E}$ \ to produce, using the theorem \ref{param. inv. ss} some holomorphic functions on \ $X \setminus Y$ \ which will be invariant by change of variable.\\

So there are only few cases where we cannot produce any interesting holomorphic of this king from an holomorphic family of frescos.

\newpage

\section{Bibliography.}

\begin{itemize}

\item{[Br.70]} Brieskorn, E. {\it Die Monodromie der Isolierten Singularit{\"a}ten von Hyperfl{\"a}chen}, Manuscripta Math. 2 (1970), p. 103-161.\\

\item {[B.93]} Daniel Barlet {\it Th\'eorie des (a,b)-modules I}, Univ. Ser. Math. Plenum (1993), p.1-43.

\item {[B.08]} Daniel Barlet {\it Sur certaines singularit\'es d'hypersurfaces II}, J. Alg. Geom. 17 (2008), p. 199-254.

\item {[B.09-a]} Daniel Barlet {\it P\'eriodes \'evanescentes et (a,b)-modules monog\`enes}, Bollettino U.M.I (9) II (2009), p.651-697.

\item {[B.09-b]} Daniel Barlet {\it Le th\`eme d'une p\'eriode \'evanescente}, preprint de l'Institut E. Cartan (Nancy) 2009  $n^0$  33,  57 pages.

\item{[B.11]} Daniel Barlet {\it Asymptotics of a vanishing period : the quotient themes of a given fresco} to appear as a preprint of the Institu E. Cartan  (january 2011).

 \item{[K.76]} Kashiwara, M. {\it b-function and holonomic systems}, Inv. Math. 38 (1976) p. 33-53.\\

\item{[M.74]} Malgrange, B. {\it Int\'egrale asymptotique et monodromie}, Ann. Sc. Ec. Norm. Sup. 7 (1974), p.405-430.\\

\end{itemize}

\end{document}